\RequirePackage{fix-cm}
\documentclass[nospthms,smallextended]{svjour3}
\smartqed  
\usepackage[T1]{fontenc}

\usepackage{mathptmx}
\usepackage{makeidx}
\usepackage{cancel}
\usepackage{comment}

\usepackage{amsmath}
\usepackage{amsfonts}
\usepackage{amsthm}

\usepackage{cite}
\usepackage{tikz}
\usepackage{graphicx,pst-all,bm}
\usepackage{fixltx2e}
\usepackage[scanall]{psfrag}
\usepackage{pgfplots}
\usepackage[title]{appendix}

\usepackage[normalem]{ulem}

\usepackage{subfig}
\usepackage{enumerate}
\usepackage{siunitx}

\usetikzlibrary{arrows}
\usetikzlibrary{fit}\usetikzlibrary{calc}
  \pgfdeclarelayer{background}
  \pgfsetlayers{background,main}

\newcommand{\cB}{\mathcal{B}}
\newcommand{\cC}{C}
\newcommand{\cD}{\mathcal{D}}

\newcommand{\cG}{\mathcal{G}}

\newcommand{\cL}{L}

\newcommand{\cN}{\mathcal{N}}

\newcommand{\cW}{W}

\newcommand{\cF}{\mathcal{F}}

\newcommand{\fT}{\mathbf{T}}
\newcommand{\cTT}{{\mathbb T}^{n,q}_{h}}
\newcommand{\cCC}{\cC([-h,\infty),\R^n)}
\newcommand{\cLL}{\cL_{\rm loc}^\infty(\R_{\ge 0},\R^q)}

\usetikzlibrary{arrows}

\newtheorem{proposition}{Proposition}[section]

\newtheorem{theorem}[proposition]{Theorem}
\newtheorem{corollary}[proposition]{Corollary}

\theoremstyle{definition}
\newtheorem{definition}[proposition]{Definition}
\newtheorem{remark}[proposition]{Remark}

\newtheorem{example}[proposition]{Example}

\newcommand{\N}{{\mathbb{N}}}

\newcommand{\R}{{\mathbb{R}}}
\newcommand{\C}{{\mathbb{C}}}

\newcommand{\T}{{\mathbb{T}}}

\newcommand{\dd}{\textrm d}
\newcommand{\dphi}{\overset{\textbf{.}}{\varphi}}
\newcommand{\half}{{\textstyle{\frac{1}{2}}}}
\newcommand{\quarter}{{\textstyle{\frac{1}{4}}}}

\DeclareMathOperator{\esup}{ess\,sup}
\DeclareMathOperator{\RE}{Re}

\newcommand{\ddts}{\tfrac{\text{\normalfont d}}{\text{\normalfont d}t}}

\newcommand{\setdef}[2]{\left\{\ #1\ \left|\ \vphantom{#1} #2 \right.\right\}}

\newcommand{\TTT}[3]{{\mathbb T}^{\,{#1},{#2}}_{{#3}}}

\renewcommand{\theenumi}{{\rm(\roman{enumi})}} 

{\left(\begin{smallmatrix}}
{\end{smallmatrix}\right)}

{\left[\begin{smallmatrix}}
{\end{smallmatrix}\right]}

\makeatletter
\renewcommand*\env@matrix[1][*\c@MaxMatrixCols c]{%
  \hskip -\arraycolsep
  \let\@ifnextchar\new@ifnextchar
  \array{#1}}
\makeatother

 \journalname{Mathematics of Control, Signals, and Systems}

\makeindex

\begin{document}
\title{Funnel control of nonlinear systems\thanks{This work was supported by the German Research Foundation (Deutsche Forschungsgemeinschaft) via the grant BE 6263/1-1.}}


\titlerunning{Funnel control of nonlinear systems}        

\author{\mbox{Thomas Berger  \ $\cdot$ \ Achim Ilchmann    \ $\cdot$ \   Eugene~P~Ryan}}


\institute{Corresponding author: Thomas Berger  \at Tel.: +49 5251 60-3779
            \and
             Thomas Berger \at
              Institut f\"ur Mathematik, Universit\"at Paderborn, Warburger Str.~100, 33098~Paderborn, Germany \\
              \email{thomas.berger@math.upb.de}
           \and
           Achim Ilchmann \at
            Institut f\"ur Mathematik, Technische Universit\"{a}t Ilmenau, Weimarer Stra{\ss}e 25, 98693~Ilmenau, Germany\\
              \email{achim.ilchmann@tu-ilmenau.de}           
           \and
           \mbox{Eugene P~Ryan}    \at
           Department of Mathematical Sciences, University of Bath,  Bath BA2 7AY, UK\\
              \email{masepr@bath.ac.uk}           
           }
\date{Received: date / Accepted: date}
\maketitle

\begin{abstract}
Tracking of reference signals is addressed in the context of a  class of nonlinear controlled systems modelled by $r$-th order functional differential equations, encompassing {\it inter alia} systems with unknown
``control direction'' and dead-zone input effects. A control structure is developed which ensures that, for every member of the underlying system class and every admissible reference signal, the
tracking error evolves in a prescribed funnel chosen to reflect transient and asymptotic {accuracy} objectives.
Two fundamental properties underpin the system class:  bounded-input bounded-output stable internal dynamics, and
a high-gain property (an antecedent of which is the concept of sign-definite high-frequency gain in the context of linear systems).

\keywords{nonlinear systems \and adaptive control \and asymptotic tracking \and funnel control \and relative degree \and functional
differential equations
}
 \subclass{
        93C10  
\and	93C23  
\and 	93C40  
}
\end{abstract}
\newpage
\noindent\textbf{Nomenclature}\\[3mm]
\boxed{
\begin{tabular}[ht]{lp{1pt}p{230pt}}
$\N$, $\N_0$
&&{{the set of positive, non-negative integers, respectively}}\\[2mm]
$\R_{\ge 0}$, $\mathbb{C}_{\ge 0}$&&the {sets} $[0,\infty)$, $\{\lambda\in \mathbb{C}\, | \RE (\lambda) \geq 0\}${, respectively}\\[2mm]
$\langle v, w\rangle$
   && the {Euclidean inner product of vectors} $v,w \in \R^n$ \\[2mm]
$\|x\|$
   && $\sqrt{\langle x, x\rangle}$, the Euclidean norm of $x \in \R^n$ \\[2mm]
$L^\infty(I,  \R^{n}) $
  &&  the Lebesgue space of measurable, essentially bounded {functions} $f\colon I\to\R^n$, where $I\subseteq\R$ is some interval\\[6mm]
$L^\infty_{\rm loc} (I, \R^{n})$ &&
 the set  of {measurable,} locally essentially bounded functions  $f\colon   I  \to \R^{n}$, where $I\subseteq\R$ is some interval\\[6mm]
$W^{k,\infty}(I,  \R^{n})$ &&
the Sobolev space of all functions
$f:I\to\R^n$ {with $k$-th order weak derivative $f^{(k)}$} and $f,f^{(1)},\ldots,f^{(k)}\in L^\infty(I,  \R^{n})$, where $I\subseteq\R$ is some interval and $k\in\N$\\[10mm]
 $C^k(V,  \R^{n}) $
  &&  the set of  $k$-times continuously differentiable functions  $f:    V  \to \R^{n}$, where $V\subseteq\R^m$ and $k\in\N_0$; \\
&&  $C(V,  \R^{n}) := C^0(V,  \R^{n})$
\end{tabular}
}

%
\section{Introduction}
Since its inception in 2002, the concept
 of funnel control has been widely investigated.  In its essence, the approach considers the
following basic question: for a given class of dynamical systems, with input~$u$ and output~$y$, and a given class of reference signals~$y_{\rm ref}$,
does there
exist a single control strategy (generating~$u$) which ensures that, for every member of the system class
and every admissible reference signal, the output~$y$ approaches the reference~$y_{\rm ref}$ with prescribed transient behaviour and prescribed
asymptotic accuracy?  The twofold objective of ``prescribed transient behaviour and asymptotic accuracy'' is encompassed by the adoption of a so-called ``performance funnel'' in which the error function $t\mapsto e(t):= y(t)-y_{\rm ref}(t)$ is required to evolve;  {see Fig.\,\ref{Fig:funnel}.}
{Underlying the present paper is a large class of systems described by {$r$-th order} functional differential
equations: we denote this class (which will be made precise in due course) by $\cN^{m,r}$, where $m\in\N$ denotes the dimension of both input and output.}
The information available for feedback to the controller is comprised of the instantaneous values of the output and its first~$r-1$ derivatives, together with the instantaneous values of the reference signal and its
first~$\hat r-1$ derivatives, where $1\leq \hat r \leq r$.   A feedback strategy is developed which assures attainment of the above twofold performance objective{{: this is the core of the main result, Theorem~\ref{Thm:FunCon-Nonl}.  We proceed to highlight the features and
distinguishing novelties of this result {\em vis \`a vis} the existing literature.}}

\subsection{\bf{Novelties and literature}}\label{Ssec:Novelties}
%

\noindent
\textit{Predecessors and relative degree}: \
The  parameter~$r$  coincides with the concept of relative degree
for  many nonlinear  examples belonging to the class { $\cN^{m,r}$.  The class is, however, of sufficient generality to encompass, not only such
examples, but also systems which do not necessarily  have  a  relative  degree as defined in, for example,~\cite{Isid95}. }Adaptive control for systems with relative degree~$r>1$ has
been an issue since the early days of high-gain adaptive control,
{ as evidenced by the} contribution~\cite{Mare84} from~1984. An early approach which takes transient behaviour into account is~\cite{MillDavi91}
in 1991, using a feedback strategy that differs { conceptually}  from the funnel
control methodology. Funnel control was introduced in~2002 by~\cite{IlchRyan02b} for
nonlinear functional systems of the form~\eqref{eq:nonlSys} with relative degree one, using a variant of the high-gain property from
Definition~\ref{Def:high-gain}. The efficacy of funnel control for systems~\eqref{eq:nonlSys} with arbitrary~$r\in\N$ was demonstrated
in~\cite{IlchRyan07} in~2007.
However, the control structure in that paper is based  on backstepping with attendant (but unavoidable)
escalating controller complexity {\it vis \`a vis} the striking simplicity of the funnel controller for relative-degree-one systems. An alternative controller was developed in~\cite{LibeTren13b} for a special class of systems with~$m=1$ and arbitrary~$r\in\N$, { termed} the bang-bang funnel controller. Since the control input switches only between two values, it is able to respect input constraints; however, it requires various feasibility assumptions and involves a complicated switching logic. A simpler control strategy for nonlinear system has been
introduced by~\cite{HackHopf13} for~$r=2$ in~2013 and by~\cite{BergLe18} for~$r\in\N$ in~2018.
~\\[-1.5ex]

\noindent
\textit{Controller complexity}: \
{ Already alluded to in the above paragraph, some explicit remarks on the issue of controller complexity
may be warranted.
 For implementation purposes, the avoidance of excessive complexity is crucial. The first approaches to funnel control for systems with arbitrary relative degree in~\cite{IlchRyan06b,IlchRyan07} showed a significant increase in controller complexity with increasing relative degree
 (a variant of the ``curse of dimensionality").}
Although these contributions have the advantage that  only
the output~--   and  not its derivatives~--  need to be  known,
{ they involve an} intrinsic backstepping procedure { which requires increasing powers of a particular gain function as the relative degree
grows.   For ``large'' relative degree, this leads to impracticality.}
Avoiding the backstepping procedure,  a low-complexity funnel controller has been developed in~\cite{HackHopf13} for relative degree two systems and in~\cite{BergLe18} for arbitrary relative degree. Nevertheless,
the control design developed in~\cite{BergLe18} involves successive derivatives of { particular} auxiliary error variables,
{ causing high-level complexity for high}  relative degree. The { relative simplicity of funnel control design
underpinning Theorem~\ref{Thm:FunCon-Nonl} helps to resolve these complexity issues.}
~\\[-1.5ex]

\noindent
\textit{Unknown control direction}: \
In the early days of high-gain adaptive control without system
identification,
 linear systems with relative degree one and positive
high-frequency gain were studied, cf.\ Section~\ref{Ssec:LinSys}.
In 1983,
\textit{Morse}~\cite{Mors83} conjectured the non-existence of a smooth adaptive controller which stabilizes every { such} system under the
{ weakened} assumption that   the high-frequency gain is not zero but its sign is unknown.
\textit{Nussbaum}~\cite{Nuss83}
showed { (constructively)} that Morse's conjecture is false.
He introduced  a { class of} sign-sensing or probing ``switching functions'' in the feedback { design}, see Section~\ref{Ssec:contr_dir}.
In the present work, we allow for a larger class of switching functions { (namely, continuous surjective maps~$[0,\infty)\to\R$, which properly contain the ``Nussbaum'' class), potentially} advantageous in applications.
~\\[-1.5ex]

\noindent
\textit{Dead-zone input}: \
A dead-zone input is a special case of input nonlinearity where the value of the nonlinearity is zero when the input is between some prescribed deadband parameters, see Sections~\ref{Ssec:InputNonl} and~\ref{Ssec:DeadZone}. A dead-zone input may appear in practical applications such as hydraulic servo valves and electronic motors, and it may severely affect the performance of a control system, see e.g.~\cite{TaoKoko96, TaoLewi01}. Several approaches have been undertaken to treat these problems, see~\cite{Na13, TaoKoko96, TaoLewi01} and the references therein. We show that the system class~$\cN^{m,r}$ encompasses a larger class of dead-zone inputs than previously considered in the literature.
~\\[-1.5ex]

\noindent
{\textit{Practical and {exact} asymptotic tracking}: \
The ``performance funnel'', which we denote by
\[
 \mathcal{F}_{\varphi}
    := \setdef{ (t,e) \in \R_{\ge 0}\times\R^m }{ \varphi(t) \, \|e\| < 1 },
\]
in which the tracking error is required to evolve, is determined by the choice of a continuous function $\varphi\colon \R_{\ge 0}\to\R_{\ge 0}$ with
requisite properties which include positivity on~$(0,\infty)$ and boundedness away from zero ``at infinity":
\[
\forall\,t >0:\ \varphi(t) >0\quad\text{and}\quad \liminf_{t\to\infty }\varphi(t) >0.
\]
For example, the unbounded function
$\varphi\colon t\mapsto e^{\alpha t}-1$, $\alpha >0$,
 is an admissible choice, in which case evolution in~$\mathcal{F}_{\varphi}$
ensures that the tracking error~$e(\cdot )$ converges to zero exponentially fast.  In particular, {{\em exact} asymptotic} tracking is achieved.  However, there
is a price to pay.  A fundamental ingredient of the funnel controller is the quantity~$\varphi (t){e}(t)$
which, in the case of unbounded~$\varphi$, inevitably leads to an ill-conditioned computation of the product of ``infinitely large'' and ``infinitesimally small''
terms.  Therefore, whilst of theoretical interest, the case of unbounded~$\varphi$ may be of limited utility in applications. }
If~$\varphi$ is bounded, then the { radius of the funnel~$t$-section $\mathcal{F}_{\varphi}\cap\big(\{t\}\times\R^m\big)$}
is uniformly bounded away from zero and so asymptotic tracking is not achieved. However, the { choice of (bounded)~$\varphi$}
is at the designer's discretion and { so {\em practical} tracking with arbitrarily small (but non-zero) prescribed asymptotic accuracy is achievable
without encountering the ill-conditioning present in the exact asymptotic tracking case.}\newline
{The assumption of {\em bounded}~$\varphi$ is widespread in {{the literature on}} funnel control. Exact asymptotic tracking
with {\em unbounded}~$\varphi$ {{was achieved}}}  in~\cite{RyanSang09} for a class of nonlinear relative degree one systems: in~\cite{IlchRyan06a} a predecessor
for linear relative degree one systems was developed utilizing the internal model principle.
Recently (and unaware of the latter results) it was observed in~\cite{LeeTren19} that
asymptotic funnel control is possible for a class of nonlinear single-input single-output systems, {{albeit more restrictive}} than the class~$\cN^{m,r}$ of the present paper.
Note also that asymptotic tracking via funnel control for systems with relative degree two has been shown by~\cite{VergDima19,VergDima20}. However, the {{radius of the}} funnel in these works {{is}} bounded away from zero and the property of exact asymptotic tracking is achieved
{{at the expense of}} a {\em discontinuous} control scheme.
\\[1ex]
\noindent
{{\textit{Parameter~${\hat r}\le r$}:}}\  {{Throughout, it is assumed that the instantaneous values of the output $y(t)$ and its first $r-1$
derivatives, together with the instantaneous value $y_{\rm ref}(t)$ of the reference signal, are available for feedback purposes.  However,
in applications, some derivatives of the reference signal may not be accessible by the controller. The parameter $\hat r\in\{1,\ldots,r\}$ quantifies the
number of derivatives that are available, and so the instantaneous information signal fed to the controller is encapsulated by the
vector $\mathbf e(t)=\big(e(t),\ldots,e^{(\hat r-1)},y^{(\hat r )} (t),\ldots, y^{(r-1)}(t)\big)$ with $e(t)=y(t)-y_{\rm ref}(t)$.
The potential to cope with non-availability of reference signal derivatives
 might be advantageous for applications.  Of course, the larger the value of $\hat r$, the more
 information is available for control use, and so it might reasonably be expected that controller ``behaviour'' improves with
 increasing $\hat r$.  This expectation is borne out by numerical simulations.}}
\\[1ex]
\noindent
\black
\textit{Prescribed Performance Control}: \
A relative of funnel control is the approach of \textit{prescribed performance control} developed by \textit{Bechlioulis and Rovithakis}~\cite{BechRovi08} in~2008.  Using so-called {performance functions} {(which admit a funnel interpretation)} and a transformation that incorporates these functions,
the original controlled system is expressed in a form for which boundedness of the states, via the prescribed performance control input, can be proved -- achieving evolution of the tracking error within the funnel defined by the performance functions. The controller presented in~\cite{BechRovi08} is not of high-gain type. Instead, neural networks are used to approximate the unknown nonlinearities of the system, {{resulting}} in a  complicated controller structure. After some developments, the complexity issue has been addressed in~\cite{BechRovi14} in~2014, where prescribed performance control is shown to be feasible for {{systems in}} pure feedback for.  {{However, the $t$-sections of the funnels corresponding to the underlying
performance functions
have radii bounded away from zero and so
exact asymptotic tracking cannot be achieved, see e.g.~\cite{BechRovi14}. Whilst funnel control and prescribed performance control are motivated by similar design objectives, the
solution methodologies are intrinsically different.}}
\\[1ex]
\noindent
\textit{Applications}: \
The new funnel control strategy has a potential impact on various applications.
Since its development in~\cite{IlchRyan02b} the funnel controller proved an appropriate tool for tracking problems in various applications such as temperature control of chemical reactor models~\cite{IlchTren04}, control of industrial servo-systems~\cite{Hack17} and underactuated multibody systems~\cite{BergOtto19}, speed control of wind turbine systems~\cite{Hack14,Hack15b}, current control for synchronous machines~\cite{Hack15a}, DC-link power flow control~\cite{SenfPaug14}, voltage and current control of electrical circuits~\cite{BergReis14a}, oxygenation control during artificial ventilation therapy~\cite{PompAlfo14}, control of peak inspiratory pressure~\cite{PompWeye15} and adaptive cruise control~\cite{BergRaue18,BergRaue20}.

\subsection{\bf System class}\label{Ssec:Intr:SysClass}

We { make precise the underlying} class~$ \cN^{m,r}$ of systems, modelled by nonlinear functional differential equations
of the form
\begin{equation}\label{eq:nonlSys}
\begin{aligned}
y^{(r)}(t)&= f\big(d(t), \fT(y,\dot{y},\dots,y^{(r-1)})(t), u(t)\big) \\
y|_{[-h,0]}&= y^0\in C^{r-1}([-h,0], \R^m),
\end{aligned}
\end{equation}
where~$h\ge 0$ quantifies the ``memory'' in the system,
$r\in\N$ is related to the concept of relative degree,
$m\in\N$ is the dimension of both the input~$u(t)$ and output~$y(t)$ at time~$t\geq 0$,
$d\in\cL^\infty(\R_{\ge 0},\R^p)$ is a ``disturbance'',
and $f\in\cC (\R^p\times\R^q \times\R^m, \R^m)$
belongs to a set of nonlinear functions characterized by a particular a high-gain property (made precise in Definition~\ref{Def:high-gain}).
The operator $\fT\colon \cCC\to \cLL$, where~$n=rm$,
belongs to the class~$\cTT$ of mappings which are causal, satisfy  a local Lipschitz  condition,  and map bounded functions to bounded
functions (made precise in Definition~\ref{Def:Operator_class}).
The most simple,  but non-trivial, prototype of the system class~$\cN^{m,r}$ are linear systems with strict relative degree~$r$ and  asymptotically stable zero  dynamics (see
Section~\ref{Ssec:LinSys}).
\begin{definition}[\textbf{Operator class}]\label{Def:Operator_class}
For $n,q\in\N$ and $h \geq 0$, the set~$\cTT$ denotes the class of operators
\[
\cTT := \setdef{\fT\colon \cCC\to \cLL }{ \text{(TP1)~--~(TP3) hold}},
\]
where~(TP1)~--~(TP3) denote the following properties.
\begin{enumerate}[\hspace{2pt}\textbf{(TP1)}]
\item[\textbf{(TP1)}]
 {\it Causality:} for all~$\zeta$, $\theta \in \cCC$ and all~$t\ge 0$,
\[
\zeta|_{[-h,t]} =\theta|_{[-h,t]} ~~\implies~~ \fT(\zeta)|_{[0,t]}=\fT(\theta)|_{[0,t]}.
\]
\item[\textbf{(TP2)}]
 {\it Local Lipschitz property:} for each~$t\ge 0$ and all $\xi\in\cC([-h,t],\R^{n})$, there exist positive constants $c_0, \delta, \tau >0$ such that, for all $\zeta_1,\zeta_2 \in \cCC$ with $\zeta_i|_{[-h,t]} = \xi$
and $\|\zeta_i(s)-\xi(t)\|<\delta$ for all $s\in[t,t+\tau]$ and $i=1,2$, we have
\[
\esup\nolimits_{s\in [t,t+\tau]}\|\fT(\zeta_1 )(s)-\fT(\zeta_2) (s)\| \leq c_0 \sup\nolimits_{s\in [t,t+\tau]}\|\zeta_1(s)-\zeta_2(s)\|.
\]
\item[\textbf{(TP3)}]
 {\it Bounded-input bounded-output (BIBO) property:}  for each $c_1 >0$  there exists $c_2 >0$ such that, for all $\zeta\in\cCC$,
\[
\sup\nolimits_{t\in[-h,\infty)}\|\zeta(t)\|< c_1 ~~\implies~~ \esup\nolimits_{t\ge 0}\|\fT(\zeta)(t)\| < c_2.
\]
\end{enumerate}
\end{definition}

Property~(TP1) is entirely natural in the context of physically-motivated controlled systems.  Property~(TP2) is
a technical condition which (in conjunction with continuity of~$f$) plays a role in ensuring well-posedness of the
initial-value problem~\eqref{eq:nonlSys} under feedback control.  Property~(TP3) is, loosely speaking, a stability condition on
the ``internal dynamics'' {of~\eqref{eq:nonlSys}}.
 For linear systems with strict relative degree, the first two conditions are trivially satisfied,
whilst  the third is equivalent to a minimum-phase assumption:
this is shown in Section~\ref{Ssec:min-ph}.

The formulation also embraces nonlinear delay elements and hysteretic effects, as we shall briefly illustrate.
\\[1ex]
{\em Nonlinear delay elements.}
For $i=0,\ldots,k$, let $\Psi_i\colon\R\times \R^m\to \R^q$ be measurable in its
first argument and locally
Lipschitz in its second argument, uniformly with respect to its first argument.  Precisely, 
 for each $\xi\in\R^m$, $\Psi_i(\cdot,\xi)$ is
measurable, and
for every compact $C\subset\R^m$, there exists a constant~$c>0$ such that
\[
\text{for a.a. $t\in \R$ \ $\forall\, \xi_1,\xi_2\in C$:} \
\|\Psi_i(t,\xi_1)-\Psi_i(t,\xi_2)\|\leq c\|\xi_1-\xi_2\|.
\]
Let $h_i >0$, $i=0,\ldots,k$, and set $h:= \max_i h_i$.
For $y\in C([-h,\infty),\R^m)$, let
\[
\fT(y)(t):=
\int_{-{h_0}}^0 \Psi_0(s,y(t+s))\, \dd s
+  \sum_{i=1}^k \Psi_i (t,y(t-h_i)), \  t\ge 0.
\]
The operator~$\fT$, so defined (which models distributed and point delays),
is of class~${{\mathbb T}^{m,q}_{h}}$; for details, see~\cite{RyanSang01}.
\\[1ex]
{\em Hysteresis.}   A large class of nonlinear operators
$\fT:C(\R_{\ge 0},\R)\to C(\R_{\ge 0},\R)$,
which includes many physically-motivated
hysteretic effects, is defined in~\cite{LogeMawb00}.  These operators are contained in the class~${{\mathbb T}^{1,1}_{0}}$.
Specific examples include relay hysteresis, backlash hysteresis, elastic-plastic hysteresis, and Preisach operators.
For further details, see~\cite{IlchRyan02a}.
\\[1ex]
Next, we introduce a high-gain property which, in effect, characterizes the class of admissible nonlinearities~$f$.

\begin{definition}[\bf High-gain property] \label{Def:high-gain}
For $p,q,m\in\N$, a function $f\in\cC (\R^p\times\R^q \times\R^m, \R^m)$
is said to have the {\it high-gain property}, if there exists~$v^*\in(0,1)$ such that,
for  every compact~$K_p\subset \R^p$ and compact~$K_q\subset\R^q$,  the (continuous) function
\[
\chi\colon\R\to\R, \
s \mapsto \min\setdef{\langle v, f(\delta,z,-s v)\rangle }{
(\delta,z)\in K_p\times K_q, 
~v\in\R^m,~v^* \leq \|v\| \leq 1
 }
\]
is such that $\sup_{s\in\R} \chi(s)=\infty$.
\end{definition}
We elucidate the  high-gain property~--
which at first sight might seem somewhat arcane~--
{in the following  two remarks,
the first of which treats the linear case.}

\begin{remark}\label{hgp-rem-lin}
Why the terminology ``high-gain property''
and how is it related to ``high-gain stabilization''? Consider   a very specific class of linear systems with no disturbance~$d(\cdot)$:

\begin{equation}\label{eq:lin}
\dot y(t)= L_1 y(t) + L_2 u(t),
\qquad  \text{for \ $L_1, L_2 \in\R^{m\times m}$.}
\end{equation}
For this system class the following implications hold.
\[
\begin{array}{c}
\text{
System \eqref{eq:lin}    has the high-gain property.}
\\[1mm]
\big\Updownarrow\\[1mm]
\text{$L_2\in\R^{m\times m}$ is sign definite,}
\\
\text{i.e.,  there exists $\sigma\in\{-1,1\}$ such that
$\sigma\langle v,L_2v\rangle >0$ for all $v\in\R^m\setminus \{0\}$.}\\[1mm]
\big\Downarrow\\[1mm]
\text{System~\eqref{eq:lin}  is \textit{high-gain stabilizable},}
\\
\text{i.e., there exist $\sigma\in\{-1,1\}$ and $k^*>0$
such that,  for all $k\ge k^*$, the control}
\\
\text{$u(t)=-\sigma\,k\, y(t)$ renders system~\eqref{eq:lin}  exponentially stable.}
\end{array}
\]
%
The equivalence of the first two statements is shown in Section~\ref{Ssec:HFG}. If
 $L_2\in\R^{m\times m}$ is sign definite, then there exists $\sigma\in\{-1,1\}$ such that
$\sigma L_2$ is positive definite and so
\[
\exists \,  k^*>0\ \forall \, k\ge k^*\ \forall\, \lambda\in \C: \ \
\det\big(\lambda I_m - k^{-1} L_1 + \sigma L_2\big) = 0\ \implies\
\text{Re}\, \lambda < 0,
\]
whence the conclusion that, for a sufficiently high value of the \emph{gain
parameter}~$k >0$, the linear control~$u(t)=-\sigma k y(t)$ renders the
system exponentially stable.
The reverse implication does not hold. As a counterexample,
consider~\eqref{eq:lin} with
\[
L_1= \left[\begin{matrix}   0&-1\\ 1&  0 \end{matrix}\right] ~~\text{and}~~L_2 =
\left[\begin{matrix}   1&0\\ 0&  0 \end{matrix}\right],
\]
which, under the control~$u(t)=- k y(t)$, takes the form
\[
    \dot y(t) =  (L_1-kL_2) y(t) = \begin{bmatrix} - k & ~-1\\ ~~1&~~~0\end{bmatrix} y(t),
\]
which is exponentially stable for all $k >0$ and so
\eqref{eq:lin}  is high-gain stabilizable. However,~$L_2$ is not sign definite.
\end{remark}

 The high-gain property in Definition~\ref{Def:high-gain} extrapolates the
 above  observations to a nonlinear setting.

\color{black}

\begin{remark}\label{hgp-rem}
\begin{enumerate}[(a)]
\item\label{item:13b}
The high-gain property holds { for $f \in\cC (\R^p\times\R^q \times\R^m, \R^m)$ if, and only if, there exists
$v^*\in(0,1)$ such that,
for every compact $K_p\subset \R^p$ and compact $K_q\subset\R^q$,
at least one of the following
two  properties is  true for the continuous function~$\chi$ defined as in Definition \ref{Def:high-gain}:}
\begin{equation}\label{eq:chi_+- }
\text{(i)}    \quad  \sup_{s >0}\chi(s) =\infty
\qquad \text{or}\qquad
\text{(ii)}      \quad \sup_{s <0} \chi(s)  =\infty .
\end{equation}
{ If (i) (respectively, (ii)) holds for every such pair $(K_p,K_q)$, then we say that $f$ has the {\it negative-definite high-gain property}
(respectively, {\it the positive-definite high-gain property}).}

\item
{ That a function $f$ may have both the negative-definite and positive-definite high-gain properties is illustrated by following example.
Let $m=1$ and let $f$ (with trivial $(\delta,z)$ dependence) be given by}
\[
f(\delta,z,u) = u \ \sin\big(\ln(1+|u|)\big),
\quad (\delta,z,u)\in\R^p\times\R^q\times\R,
\]
which has the set of zeros $\{u_k, -u_k\}$ with
\[
    u_k = e^{k\pi} - 1,\quad k\in\N_0.
\]
Define the sequence $(s_k)$ by
\[
    s_k := \half (u_{k+1} - u_k) =
    \half e^{k\pi} (e^\pi - 1) > 0 ,  \quad k\in\N.
\]
Noting that $4e^{{\pi}/{2}} < e^\pi -1$, we have
\[
\ln \big(1+\half s_k \big)=\ln\big(e^{k\pi}\big(e^{-k\pi}+\quarter(e^\pi -1)\big)\big) > k\pi+\tfrac{\pi}{2}.
\]
Also,
\[
\ln \big(1+ s_k \big)=\ln\big(e^{k\pi}\big(e^{-k\pi}+\half(e^\pi -1)\big)\big)
< \ln\big(e^{k\pi}( e^\pi/2)\big)=(k+1)\pi-\ln 2.
\]
Therefore, for all $v\in \R$ with $\half \le |v| \le 1$ we have 
\[
k\pi + \tfrac{\pi}{2} < \ln(1+s_k|v|) < (k+1)\pi -\ln 2.
\]
It follows that
\[
0 < \sin (\pi-\ln 2) < \left\{\begin{array}{ll}+\sin \big(\ln (1+s_{k}|v|)\big),& \text{$k$ even}\\-\sin\big(\ln (1+s_{k}|v|)\big),
&\text{$k$ odd}\end{array}\right\} < 1.
\]
Set ${v^*} =\half$.  Then we find that
\[
 \chi(-s_{2k})
    ={ \min_{\frac{1}{2}\le |v| \le 1}  s_{2k} }\ v^2 \ \sin\big(\ln(1+s_{2k}|v|)\big) \\
    >  \tfrac14 \, s_{2k} \sin (\pi-\ln 2)
\]
and
\[
 \chi(s_{2k+1})
    = \min_{\frac{1}{2}\le |v| \le 1}  -s_{2k+1} \ v^2 \ \sin\big(\ln(1+s_{2k+1}|v|)\big) \\
    >  \tfrac14 \, s_{2k+1}\sin (\pi-\ln 2).
\]
Since $\sin (\pi-\ln 2) >0$, it follows that $\sup_{s >0}\chi (s)=\infty =\sup_{s <0}\chi(s)$.  Therefore,
{ $f$ has both the negative-definite and positive-definite high-gain properties.}
%

\item
For linear systems with strict relative degree,
we will show in Section~\ref{Ssec:HFG}
that~(i) (respectively,~(ii))
 is equivalent to
the high-frequency gain matrix being negative definite (respectively, positive definite).

\item
If it is known in advance
that the  negative-definite (respectively, positive-definite)
high-gain property holds, then the controller structure
can be simplified considerably as we will discuss in Remark~\ref{hgp-rem-rev}.
\end{enumerate}
\end{remark}

Now we are in a position to define the general system class to be considered.

\begin{definition}[\textbf{System class}]\label{Def:System_class}
For $m, r\in\N$ we say that system~\eqref{eq:nonlSys} belongs to the
\textit{system class}~$\cN^{m,r}$, written  {$(d,f,\fT)\in\cN^{m,r}$},
if, for some $p,q\in\N$ and $h\ge 0$ {the following hold}: $d\in\cL^\infty(\R_{\ge 0},\R^p)$,  $f\in\cC (\R^p\times\R^q \times\R^m, \R^m)$
has the high-gain property, and the operator $\fT$ is of class {$\T^{rm,q}_h$}.
\end{definition}
We emphasize that the system class $\cN^{m,r}$ is parameterized only by two integers, namely, $m$ (which denotes the
common dimension of the input and output spaces) and~$r$ (which is related to the concept of relative degree).
In particular, the class~$\cN^{m,r}$ encompasses systems with arbitrary state space dimension, including systems with infinite-dimensional
internal dynamics, see e.g.~\cite{BergPuch20a}:  we will elaborate further on this in Section~\ref{Sec:Concl}.

\subsection{\bf Control objectives}\label{Ssec:ContrObj}
\begin{figure}[h!t]
\begin{center}
\hspace{4cm}
\includegraphics[trim=0cm 3cm 0cm 4cm,clip=true,width=0.9\textwidth]{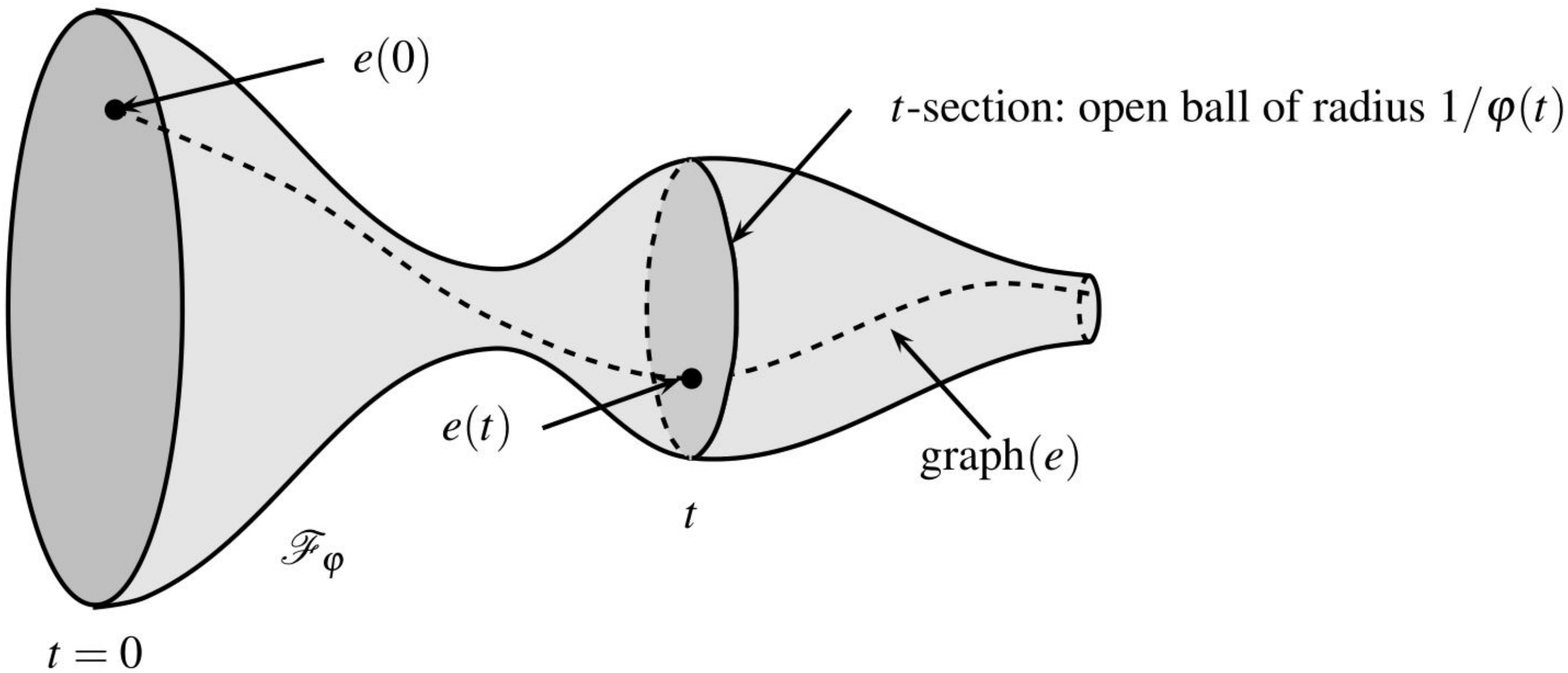}
\vspace{-1cm}
\end{center}
\caption{Performance funnel $\cF_\varphi$.} \label{Fig:funnel}
\end{figure}

%
The control problem to be addressed is to determine an output derivative feedback strategy which ensures that, for every system of class~\eqref{eq:nonlSys} and any reference signal $y_{\rm ref}\in\cW^{r,\infty}(\R_{\ge 0},\R^m)$, the output~$y$ approaches the reference~$y_{\rm ref}$ with prescribed transient behaviour and asymptotic accuracy.
This objective is reflected in the adoption of a so-called ``performance funnel'', defined by
\begin{equation} \label{eq:funnel}
   \mathcal{F}_{\varphi}
    := \setdef{ (t,e) \in \R_{\ge 0}\times\R^m }{ \varphi(t) \, \|e\| < 1 },
\end{equation}
in which the error function $t \mapsto e(t):=y(t)-y_{\text{\rm ref}}(t)$ is required to evolve;
see Fig.\,\ref{Fig:funnel}.

The funnel is shaped~-- through the choice of its boundary (determined by the reciprocal of $\varphi$)~-- in accordance with the specified transient behaviour and asymptotic accuracy;
$\varphi$  is assumed to belong to the set
\begin{equation*}\label{Phi-r}
\Phi:=\setdef{\hspace{-0.8ex}  \varphi\in {\mathrm{AC}}_{\mathrm{loc}}(\R_{\ge 0},\R_{\ge 0})}
{
\begin{array}{l}
\forall\,t >0:\ \varphi (t) >0~,
\quad \liminf_{t\to\infty}\varphi(t) >0,\\[1ex]
\exists\, c>0 :~|\dphi(t)|\leq c\big(1+\varphi(t)\big)~\text{for a.a.}\ t\geq 0
\end{array}
},
\end{equation*}
where ${\mathrm{AC}}_{\mathrm{loc}}(\R_{\ge 0},\R_{\ge 0})$ denotes the set of locally absolutely continuous functions $f: \R_{\ge 0}\to\R_{\ge 0}$.
Note that, for $t >0$, the funnel $t$-section $\cF_\varphi \cap \big(\{t\}\times\R^m\big)$ is the open ball in $\R^m$ of radius $1/\varphi (t)$.

While it is often convenient to adopt a monotonically
shrinking funnel (through the choice of a monotonically increasing function $\varphi$), it might be advantageous to
widen the funnel over some later time intervals to accommodate, for instance, periodic disturbances or strongly varying
reference signals.

\subsection{\bf Funnel control structure}\label{Sec:ContrStruc}
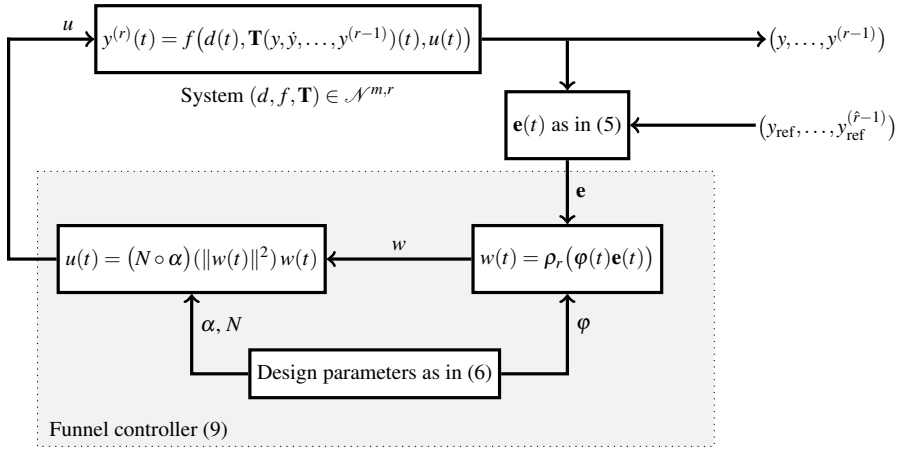
\begin{figure*}[h!t]
\centering
\resizebox{\textwidth}{!}{
   \begin{tikzpicture}[very thick,scale=0.7,node distance = 9ex, box/.style={fill=white,rectangle, draw=black}, blackdot/.style={inner sep = 0, minimum size=3pt,shape=circle,fill,draw=black},blackdotsmall/.style={inner sep = 0, minimum size=0.1pt,shape=circle,fill,draw=black},plus/.style={fill=white,circle,inner sep = 0,very thick,draw},metabox/.style={inner sep = 3ex,rectangle,draw,dotted,fill=gray!20!white}]
 \begin{scope}[scale=0.5]
    \node (sys) [box,minimum size=7ex]  {$y^{(r)}(t)= f\big(d(t), \fT(y,\dot{y},\dots,y^{(r-1)})(t), u(t)\big)$};
    \node [minimum size=0pt, inner sep = 0pt,  below of = sys, yshift=3ex] {System $(d,f,\fT)\in  \cN^{m,r}$};
    \node(fork1) [minimum size=0pt, inner sep = 0pt,  right of = sys, xshift=20ex] {};
    \node(end1)  [minimum size=0pt, inner sep = 0pt,  right of = fork1, xshift=18ex] {$\big(y,\ldots,y^{(r-1)}\big)$};

   \draw[->] (sys) -- (end1) node[pos=0.4,above] {};

  \node(FC0) [box, below of = fork1,yshift=0ex,minimum size=7ex] {{${\mathbf e}(t)$ as in~\eqref{eq:fback-quantities}}};
    \node(FC1) [box, below of = FC0,yshift=-5ex,minimum size=7ex] {$w(t)=\rho_r\big(\varphi(t){\mathbf e}(t)\big)$};
    \node(phi) [minimum size=0pt, inner sep = 0pt,  right of = FC0, xshift=18ex] {$\big(y_{\rm ref},\ldots,y_{\rm ref}^{(\hat r-1)}\big)$};
    \draw[->] (fork1) -- (FC0) {};
    \draw[->] (phi) -- (FC0) node[midway,above] {};
    \draw[->] (FC0) -- (FC1) node[midway,right] {${\mathbf e}$};
    \node(FC2) [box, left of = FC1,xshift=-30ex,minimum size=7ex] {$u(t) =\big(N\circ \alpha\big)(\|w(t)\|^2) \, w(t)$};
   \node(alpN) [minimum size=0pt, inner sep = 0pt,  below of = FC2, yshift=-5ex] {};
   \draw[->] (FC1) -- (FC2) node[midway,above] {$w$};
   \node (DP)[box, right of = alpN,xshift=10ex,yshift=2ex,minimum size=5ex]{Design parameters as in~\eqref{eq:fcts-FC}};
   \draw[->] (DP) -| (FC2) node[pos=0.8,right] {$\alpha$, $N$};
   \draw[->] (DP) -| (FC1) node[pos=0.8,right] {$\varphi$};
   \node(fork2) [minimum size=0pt, inner sep = 0pt,  left of = sys, xshift=-20ex] {};
   \draw (FC2) -| (fork2.north) {};
   \draw[->] (fork2.west) -- (sys) node[pos=0.7,above] {$u$};

   \node [minimum size=0pt, inner sep = 0pt,  below of = alpN, yshift=5ex, xshift=-5.5ex] {Funnel controller~\eqref{eq:FC}};
\end{scope}
\begin{pgfonlayer}{background}
      \fill[lightgray!20] (-4.7,-2.5) rectangle (8,-7.7);
      \draw[dotted] (-4.7,-2.5) -- (8,-2.5) -- (8,-7.7) -- (-4.7,-7.7) -- (-4.7,-2.5);
  \end{pgfonlayer}
  \end{tikzpicture}
}
\caption{Construction of the funnel controller~\eqref{eq:FC} depending on its design parameters.\\[-4mm]}
\label{Fig:funnel-controller}
\end{figure*}

%
We outline the  design  of  funnel control for  any system~\eqref{eq:nonlSys} of class $\cN^{m,r}$.\\[1ex]
{\it Information available for feedback.}\ \
Throughout, it is assumed that the instantaneous value of the output $y(t)$ and its first $r-1$ derivatives $\dot y(t),\ldots,y^{(r-1)}(t)$
are available for feedback.
Admissible reference signals are functions $y_{\textrm{ref}} \in W^{r,\infty}(\R_{\ge 0},\R^m)$.
The instantaneous reference value $y_{\textrm{ref}}(t)$ is assumed to be accessible to the controller
and, if $r\geq 2$, then, for some $\hat r\in\{1,\ldots,r\}$, the derivatives $\dot y_{\textrm{ref}}(t),\ldots,y^{(\hat r-1)}_{\textrm{ref}}(t)$ are also
accessible for feedback.  In summary,
 for some $\hat r\in \{1,\ldots,r\}$, the following
instantaneous vector is available
for feedback purposes:
\begin{equation}\label{eq:fback-quantities}
{\mathbf e(t)}=\big(e^{(0)}(t),\ldots, e^{(\hat r-1)},
y^{(\hat r)}(t),\ldots,y^{(r-1)}(t)\big)\in\R^{rm},
\quad
{e(t):=y(t)-y_{\textrm{ref}}(t),}
\end{equation}
with the notational convention that $e^{(0)}\equiv e$.\\[1ex]
{\it Feedback strategy.}\ \
Preliminary ingredients in the feedback construction,
{called {\it funnel control design parameters},} are:
\begin{equation}\label{eq:fcts-FC}
\left.
\begin{array}{l}
{\varphi\in\Phi,~~\text{bounded if \ $\hat r < r$},}\\[.6ex]
N\in C(\R_{\ge 0},\R),~~\text{a surjection,}
\\[.6ex]
\alpha\in C^1( [0,1), [1,\infty)),~~\text{a bijection.}
\end{array}\right\}
\end{equation}
{These functions are open to choice.} For notational convenience, we {define}
\begin{equation}\label{eq:fcts-FC2}
\gamma\colon\cB\to\R^m, \ w\mapsto  \alpha (\|w\|^2)\, w, \qquad \text{where \  $\cB :=\setdef{w\in\R^m }{ \|w \| < 1}$.}
\end{equation}
Next, we introduce continuous maps $\rho_k\colon\cD_k\to \cB$, $k=1,\ldots,r$, recursively as follows:
\begin{equation}\label{eq:fcts-FC3}
\left.
\begin{array}{l}
\cD_1:=\cB,\quad \rho_1\colon\cD_1\to\cB,~\eta_1\mapsto\eta_1,
\\[0.6ex]
\cD_k:= \setdef{(\eta_1,\ldots,\eta_k)\in\R^{km}}{\begin{array}{l} (\eta_1,\ldots,\eta_{k-1})\in\cD_{k-1},\\
 \eta_k+\gamma(\rho_{k-1}(\eta_1,\ldots,\eta_{k-1}))\in\cB\end{array}},
\\[1.9ex]
\rho_k\colon\cD_k\to\cB,~~(\eta_1,\ldots,\eta_k)\mapsto \eta_k +\gamma (\rho_{k-1}{ (\eta_1,\ldots,\eta_{k-1})}).
\end{array}
\right\}
\end{equation}
Note that each of the sets~$\cD_k$ is non-empty and open.
With reference to Fig.~\ref{Fig:funnel-controller}, and with~$\mathbf{e}$ and~$\rho_r$ defined by~\eqref{eq:fback-quantities}
and~\eqref{eq:fcts-FC3},
the \textit{funnel controller} is given  by
\begin{equation}\label{eq:FC}
\framebox{$
u(t)=\big(N\circ \alpha\big)(\|w(t)\|^2) \, w(t),\qquad w(t):=\rho_r\big(\varphi (t)\mathbf{e}(t)\big).$
}
\end{equation}
{ Note the striking simplicity of the control~\eqref{eq:FC}:
proportional feedback of the information vector~$w(t)$, with scalar gain.  Further comments on its distinctive features
are expounded in
{Example~\ref{Ex:design_pars} and} Remarks~\ref{Rem:FC} {\&}~\ref{hgp-rem-rev}
 below.}

\captionsetup[subfloat]{labelformat=empty}

\begin{example} \label{Ex:design_pars}
Choosing the design parameter triple
\[
{\varphi\in\Phi,} \quad N\colon s\mapsto  s \sin ( s),\quad
\alpha\colon s\mapsto 1/(1-s),
\]
{(with $\varphi$ bounded if $\hat r < r$),}
the feedback becomes
\[
u(t) \ = \
\big(1-\|w(t)\|^2\big)^{-1} \sin\left(  \big(1-\|w(t)\|^2\big)^{-1}\right) \cdot w(t),
\]
where the signal~$w(t)$ is, for example,
\[
w(t)=\begin{cases}
\varphi(t)e(t), & \text{if $r=1=\hat r$,}
\\
\varphi(t)\dot y(t)+\gamma (\varphi(t)e(t))
, & \text{if $r=2,\,\hat r=1$,}
\\
\varphi(t)\dot e(t)+\gamma (\varphi(t)e(t))
, & \text{if $r=2,\,\hat r=2$,}
\\
\varphi(t)\ddot y(t)+\gamma \big(\varphi(t)\dot y(t)+\gamma(\varphi(t)e(t))\big), & \text{if $r=3,\,\hat r =1$},
\end{cases}
\]
with~$\gamma$ given by~\eqref{eq:fcts-FC2}.
\end{example}

\begin{remark}\label{Rem:FC}
\begin{enumerate}[(a)]
\item 
The  intermediate  signal~$w(t)$  in~\eqref{eq:FC}
is a feedback~-- via the function~$\gamma$~-- of the available information, given by~\eqref{eq:fback-quantities},
``weighted'' by~$\varphi(t)$.
\item %

We point out that the complexity of the controller is  much lower than in previous approaches such as~\cite{BergLe18}, where successive derivatives of
auxiliary error variables\footnote{The auxiliary error variables are given by~$e_i(t)$ in equation~(5) of~\cite{BergLe18} for $i=0,\ldots,r-1$.}  need to be calculated before implementation. This complicates the feedback structure {for larger values of the parameter~$r$.}
In~\eqref{eq:FC} all required signals are explicitly given by the recursion in~\eqref{eq:fcts-FC3} and can be implemented directly.
\item 
The  parameter~$\hat r\in\{1,\ldots,r\}$
specifies the number of  derivatives  of~$y_{\textrm{ref}}$
available for feedback.
With increasing~$\hat r$, more information becomes accessible and so, not unreasonably, it might be expected that, loosely speaking,
controller performance improves: this expectation is  borne out by numerical simulations in Section~\ref{Sec:Sim}.
\item
Note that, if~$\hat r=r$, then {polynomial or} exponentially increasing funnel functions~$\varphi$ are admissible. For example, the {choices $\varphi\colon t\mapsto a t^\ell$ or} $\varphi\colon t\mapsto e^{at}-1$,
$a >0$, {$\ell\in\N$,}
ensure {polynomial/exponential} decay (to zero) of the tracking error $t\mapsto e(t)=y(t)-y_{\rm ref}(t)$.
 If~$\hat r < r$, then boundedness of~$\varphi$ is required.
As an exemplar in this case, the choice $\varphi\colon t\mapsto\min \{e^{at}-1\,,\,b\}$, $a,b >0$, ensures that the tracking error approaches the ball of
(arbitrarily small) radius~$b^{-1}$ exponentially fast and resides in that ball for all~$t \geq a^{-1}\ln (1+b)$.

\item
Funnel control presents an anomaly:  its performance might seem to  contradict the {\it internal model principle} which asserts that ``a regulator is
structurally stable only if the controller~[\ldots]   incorporates~[\ldots]
 a suitably reduplicated model of the dynamic structure of the exogenous signals which
the regulator is required to process''~\cite[p.\,210]{Wonh79}.   Diverse sources echo this principle~-- one such source is noted in\cite{HuIsi18}:
a young Mark Twain,
when apprenticed to a
Mississippi river pilot, recorded the latter's advice on navigating the river in the words ``you can always steer by the shape that's in your head, and
never mind the one that's before your eyes''~\cite[Ch.VIII]{Twain}.   But the funnel controller has no ``shape'' in its ``head'', it operates only on what is
before its eyes.  It does not incorporate ``a suitably reduplicated
model~[\ldots]  of the exogenous signals''.  How is this anomaly to be resolved?
The internal model principle applies in the context of {\it exact} asymptotic tracking of reference signals.
In the case of a {\it bounded} funnel function~$\varphi$, only {\it approximate} tracking, with non-zero prescribed asymptotic accuracy, is assured
{ in which case the anomaly is spurious}.

\item
But what of the case of an unbounded funnel function~$\varphi$, which is permissible whenever~$\hat r =r$?    In this case, exact asymptotic tracking is achieved.
({ See also the paragraph `Practical and exact asymptotic tracking' in Section~\ref{Ssec:Novelties}.})
Returning to the control-theoretic origins of the internal model principle, summarised in~\cite[p.\,210]{Wonh79} as ``every good regulator must incorporate a
model of the outside world'', we regard the term ``good regulator'' as  most pertinent.  A fundamental ingredient of the funnel controller is the quantity~$\varphi (t){\mathbf{e}}(t)$
which, in the case of unbounded~$\varphi$, inevitably leads to an ill-conditioned computation of the product of ``infinitely large'' and ``infinitesimally small''
terms.    Such a controller cannot be deemed ``good''.   Whilst of theoretical interest, the case of unbounded~$\varphi$ is of limited practical utility.
\end{enumerate}
\end{remark}
\begin{remark}\label{hgp-rem-rev}
We comment on the function $N\in C(\R_{\ge 0},\R)$ in~\eqref{eq:fcts-FC}.
\begin{enumerate}
\item[(a)]
Note that~$N$ is a surjection if, and only if,
\begin{equation}\label{Nprops}
\limsup_{s\to\infty} N(s)=+\infty\quad\text{and}\quad \liminf_{s\to\infty} N(s)=-\infty.
\end{equation}
These two conditions are a generalization of { the so-called \textit{Nussbaum properties} (to be discussed further in Section~\ref{Ssec:contr_dir}
below)}.
Reiterating Remark~\ref{hgp-rem}, the high-gain property implies that, for every pair $(K_p,K_q)$ of compact sets, {at least} one of the
{conditions in~\eqref{eq:chi_+- } must hold.}
{ In the
absence of any further {\em a priori} knowledge pertaining to these two possibilities}, the role of the function~$N$ is to provide the controller with
{ a ``probing'' capability which implicitly accommodates} each possibility.
\item[(b)]
If it is known
{ {\em a priori} that~$f$ has the negative-definite high-gain property}, then
{ $N$ may be replaced by any continuous surjection $[0,\infty)\to [0,\infty)$, the simplest example being the identity map~$s\mapsto s$
in which case} the feedback in~\eqref{eq:FC} takes the form
$u(t) = \alpha (\|w(t)\|^2) \, w(t)$.
\\
Similarly, { if it is known {\em a priori} that~$f$ has the positive-definite high-gain property}, then
{ $N$ may be replaced by any continuous surjection $[0,\infty)\to (-\infty,0]$, the simplest example being the  map $s\mapsto -s$
in which case} the feedback takes the form
$u(t) =  -\alpha (\|w(t)\|^2) \, w(t)$.

\item[(c)]
As the example in Remark~\ref{hgp-rem}\,(c) shows, it is also possible that~(i) and~(ii) in~\eqref{eq:chi_+- } hold simultaneously, in which case both
{ of the above simplified} feedback laws are feasible.
To illustrate this, consider the { scalar} system
\[
    \dot x(t) = u(t)\, \sin\big(\ln(1+|u(t)|)\big),\quad x(0) = 1,
\]
under control~\eqref{eq:FC} with $\alpha:s\mapsto 1/(1-s)$ and $N$ { replaced by}  $N: s \mapsto \sigma s$, where $\sigma\in\{-1,1\}$, that is
\[
    u(t) = \frac{\sigma w(t)}{1 - w(t)^2},\quad w(t) = \varphi(t) \big(y(t) - y_{\textrm{ref}}(t)\big).
\]
We choose $\varphi(t) = t^2$, $y_{\textrm{ref}}(t) = \sin t$ for $t\ge 0$ and perform the simulation\footnote{All simulations in the paper
are MATLAB generated (solver: {\tt ode45}, rel.\ tol.: $10^{-14}$, abs.\ tol.: $10^{-10}$).} over the time interval $[0,10]$. The results are shown in Fig.~\ref{fig:sim-N}, where the tracking error and input function for $\sigma=-1$ are depicted in Figs.~\ref{fig:simN-e1} and~\ref{fig:simN-u1}, and for $\sigma=1$ in Figs.~\ref{fig:simN-e2} and~\ref{fig:simN-u2}, resp. In the latter case, the input exhibits a sharp increase when the tracking error approaches the funnel boundary, and it stays within the interval $[20,25]$ thereafter, while for $\sigma=-1$ the input stays within the interval $[-1.5,1.5]$. This suggests that the system structure allows the input to ``probe'' for an appropriate interval of control values, independent of the sign of~$\sigma$.
\captionsetup[subfloat]{labelformat=empty}
\begin{figure}[h!tb]
  \centering
  \subfloat[Fig.~\ref{fig:sim-N}a: Funnel and tracking error for $\sigma = -1$]
{
\centering
\hspace{-2mm}
  \includegraphics[width=0.53\textwidth]{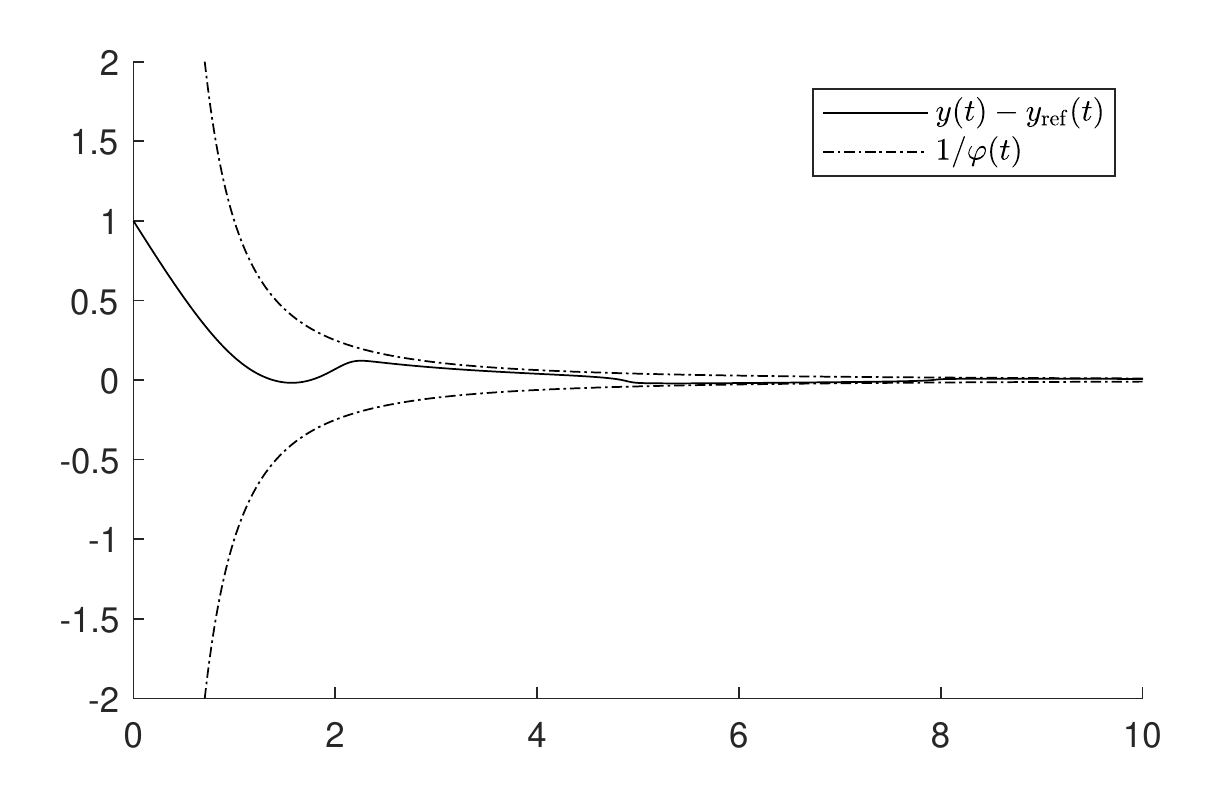}
\label{fig:simN-e1}
}
\subfloat[Fig.~\ref{fig:sim-N}b: Funnel and tracking error for $\sigma = +1$]
{
\centering
\hspace{-5mm}
  \includegraphics[width=0.53\textwidth]{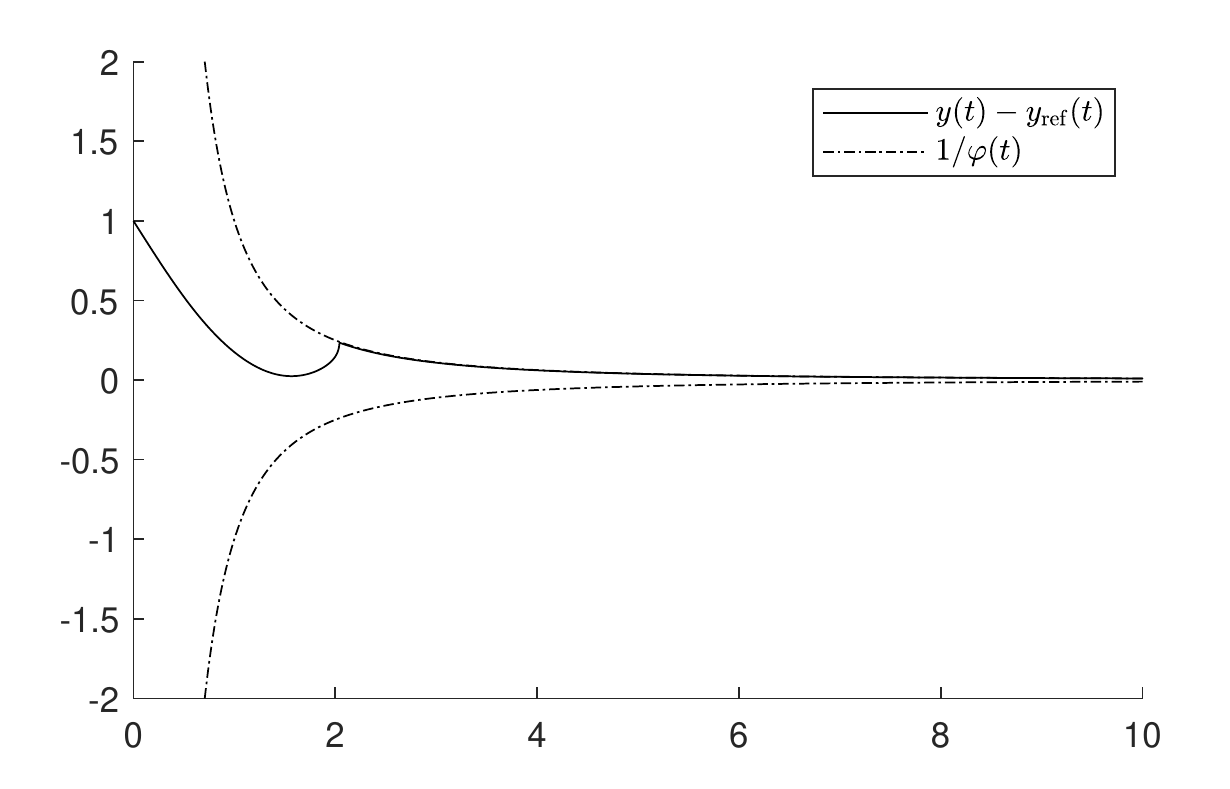}
\label{fig:simN-e2}
}\\
\vspace{5mm}
\subfloat[Fig.~\ref{fig:sim-N}c: Input function for $\sigma = -1$]
{
\centering
\hspace{-2mm}
  \includegraphics[width=0.53\textwidth]{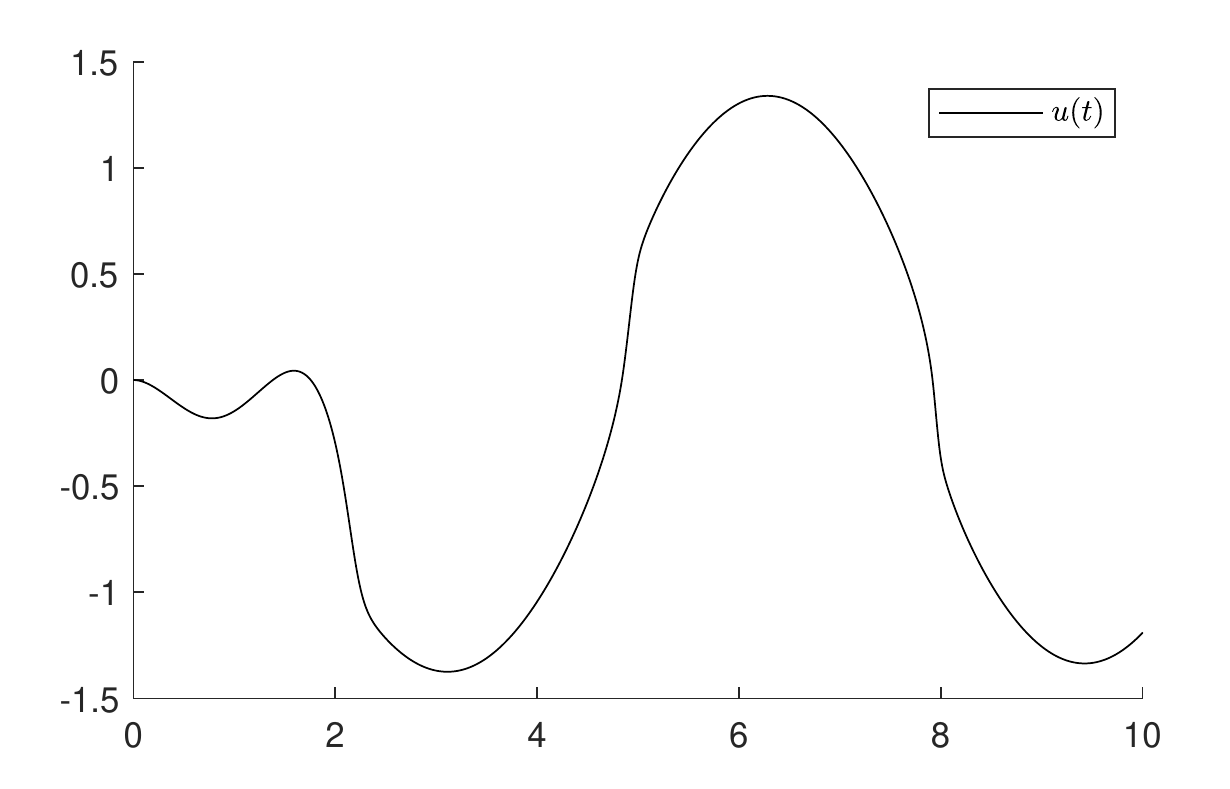}
\label{fig:simN-u1}
}
\subfloat[Fig.~\ref{fig:sim-N}d: Input function for $\sigma = +1$]
{
\centering
\hspace{-5mm}
  \includegraphics[width=0.53\textwidth]{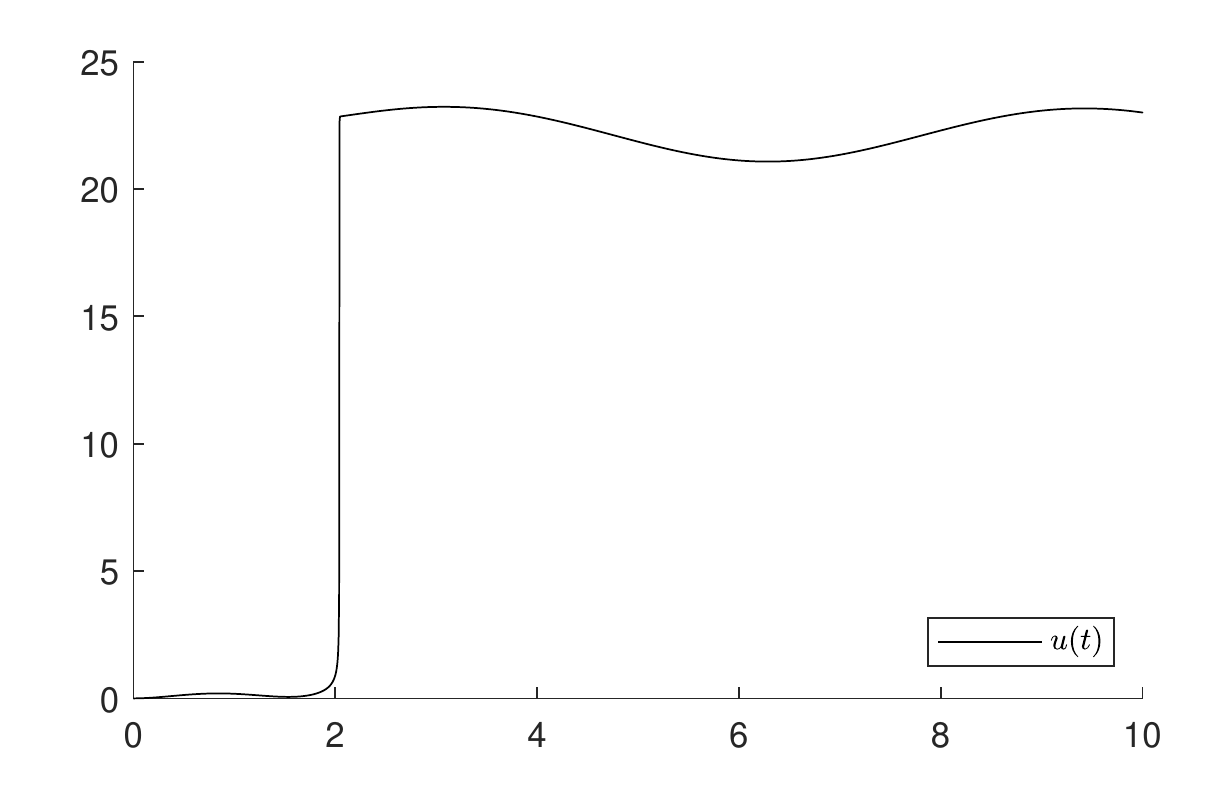}
\label{fig:simN-u2}
}
\caption{Simulation of the example from Remark~\ref{hgp-rem}\,(a) under control~\eqref{eq:FC} with $N:s \mapsto \sigma s$.}
\label{fig:sim-N}
\end{figure}
\end{enumerate}
\end{remark}

\subsection{\bf Funnel control~-- main result}\label{Sec:Feasib}
%
If the funnel controller~\eqref{eq:FC}  is applied to a system~\eqref{eq:nonlSys}, then
the first issue is to prove the existence of solutions of the closed-loop initial-value problem
 and to establish the efficacy of the control. We stress that the proof is quite delicate~-- even in the case of linear systems of the form~\eqref{eq:ABC}.
The reason is that the function~$\alpha$ used in the feedback~\eqref{eq:FC} introduces a potential singularity on the right hand side of the closed-loop differential equation.

By a \textit{solution} of~\eqref{eq:nonlSys},~\eqref{eq:FC} on~$[-h,\omega)$ we mean a function $y\in C^{r-1}([-h,\omega),\R^m)$, $\omega\in(0,\infty]$, with $y|_{[-h,0]} = y^0$ such that $y^{(r-1)}|_{[0,\omega)}$ is locally absolutely continuous and satisfies the differential equation in~\eqref{eq:nonlSys} with~$u$ defined in~\eqref{eq:FC} for almost all $t\in[0,\omega)$; $y$ is
said to be \textit{maximal}, if it has no right extension that is also a solution.

We are now in the position to present the main result
for systems belonging to the system class~$ \cN^{m,r}$.

\begin{theorem}\label{Thm:FunCon-Nonl}
Consider system~\eqref{eq:nonlSys} with $(d,f,\fT)\in \cN^{m,r}$, $m,r\in\N$, and initial data $y^0\in C^{r-1}([-h,0],\R^m)$. { Let $y_{\mathrm{ref}}\in W^{r,\infty}(\R_{\ge 0},\R^m)$, with associated parameter $\hat r\in\{1,\ldots,r\}$,} be arbitrary. Choose the triple~$(\alpha,N,\varphi)$ of funnel control design parameters as in~\eqref{eq:fcts-FC}. Assume that the instantaneous vector~${\mathbf{e}}(t)$, given by~\eqref{eq:fback-quantities}, is available for feedback and the following holds:
\begin{equation}\label{ic}
\varphi(0){\mathbf{e}}(0)\in \cD_r,
\end{equation}
(trivially satisfied if $\varphi(0)=0$).\\
%
Then the funnel control~\eqref{eq:FC} applied to~\eqref{eq:nonlSys} yields an initial-value problem which has a solution, every solution can be maximally extended and every maximal solution
$y:\left[-h,\omega\right)\rightarrow \R^m$ has the properties:
\begin{enumerate}
\item \label{item-main-1}
$\omega=\infty$ (global existence);
\item \label{item-main-2}
$u\in\cL^\infty (\R_{\ge0},\R^m)$, $y\in\cW^{r,\infty}([-h,\infty),\R^m)$;
\item \label{item-main-3}
the tracking error $e\colon \R_{\ge 0}\to \R^m$ as in~\eqref{eq:fback-quantities} evolves in the funnel~$\mathcal{F}_{\varphi}$  and is uniformly bounded away
from the funnel boundary
\[
\partial \cF_{\varphi} =\setdef{(t,\zeta)\in\R_{\ge 0}\times\R^m}{\varphi(t)\|\zeta\|=1}
\]
in the sense that there exists $\varepsilon \in (0,1)$ such that
$\varphi(t)\|e(t)\|\leq \varepsilon$ for all $t\geq 0$.
\item \label{item-main-4}
If { $\hat r=r$} and~$\varphi$ is unbounded, then $e^{(k)}(t)\to 0$ as $t\to\infty$, \  { $k=0,\ldots, r-1$}.
\item \label{item-main-5}
If the system is known to satisfy the
 negative-definite (respectively, positive-definite) high-gain property
 (see  Remark~\ref{hgp-rem}\,\eqref{item:13b}),
 then the feedback~\eqref{eq:FC} may be simplified { by substituting the identity map $s\mapsto s$ (respectively,
the map $s\mapsto -s$) for $N$} and Assertions  {\rm (i)--(iv)} remain valid.
\end{enumerate}
\end{theorem}

\noindent The proof is relegated to Appendix~\ref{Sec:App}.

When interpreted in specific cases, the initial condition constraint~\eqref{ic} becomes more transparent.
For example, in the relative-degree-one case $r=1=\hat r$, it is simply the requirement that
$\varphi (0)\|e(0)\|<1$, where $e(0)=y^0(0)-y_{\textrm{ref}}(0)$ and, in the case $r=2$, it is equivalent to the same requirement augmented by
\[
\|\varphi(0)z +\gamma (\varphi(0)e(0))\| < 1,\ \ \text{with}\ z=\begin{cases}\dot{y}^0(0)-\dot{y}_{\textrm{ref}}(0), & {\text{if}}~ \hat r=2,\\
\dot y^0(0), & \text{if}~ \hat r=1.\end{cases}
\]

In some specific circumstances, computation of {\it a priori} bounds on the evolution of the tracking error $e$ and (some of) its derivatives is possible.
We highlight one such circumstance.  Assume that $\hat r \geq 2$ and $\varphi\in\Phi$ is such that $\varphi (0)>0$.  Define
\[
\mu_0 := \esup_{t\ge 0}\big(|\dphi (t)|/\varphi(t)\big).
\]
Let $\alpha^\dag\in C^1(\R_{\ge 0},[0,1))$
denote the inverse of the continuously differentiable bijection $[0,1)\to\R_{\ge 0}$, $s\mapsto s\alpha(s)$
and, for notational convenience, introduce the continuous function
\[
\tilde\alpha\colon~[0,1)\to\R_{\ge 0},\ s\mapsto 2s\alpha^\prime (s)+\alpha (s).
\]
Define $(\mu_k,e_k^0,c_k)$, $k=1,\ldots,\hat r-1$, recursively as follows:
\begin{equation}\label{ckbounds}
\left.
\begin{aligned}
e_1^0&:= \varphi (0)e(0),
~~c_1:= \max\{\|e_1^0\|^2,\alpha^\dag (1+\mu_0)\}^{1/2} < 1,
~~\mu_1:=1+\mu_0c_1,\ \ \\[1ex]
\mu_k&:=1+\mu_0\big(1+c_{k-1}\alpha (c_{k-1}^2)\big)
+\tilde\alpha (c_{k-1}^2)\big(\mu_{k-1}+c_{k-1}\alpha (c_{k-1}^2)\big),\\[1ex]
e_k^0 &:= \varphi (0)e^{(k-1)}(0)+\alpha (\|e_{k-1}^0\|^2) e_{k-1}^0,
\\[1ex]
c_k&:= \max\{\|e_k^0\|^2,\alpha^\dag (\mu_k)\}^{1/2}  < 1.
\end{aligned}
\right\}
\end{equation}
We emphasize that the constants $c_k$ are determined by the design parameters~$\varphi$ and~$\alpha$, together with {the} known
initial data: $y(0),\ldots,y^{(\hat r-1)}(0)$ and $ y_{\textrm{ref}}(0),\ldots,y_{\textrm{ref}}^{(\hat r-1)}(0)$.

\begin{corollary}\label{Prop:FunCon-Bounds}
Let all hypotheses of Theorem \ref{Thm:FunCon-Nonl} hold.  Assume, in addition, that
\[
\hat r \geq 2,\quad \varphi (0) >0\quad \text{ and}~\alpha^\prime ~\text{is monotonically non-decreasing.}
\]
Then, for every maximal solution $y:[-h,\infty)\to\R^m$ of the feedback system \eqref{eq:nonlSys}~\&~\eqref{eq:FC}, the tracking error~$e=y-y_{\rm ref}$
and its first $\hat r-2$ derivatives
satisfy{, for all~$k=1,\ldots,\hat r-2$ and all~$t\ge 0$,}
\[
\|e(t)\|\leq\varphi(t)^{-1}{c_1}, \quad
{\|e^{(k)}(t)\|\leq \varphi(t)^{-1}\big({c_{k+1}}+c_{k}\alpha (c_{k}^2)\big).}
\]
where the constants $c_k$ are given by \eqref{ckbounds}.
\end{corollary}
\noindent The proof is relegated to Appendix~\ref{Sec:App}.

{Note that} these findings are much simpler than the complicated bounds derived in~\cite[Prop.~3.2]{BergLe18}.

\begin{example}
Assume $\hat r=3$, $\varphi\colon t\mapsto a+bt$, $a,b >0$, and $\alpha\colon s\mapsto 1/(1-s)$.  In this case, we have
$\mu_0 =b/a$  and $\alpha^\dag\colon s\mapsto s/(1+s)$.   Therefore, {for all $t\ge 0$},
\[
e_1^0 = a e(0),\quad
c_1 = \max\{\|e_1^0\|^2\,,\, (1+\mu_0)/(2+\mu_0)\}^{1/2}\quad\text{and}\quad
\|e(t)\|\leq \frac{{c_1}}{a+bt}.
\]
Furthermore,  $\mu_1=1+\mu_0c_1$,
\begin{align*}
&\tilde\alpha\colon s\mapsto (1+s)/(1-s)^2,\quad \mu_2 = 1+ \mu_0\big(1+c_1\alpha (c_{1}^2)\big)+\tilde\alpha (c_{1}^2)\big(\mu_1
+c_1\alpha (c_{1}^2)\big),\\
&e_2^0=a\dot e(0)+\big(1-\|e_1^0\|^2\big)^{-1}e_1^0,\quad
c_2 = \max\{\|e_2^0\|^2\,,\, \mu_2/(1+\mu_2)\}^{1/2}\\
&\text{and}\quad {\forall\, t\ge 0:\  \|\dot e(t)\|\leq \frac{c_2+{c_1}/\big({1-c_1^2}\big)}{a+bt}.}
\end{align*}
\end{example}

\section{\bf { Examples}}\label{Sec:Examples}

In this section, we show that the class ~$\cN^{m,r}$ encompasses the prototype of linear multi-input multi-output systems with strict relative degree~$r\in\N$ and asymptotically stable zero dynamics, see Subsection~\ref{Ssec:LinSys}.
Furthermore, the issues of control directions are discussed in Subsection~\ref{Ssec:contr_dir} and input nonlinearities in Subsection~\ref{Ssec:InputNonl}; a special case of the latter is a so called dead-zone input which is discussed in Subsection~\ref{Ssec:DeadZone}.

\subsection{\bf The prototypical linear system class}\label{Ssec:LinSys}
%
As a concrete example we consider linear, finite-dimensional systems of the form
\begin{equation} \label{eq:ABC}
 \left.
 \begin{array}{l}
  \dot{x}(t) = A\,  x(t) + B\, u(t), \quad x(0)=x^0\in\R^n ,
  \\[1mm]
    y(t) =C\, x(t)
 \end{array}
\right\}
\end{equation}
where~$(A,B,C)\in\R^{n\times n} \times \R^{n\times m} \times \R^{m\times n}$,  $m\leq n$,
and discuss its relationships to Properties~(TP1)--(TP3) and the high-gain property.

\subsubsection{\bf Strict relative degree}
%
We show that   system~\eqref{eq:ABC} can be equivalently written in the form~\eqref{eq:nonlSys},
if  system~\eqref{eq:ABC}
has {\it (strict) relative degree} $r  \in \N$, that is
\[
CA^kB=0,~~k=0, \ldots ,r - 2
\qquad \text{and} \qquad
\Gamma :=CA^{r -1}B \quad  \text{is invertible.}
\]
It is shown in~\cite{IlchRyan07} that under this assumption there exists a state space transformation
\[
z = \begin{pmatrix}\xi\\\eta\end{pmatrix} = Ux,\quad \xi = \begin{pmatrix} \xi_1\\ \vdots\\ \xi_r\end{pmatrix},\quad U\in\R^{n\times n}\ \text{invertible},
\]
which transforms~\eqref{eq:ABC} into \textit{Byrnes-Isidori form}
\begin{equation*}\label{eq:ABC-BIF}
    \dot z(t) = \widetilde A z(t) + \widetilde B u(t),\quad y(t) = \widetilde C z(t),
\end{equation*}
where
\begin{equation}\label{eq:BIF}
    (\widetilde A,\widetilde B,\widetilde C)=(UAU^{-1}, UB, CU^{-1})
\end{equation}
with
\begin{align*}
\widetilde A &=
\begin{bmatrix}
0 & I_m & 0 & \cdots & 0 & 0 \\
0 & 0 & I_m &        &  & 0 \\
\vdots & & \ddots & \ddots & & \vdots \\
0 & 0 & \cdots & 0 & I_m & 0 \\
R_1 & R_2 & \cdots & R_{r -1} & R_{r } & S \\
P & 0 & \cdots & 0 & 0 & Q
\end{bmatrix},\quad
\widetilde B =\begin{bmatrix} 0_{m\times m}\\  \vdots \\0_{m\times m} \\\Gamma\\0_{(n-rm)\times m}\end{bmatrix},\\
\widetilde C &=\begin{bmatrix} I_m, & 0_{m\times m}, &\ldots\ , &0_{m\times m}, & 0_{m\times (n-rm)}\end{bmatrix}.
\end{align*}
In the new coordinates,  the system representation of~\eqref{eq:ABC}  becomes
{
\begin{equation}\label{eq:ABC_BIform}
\left.
\begin{array}{l}
\dot\xi_k (t)= \ \xi_{k+1}(t),~~k=1,\ldots, r-1,\\[1ex]
\dot\xi_r(t)= \ \sum_{k=1}^r R_k\xi_k(t)+S\eta (t) +\Gamma u(t),\\[1ex]
              \dot\eta(t)=  P\xi_1(t)+Q\eta (t),
\end{array}\right\}\quad\text{with output}~~ y(t)= \xi_1(t).
\end{equation}
}
With the {third} equation in~\eqref{eq:ABC_BIform}, the so-called internal  dynamics, we may associate a linear operator

\begin{equation}\label{eq:L}
L\colon  y(\cdot) \mapsto \left( t \mapsto
\int_0^t  {\rm e}^{Q(t-\tau)} Py(\tau)\, \dd \tau \right).
\end{equation}
With initial data $\eta(0)= {\eta^0=} [0, I_{n-rm}] U x^0$ and
$d(\cdot) := {\rm e}^{Q \cdot} \eta^0$, we find that
\[
 \eta(t) = d(t) +L(y)(t).
\]
Introducing the (linear) operator
\begin{equation}\label{eq:T}
\begin{array}{rcl}
\fT \colon \cC(\R_{\ge 0},\R^{rm })&\to&
\cL^\infty_{\text{loc}}(\R_{\ge 0},\R^m),\\[1ex]
\zeta =(\zeta_1,\ldots,\zeta_r)
&\mapsto&
\Big(  t\mapsto \sum_{k=1}^{r}  R_k  \zeta_k(t) +{S}  L(\zeta_1)(t) \Big),
\end{array}
\end{equation}
it follows from~\eqref{eq:ABC_BIform}  that~\eqref{eq:ABC} is equivalent to the functional differential system
\begin{equation}\label{eq:ABC-functional}
\left.\begin{array}{l}y^{(r )}(t)  =  {S} d(t)+\fT(y,\ldots,y^{(r-1)})(t)+\Gamma u(t)
\\[1ex]
y(0)=Cx^0, \ \ldots \  ,y^{(r-1)}(0)=CA^{r-1}x^0.
\end{array}\right\}
\end{equation}
It is easy to see that the operator~$\fT$ satisfies properties~(TP1) and~(TP2) from Definition~\ref{Def:Operator_class}. The following section is devoted to~(TP3).

\subsubsection{\bf Minimum phase}\label{Ssec:min-ph}
%
Suppose that system~\eqref{eq:ABC} has strict relative degree $r\in\N$.
Then  the BIBO property~(TP3)
of the operator~$\fT$ in~\eqref{eq:ABC-functional}
is closely related to
system~\eqref{eq:ABC} having \textit{asymptotically stable   zero dynamics}, i.e.,
\begin{equation}\label{eq:asy-st-ZD}
    \forall\,\lambda\in\C_{\ge 0}:\ \det \begin{bmatrix} \lambda I -A& B \\ C & 0\end{bmatrix}\neq 0.
\end{equation}
This concept (also closely related to the \emph{minimum phase} property in the literature,
cf.~\cite{IlchWirt13})
is extensively studied since its relevance has been revealed in classical works such as~\cite{ByrnWill84, Mare84}.
%
To be precise, assume that the transfer function $C(sI-A)^{-1}B \in\R(s)^{m\times m}$ of $(A,B,C)$ is invertible over~$\R(s)$, then we have the following:\\[-10mm]
\begin{center}
\[
\begin{array}{lcl}
\begin{array}{l}
\text{$(A,B,C)$ satisfies~\eqref{eq:asy-st-ZD}}
\end{array}
& \quad \stackrel{\text{\cite[Cor.~3.3]{Berg14c}}}{\Longleftrightarrow}  &  \quad
\begin{array}{l}
\text{$(A,B,C)$ stabilizable \& detectable,}\\
\text{$C(sI-A)^{-1}B$ has no zeros in $\C_{\ge 0}$}
\end{array}
\\[4mm]
&& \hspace*{3cm}\Updownarrow\text{\cite[Cor.~2.8]{Berg14c}}\\[2mm]
\begin{array}{l}
\text{$(A,B,C)$ stab.\ \& det.,}\\
\text{$\fT$ satisfies~(TP3)}
\end{array}
& \quad \stackrel{\text{\cite[Thm.~3.21]{TrenStoo01}}}{\Longleftrightarrow}  & \quad
\begin{array}{l}
\text{$(A,B,C)$ stabilizable \& detectable,}\\
\text{$S(sI-Q)^{-1}P$ has no poles in $\C_{\ge 0}$}
\end{array}
\end{array}
\]
\end{center}
%
%
For the last equivalence above we note that by~\cite[Thm.~3.21]{TrenStoo01} it is straightforward that $S(sI-Q)^{-1}P$ having no poles in $\C_{\ge 0}$ is equivalent to $(Q,P,S)$ being externally stable or, in other words, the operator~$L$ from~\eqref{eq:L} satisfies~(TP3). It is easily seen that this is the same as $\fT$ satisfying~(TP3).

\subsubsection{\bf Sign-definite high-frequency  gain matrix}\label{Ssec:HFG}
%
We show that system \eqref{eq:ABC} satisfies the {\it high-gain property} (recall Definition~\ref{Def:high-gain}) if, and only if,
the {high-frequency gain matrix}~$\Gamma =CA^{r -1}B$
is sign definite.  Otherwise stated, we seek to establish {the following equivalence}:
{
\[
{\text{(a)\quad \eqref{eq:ABC} has the high-gain property}}
\qquad \Longleftrightarrow \qquad
\text{(b)}\quad  \forall\,v\in\R^m\backslash\{0\} : \ v^\top \Gamma v \neq 0.
\]
}
(a) $\implies$ (b): \
Assume (a).  {Let $v^*\in(0,1)$ be given and choose $K_p=\{0\}$, $K_q=\{0\}$. Write $A_m:=\setdef{v\in\R^m}{v^* \leq \|v\| \leq 1}$.}
{Suppose (b) is false.  Then there exists $\hat v\in A_m$ such that $\hat v^\top \Gamma \hat v=0$, thus
\[
    \forall\, s\in\R:\ \chi (s)=\min_{v\in A_m} \big(-sv^\top\Gamma v\big) \le -s\hat v^\top \Gamma \hat v=0,
\]
which contradicts~(a).}
\\[1ex]
(b) $\implies$ (a): \
Assume (b).   Then there exists $\sigma\in\{-1,1\}$ such that $\sigma\Gamma$ is positive definite.
{Let $G:=(\sigma/2)(\Gamma+\Gamma^\top)$
denote the symmetric part of~$\sigma\Gamma$
and let $\lambda_*>0$ be the smallest eigenvalue of~$G$.}
{Set $v^*=\half$,} choose compact $K_p\subset \R^p$ and $K_q\subset\R^q$ and define
\[
c_1:= \min\setdef{ v^\top (\delta +z)}{(\delta,z,v)\in K_p\times K_q\times A_m}.
\]
Then,
\[
\forall\, s\in\R:\ \chi (s)-c_1 \geq\min_{v\in A_m} \big(-sv^\top \Gamma v\big) =\min_{v\in A_m} \big(-s\sigma v^\top G v\big).
\]
Let $(s_n)$ be a real sequence
with~$\sigma s_n <0$ for all~$n\in\N$ and $\sigma s_n\to -\infty$ as~$n\to\infty$.
It follows that
\[
\forall\, n\in\N\ \forall\, v\in A_m: \
-\sigma s_n v^\top Gv \geq -\sigma s_n \lambda_*\|v\|^2 \geq -\frac{\sigma s_n \lambda_*}{4}
\]
and so we have
\[
\forall\, n\in\N:\ \chi(s_n)\geq c_1 -\frac{\sigma s_n \lambda_*}{4}.
\]
Therefore, $\chi(s_n)\to\infty$ as $n\to\infty$ and so (a) holds.

\subsection{\textbf{Known and unknown   control directions}}\label{Ssec:contr_dir}
%
For linear systems~\eqref{eq:ABC} with relative degree $r\in\N$
the notion of ``control direction'' is captured  by  the sign of the high-frequency gain matrix $\Gamma = CA^{r-1} B$ as discussed in Section~\ref{Ssec:HFG}. More precisely, if
{$\sigma \Gamma $ is positive definite}
for some $\sigma\in \{-1,1\}$, then~$\sigma$ is called the \textit{control direction}. If~$\sigma$ is \textit{known} and the system~\eqref{eq:ABC} has asymptotically stable zero dynamics, see~\eqref{eq:asy-st-ZD}, then it can be shown that the ``classical high-gain adaptive feedback''
\begin{equation}\label{eq:dotk=y2}
u(t)= - \sigma k(t) y(t),\quad
\dot k(t)= \|y(t)\|^2,
\end{equation}
{with $k(0) = k^0 \ge 0$}, applied to~\eqref{eq:ABC} yields a closed-loop system, where for any solution~$(x,k)$ we have that $x(t)\to 0$ as $t\to\infty$ and $k(\cdot)$ is bounded; see~\cite{ByrnWill84,Mare84,Mors83}.

For the case of \textit{unknown control direction}~$\sigma$, the adaptive stabilization was an obstacle over many years.
\textit{Morse}~\cite{Mors83} conjectured the non-existence of a smooth adaptive controller which stabilizes every linear single-input single-output system~\eqref{eq:ABC}, i.e.~$m=1$, under the assumption that~$\Gamma\neq 0$.
It was shown by \textit{Nussbaum} in~\cite{Nuss83} that this conjecture is false:
One has to incorporate a {``sign-sensing function''} in the feedback law~\eqref{eq:dotk=y2}
so that it becomes
\begin{equation}\label{eq:dotk=y2-unknown}
u(t)= -N(k(t)) y(t),\quad  \dot k(t)= \|y(t)\|^2,
\end{equation}
where the smooth function~$N\colon\R_{\ge 0}\to\R$
satisfies the so-called \textit{Nussbaum property}
\begin{equation}\label{eq:Nussbaum}
\forall\, k^0\ge 0: \ \
\sup_{k>k^0} \frac{1}{k-k^0} \int_{k^0}^k N(\kappa)\, {\rm d}\kappa = \infty
\quad
  \text{and}\quad  \inf_{k>k^0} \frac{1}{k-k^0} \int_{k^0}^k N(\kappa)\, {\rm d}\kappa = -\infty,
\end{equation}
see, for example, \cite{GeHong04, GeWang02, GeWang03, JianMare04, Ye01}.  Loosely speaking, when incorporated in the control design, ``Nussbaum'' functions provide a
mechanism that can ``probe'' in both control directions.

The present paper utilizes a larger class of ``probing'' functions: in particular, the proposed control design
permits the adoption of
any continuous function~$N: \R_{\ge 0}\to\R$ {which is surjective or, equivalently, satisfies~\eqref{Nprops}}. Properties~\eqref{eq:Nussbaum} imply properties~\eqref{Nprops}, but the reverse implication is false: for example,
the function $s\mapsto N(s) = s\sin s$ exhibits properties~\eqref{Nprops}, but fails to exhibit the Nussbaum properties~\eqref{eq:Nussbaum}.

\subsection{\textbf{Input nonlinearities}}\label{Ssec:InputNonl}
In addition to accommodating  the issue of (unknown)
control direction (cf.\ Section~\ref{Ssec:contr_dir}), the generic formulation~\eqref{eq:nonlSys} with associated high-gain property encompasses a wide variety of input nonlinearities.  Consideration of a scalar
system of the simple form
\begin{equation}\label{eq:f1f2}
 \dot y(t)=f_1(y(t))+f_2(y(t)) \ \beta (u(t))
\end{equation}
with $f_1\in C(\R,\R)$, $f_2\in C(\R,\R\backslash\{0\})$ and $\beta\in C(\R,\R)$, will serve to illustrate this variety.  The assumption that $f_2$ is a non-zero-valued continuous
function ensures a well-defined control direction (unknown to the controller). Without loss of generality, we may assume that~$f_2\in C(\R,\R_{>0})$;
if~$f_2$ is negative-valued, then, in~\eqref{eq:f1f2}, simply replace~$f_2$ by~$-f_2$ and~$\beta$ by~$-\beta$.
We impose the following conditions on~$\beta\in C(\R,\R)$:
\begin{equation}\label{betaprops}
\beta {\text{ is surjective, with}}~~|\beta (\tau)| \to \infty ~~{\text{as}}~|\tau|\to\infty,
\end{equation}
which is equivalent to the requirement that one of the following conditions hold:
\[
   \lim_{\tau\to\pm \infty} \beta(\tau) = \pm \infty\quad \text{or} \quad \lim_{\tau\to\pm \infty} \beta(\tau) = \mp \infty.
\]
We proceed to show that system~\eqref{eq:f1f2} has the high-gain property. {Set $v^*=\half$,}
let~$K_1\subset \R$ be compact and define
\[
A_1:= \big[-1,-\half\big]\cup \big[\half,1\big], \qquad
c_1:=\min\setdef{v\,f_1(z)}{ (z,v)\in K_1\times A_1  }\in\R.
\]
Consider the function
\[
\chi\colon~\R\to\R,~~s\mapsto\min\setdef{v\big(f_1(z)+f_2(z)\beta(-sv)\big) }{ (z,v)\in K_1\times A_1}.
\]
Then
\begin{equation}\label{eq:c1-min}
\forall \, s \in \R : \
\chi(s)\geq
c_1 +   \min\setdef{vf_2(z)\beta (-sv) }{ (z,v)\in K_1\times A_1}.
\end{equation}
Let~$M>0$ be arbitrary.  To conclude that the high-gain property holds, it suffices to show that there exists~$s\in\R$  such that
\[
\forall\,(z,v)\in K_1\times A_1: \
v f_2(z)\beta(-sv) > M.
\]
Define
\[
c_2:= \min_{z\in K_1} f_2(z) >0\qquad\text{and}\qquad c_3:= 2M/c_2.
\]
By properties of~$\beta$, there exist~$\sigma\in\{-1,1\}$ and~$c_4 >0$ such that
\[
\forall\,\tau > c_4 : \
\min\big\{\beta (\sigma\tau)\,,\,-\beta (-\sigma\tau)\big\}>c_3.
\]
Let~$(z,v)\in K_1\times A_1$ be arbitrary. Fix~$s\in\R$ such that $\sigma s<  -2c_4$ and so~$|sv| > c_4$.
Then
\[
v f_2(z)\beta(-sv) = \left\{\begin{array}{ll}|v|f_2(z)\beta(\sigma|sv|), & \text{if $ v>0$}
\\[1ex]
|v|f_2(z)\big(-\beta (-\sigma |sv|)\big), & \text{if $v < 0$}\end{array}\right\}
> \frac{c_2c_3}{2} = M.
\]
Therefore, the high-gain property holds.

\subsection{{\bf Dead-zone input}} \label{Ssec:DeadZone}
An important example of a nonlinearity
$\beta=D$ with properties~\eqref{betaprops} is a so-called \emph{dead-zone input}  of the form
\[
   D : \R\to \R, \quad v\mapsto
    D(v) = \left\{
    \begin{array}{rcl} D_r(v), && v\ge b_r,\\ 0,&& b_l<v<b_r,\\ D_l(v), && v\le b_l\end{array}\right.
\]
with unknown deadband parameters~$b_l<0<b_r$ and unknown functions $D_l,D_r\in\cC(\R,\R)$ which satisfy,
 for unknown $\sigma\in\{-1,1\}$,
\[
D_l(b_l) = D_r(b_r) = 0 \qquad \text{and} \qquad
\lim_{s\to\infty} \sigma D_r(s) = \infty, \quad \lim_{s\to-\infty} \sigma D_l(s) = -\infty.
\]
Note that the above assumptions allow for a much larger class of
functions~$D_l, D_r$ compared to e.g.~\cite{Na13}, where assumptions on their derivatives are used.
In particular,~in the present context, $D_l$ and~$D_r$ need not be differentiable or monotone.

\section{Simulations}\label{Sec:Sim}

We compare the controller~\eqref{eq:FC} to the controller presented in~\cite{BergLe18} and, to this end, consider the simulation examples presented therein.

\subsection{\bf Mass-on-car system}\label{Ssec:Mass_on_car}

\begin{figure}[htp]
    \begin{center}
    \includegraphics[trim=2cm 4cm 5cm 15cm,clip=true,width=6.5cm]{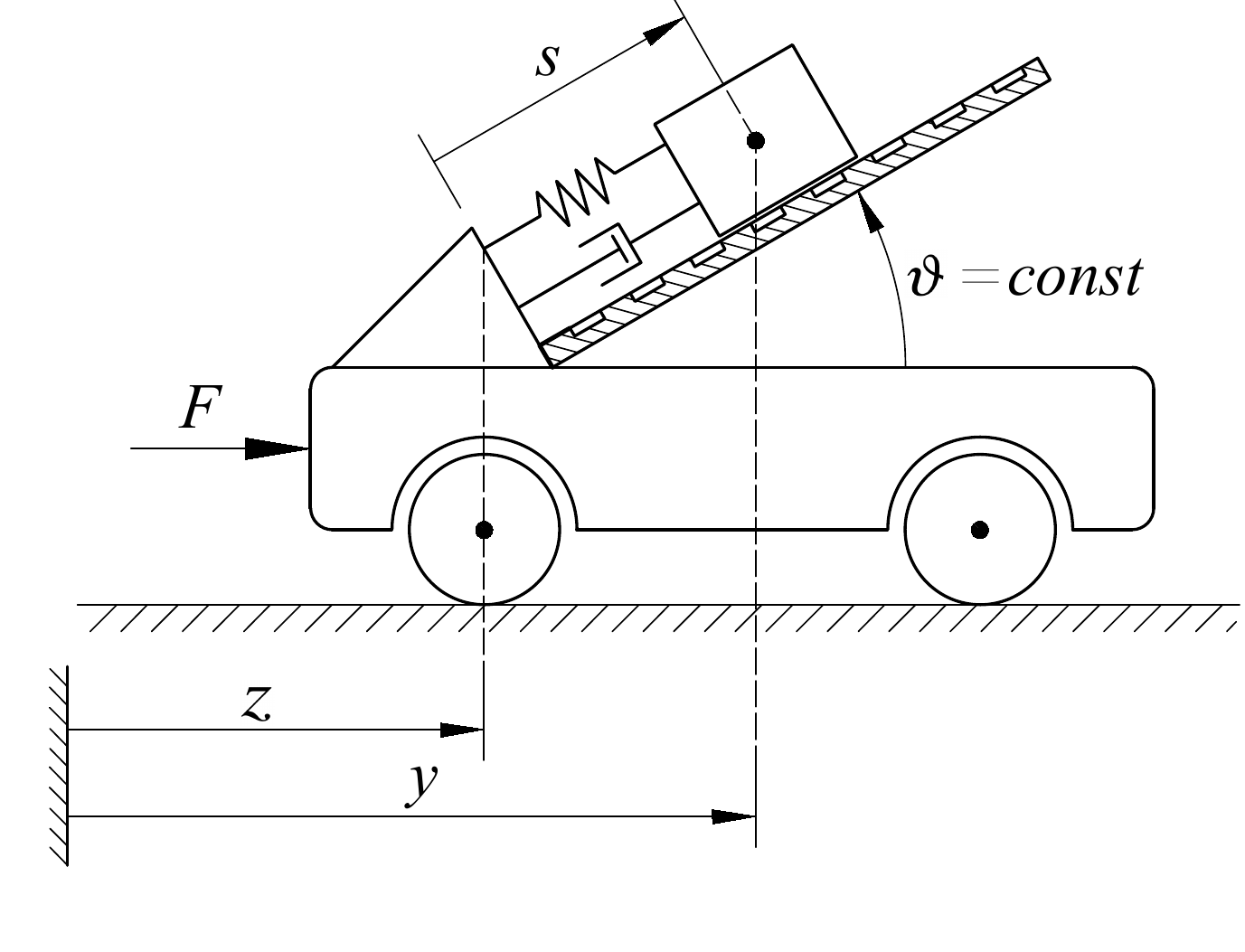}
    \end{center}
    \vspace*{-3mm}
    \caption{Mass-on-car system.}
    \label{Mass.on.car}
    \end{figure}

To illustrate the controller~\eqref{eq:FC}, we consider a mass-spring system mounted on a car from~\cite{SeifBlaj13}, see Fig.~\ref{Mass.on.car}. The mass~$m_2$ (in~\si{\kilo\gram})  moves on a ramp inclined by the angle~{$\vartheta \in [0,\frac{\pi}{2})$} (in \si{\radian}) and mounted on a car with mass~$m_1$ (in \si{\kilo\gram}), for which it is possible to control the force with~$u=F$ (in \si{\newton}) acting on it. The equations of motion for the system are given by
\begin{equation}\label{mass.on.car.equ}
\begin{bmatrix}
m_1+m_2&m_2\cos \vartheta\\
m_2\cos \vartheta&m_2
\end{bmatrix} \begin{pmatrix} \ddot{{ z}}(t)\\ \ddot{s}(t) \end{pmatrix} +\begin{pmatrix}
0\\
ks(t)+d\dot{s}(t)
\end{pmatrix}=\begin{pmatrix}
u(t)\\
0
\end{pmatrix},
\end{equation}
where $t$ is current time (in \si{\second}),  ${ z}$ (in \si{\metre}) is the horizontal car position and~$s$ (in \si{\metre}) the relative position of the mass on the ramp. The constants~$k { >0}$ (in \si{\newton\per{\metre}}), $d{ >0}$ (in \si{\newton\second\per{\metre}}) are the coefficients of the spring and damper, respectively. The output
 $y$ (in \si{\metre}) of the system is given by the horizontal position of the mass on the ramp,
\[
y(t)={ z}(t)+s(t)\cos \vartheta.
\]
{
Writing $\mu:=m_2\big(m_1+m_2 \sin^2 \vartheta\big)$, $\mu_1:=m_1/\mu$ and $\mu_2:=m_2/\mu$, it is readily verified that this system takes
the form \eqref{eq:ABC}, with
\[
x(t):=\begin{pmatrix}z(t)\\\dot z(t)\\s(t)\\\dot s(t)\end{pmatrix},~~
A:=\begin{bmatrix}0&~~1&0&0\\0&~~0&\mu_2k\cos\vartheta & \mu_2 d\cos\vartheta\\0&~~0&0&1\\0&~~0&~-(\mu_1+\mu_2)k &~-(\mu_1+\mu_2)d
\end{bmatrix},~~
B:=\begin{bmatrix}0\\\mu_2\\0\\-\mu_2\cos\vartheta\end{bmatrix}
\]
and $C:= \begin{bmatrix} 1 &0&\cos\vartheta &0\end{bmatrix}$.
Observe that
\[
CB=0,\quad CAB=\mu_2 \sin^2 \vartheta,\quad CA^2B=d\mu_1\mu_2\cos^2\vartheta
\]
and so the relative degree $r$ of the system is given by
\[
r=\left\{\begin{array}{ll}2,&~~\text{if}~\vartheta\in (0,\frac{\pi}{2})\\[1ex] 3, &~~\text{if}~\vartheta=0.\end{array}\right.
\]
Moreover, $CA^{r-1}B >0$ and so the positive-definite high-gain property holds.  Furthermore, a straightforward (if tedious) calculation
reveals that the {eigenvalues} of the matrix $Q\in \R^{(4-r)\times(4-r)}$ in the Byrnes-Isidori form \eqref{eq:ABC_BIform} are given by
\[
\lambda := -k/d~~~\text{in the case $r=3$}
\]
or, in the case $r=2$, by
\[
\lambda _\pm := -(\tilde d/2)\pm\sqrt{(\tilde d/2)^2 -\tilde k},\quad \tilde d:= d/(m_2 \sin^2 \vartheta),~~
 \tilde k:= k/(m_2 \sin^2 \vartheta).
 \]
Thus, in each case, the zero dynamics are asymptotically stable  and so
 property (TP3) holds for the associated operator ${\mathbf T}$ given by \eqref{eq:T}.
 Therefore, the system is of class $\cN^{1,r}$ to which the funnel control \eqref{eq:FC} applies.
}
Invoking Assertion\,(v) of Theorem~\ref{Thm:FunCon-Nonl}, the function~$N$ in~\eqref{eq:FC} may be { substituted by the map} $s\mapsto -s$.

For the simulation, we choose the parameters
$m_1=4$, $m_2=1$, $k=2$, $d=1$,
the initial values $x(0)= s(0) = 0$, $\dot{x}(0) = \dot s= 0$ and the reference trajectory $y_{\rm ref} \colon t\mapsto \cos t$.
We emphasize that the {\it function} $y_{\rm ref}(\cdot)$ is not available {\em a priori} to the controller: all that is available is the function value
 at the current time~$t$ together with the values of its first~$\hat r-1$ derivatives, $y^{(i)}_{\rm ref}(t)$, $i=0,\ldots,\hat r-1$.
We consider two cases.

\vspace{2mm}
\noindent
{\bf Case 1:} If $0<\vartheta<\frac{\pi}{2}$, then system~\eqref{mass.on.car.equ} has relative degree~$r=2$,
and the funnel controller~\eqref{eq:FC} with~$\hat r = r = 2$ is
\[
u(t) = -\alpha\big(w(t)^2\big)\, w(t),~~\text{with}~~  w(t) = \varphi(t) \dot e(t) + \alpha\big(\varphi(t)^2 e(t)^2\big)\, \varphi(t) e(t),
\]
where~$\alpha(s) = 1/(1-s)$ for~$s\in[0,1)$. The controller presented in~\cite{BergLe18}  takes the form
\begin{equation}\label{eq:FC-BLR-2}
  u(t) = -\alpha\big(\varphi_1(t)^2 w_1(t)^2\big)\, w_1(t),~~\text{with}~~w_1(t) = \dot e(t) + \alpha\big(\varphi(t)^2 e(t)^2\big)\, e(t),
  \end{equation}
where~$\varphi_1$ is a second funnel function, chosen appropriately, cf.~\cite{BergLe18}. Note that $w(t)=\varphi(t) w_1(t)$. As {simulations show}, the performance of the controller~\eqref{eq:FC-BLR-2} can be improved compared to the simulations in~\cite{BergLe18}, by choosing $\varphi_1 = \varphi$. As in~\cite{BergLe18}, we  set~$\varphi(t) = (5 e^{-2t} + 0.1)^{-1}$ for $t\geq 0$.
\captionsetup[subfloat]{labelformat=empty}
\begin{figure}[h!tb]
  \centering
  \subfloat[Fig.~\ref{fig:sim-rd2}a: Funnel and tracking errors]
{
\centering
\hspace{-2mm}
  \includegraphics[width=0.53\textwidth]{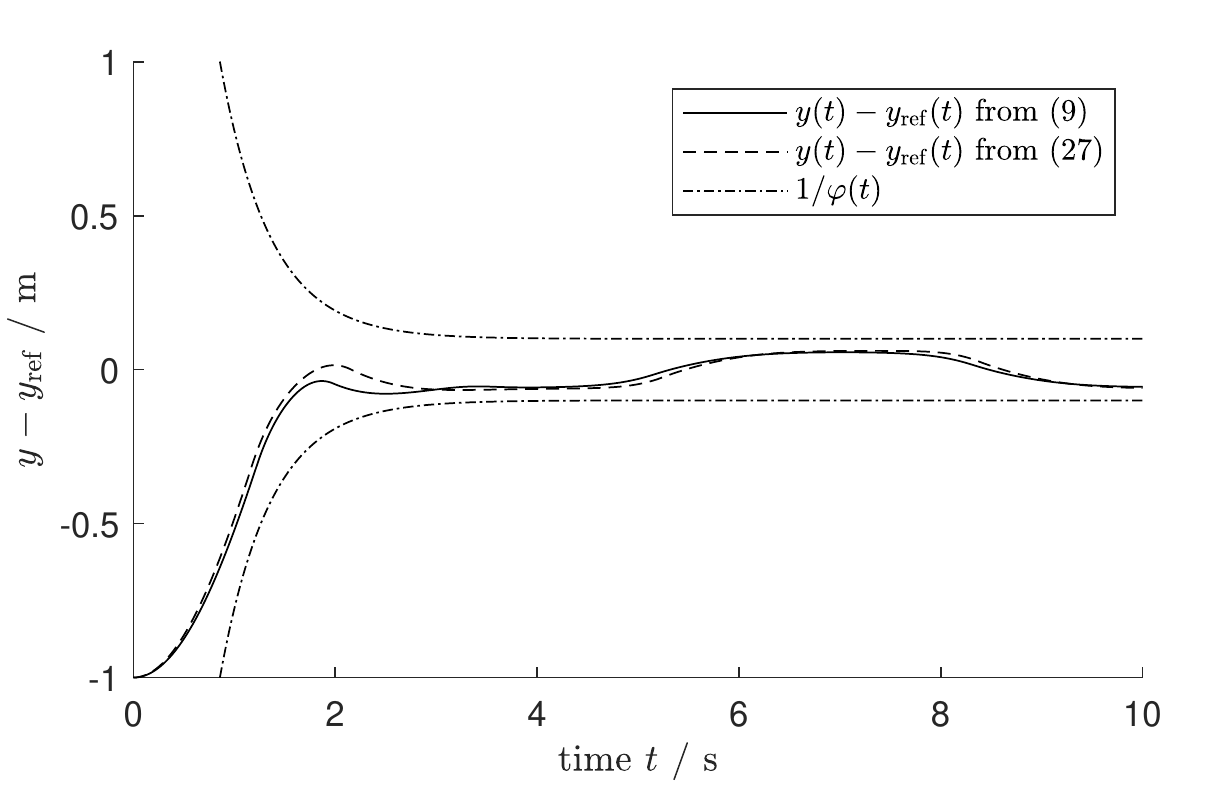}
\label{fig:sim-e}
}
\subfloat[Fig.~\ref{fig:sim-rd2}b: Input functions]
{
\centering
\hspace{-5mm}
  \includegraphics[width=0.53\textwidth]{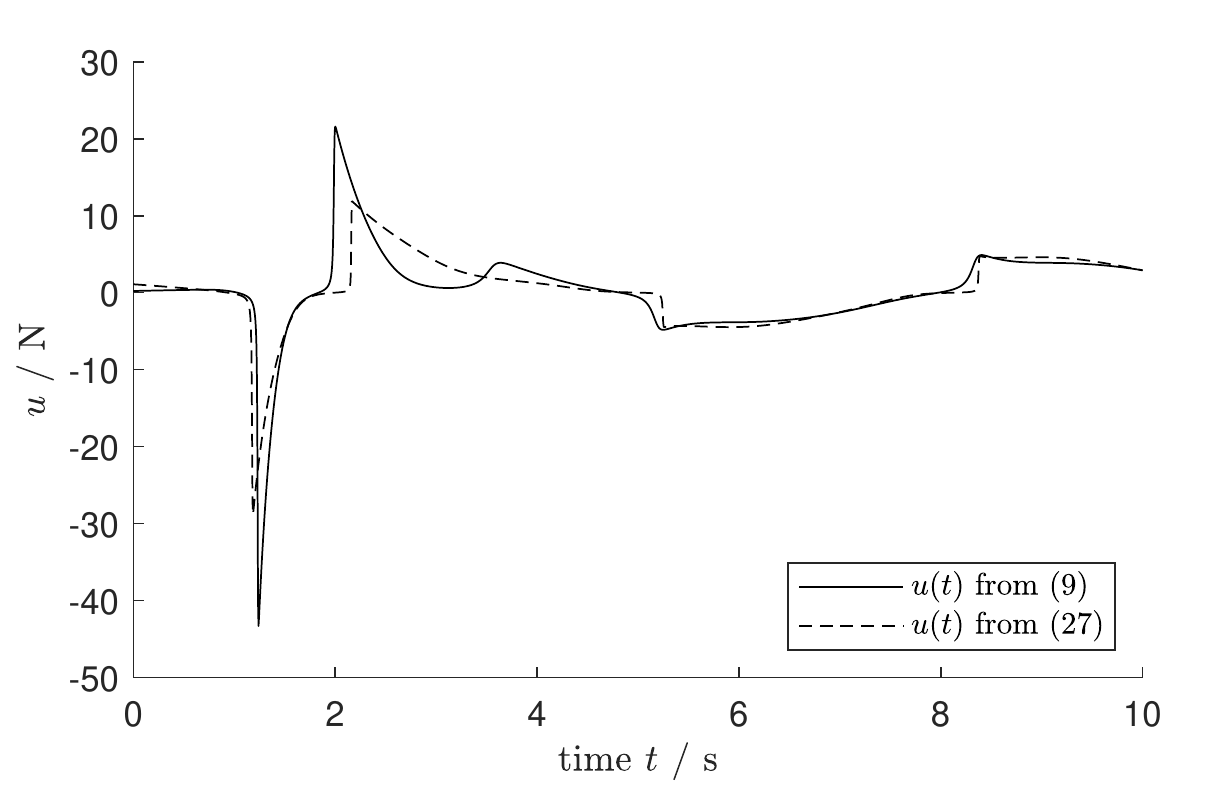}
\label{fig:sim-u}
}
\caption{Simulation, under controllers~\eqref{eq:FC} and~\eqref{eq:FC-BLR-2}, of system~\eqref{mass.on.car.equ} with $\vartheta=\frac{\pi}{4}$.}
\label{fig:sim-rd2}
\end{figure}
\\\noindent
The performance of the controllers~\eqref{eq:FC} and~\eqref{eq:FC-BLR-2} applied to~\eqref{mass.on.car.equ} is depicted in Fig.~\ref{fig:sim-rd2}.
Fig.~\ref{fig:sim-e} shows the tracking errors generated by the two different controllers, while Fig.~\ref{fig:sim-u} shows the respective input functions.
Comparable performance is evident, suggesting broadly similar efficacy in cases wherein both controllers are feasible.  However,
~\eqref{eq:FC} is feasible in certain situations which are outside the scope of ~\eqref{eq:FC-BLR-2}.  For example,~\eqref{eq:FC} is able to
achieve asymptotic tracking, {to address the issue of an unknown control direction} and is applicable when the instantaneous
value~$\dot y_{\rm ref}(t)$ is not available to the controller: these features form the basis of {the example in Section~\ref{Ssec:Ex-deadzone}} below.
\\[2mm]
\noindent
{\bf Case 2:} If $\vartheta=0$ and $d\ne 0$, then system~\eqref{mass.on.car.equ} has relative degree~$r=3$. Then the funnel controller~\eqref{eq:FC},  with~$\hat r = r = 3$, takes the form
\begin{align*}
  w(t) &= \varphi(t) \ddot e(t) + {\gamma}\big(\varphi(t) \dot e(t) + {\gamma}\big(\varphi(t) e(t)\big)\big),\\
  u(t) &= -{\gamma}\big(w(t)\big),
\end{align*}
where~$\gamma(s) = s\alpha(s^2)$ for~$s\in(-1,1)$. The controller presented in~\cite{BergLe18} reads
\begin{equation}\label{eq:FC-BLR-3}
\begin{aligned}
  w_1(t) &= \dot e(t) + {\alpha}\big(\varphi(t)^2 e(t)^2\big)\, e(t),\\
  w_2(t) &= \dot w_1(t) + {\alpha}\big(\varphi_1(t)^2 w_1(t)^2\big)\, w_1(t)\\
  &= \ddot e(t) + 2 {\alpha}\big(\varphi(t)^2 e(t)^2\big)^2 \big( \dot \varphi(t) \varphi(t) \|e(t)\|^2 + \varphi(t)^2 e(t)^\top \dot e(t)\big)\, e(t)\\
   &\quad + {\alpha}\big(\varphi(t)^2 e(t)^2\big)\, \dot e(t) + {\alpha}\big(\varphi_1(t)^2 w_1(t)^2\big)\, w_1(t),\\
  u(t) &= -{\alpha}\big(\varphi_2(t)^2 w_2(t)^2\big)\, w_2(t),
\end{aligned}
\end{equation}
where~$\varphi_1, \varphi_2$ are appropriate additional funnel functions, cf.~\cite{BergLe18}.
Here, we choose $\varphi_1 = \varphi_2=\varphi$, with $\varphi(t) = (3 e^{-t} + 0.1)^{-1}$ for $t\ge 0$ and
compare the controller~\eqref{eq:FC} with~\eqref{eq:FC-BLR-3}.

\captionsetup[subfloat]{labelformat=empty}
\begin{figure}[h!tb]
  \centering
  \subfloat[Fig.~\ref{fig:sim-rd3}a: Funnel and tracking errors]
{
\centering
\hspace{-2mm}
  \includegraphics[width=0.53\textwidth]{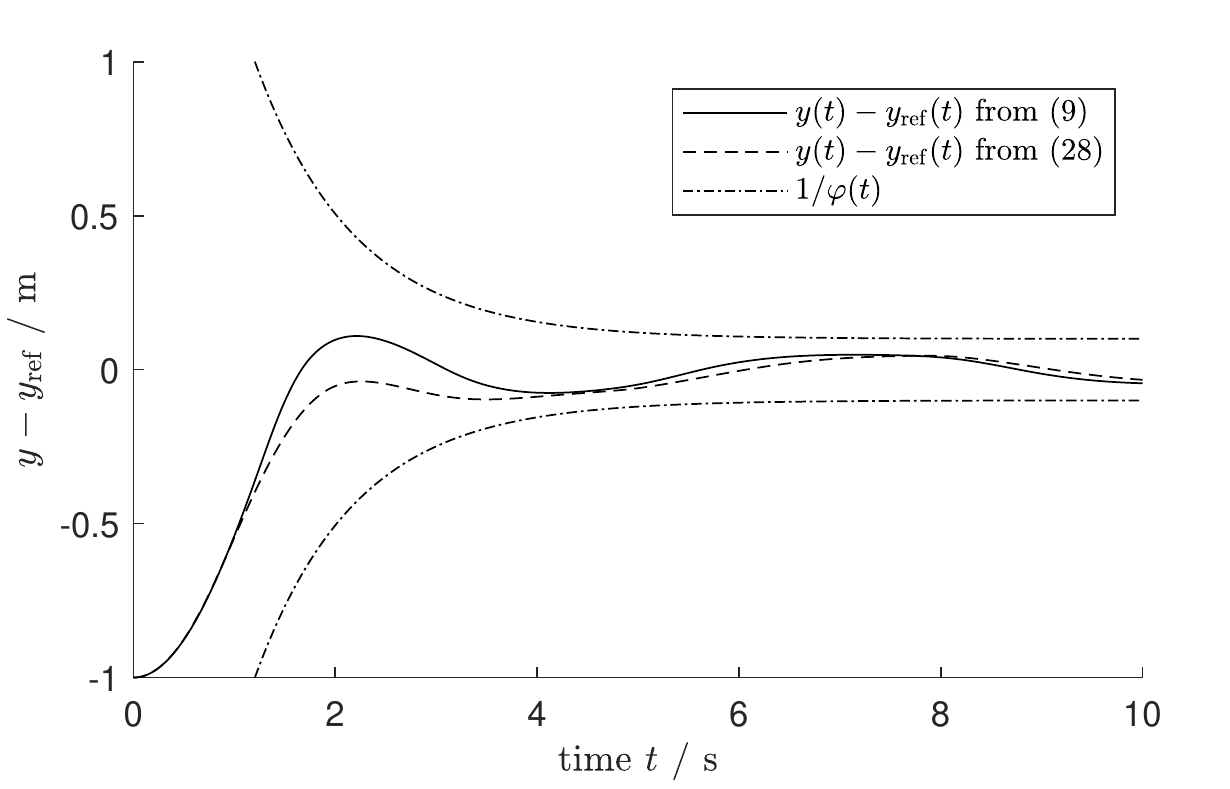}
\label{fig:sim2-e}
}
\subfloat[Fig.~\ref{fig:sim-rd3}b: Input functions]
{
\centering
\hspace{-5mm}
  \includegraphics[width=0.53\textwidth]{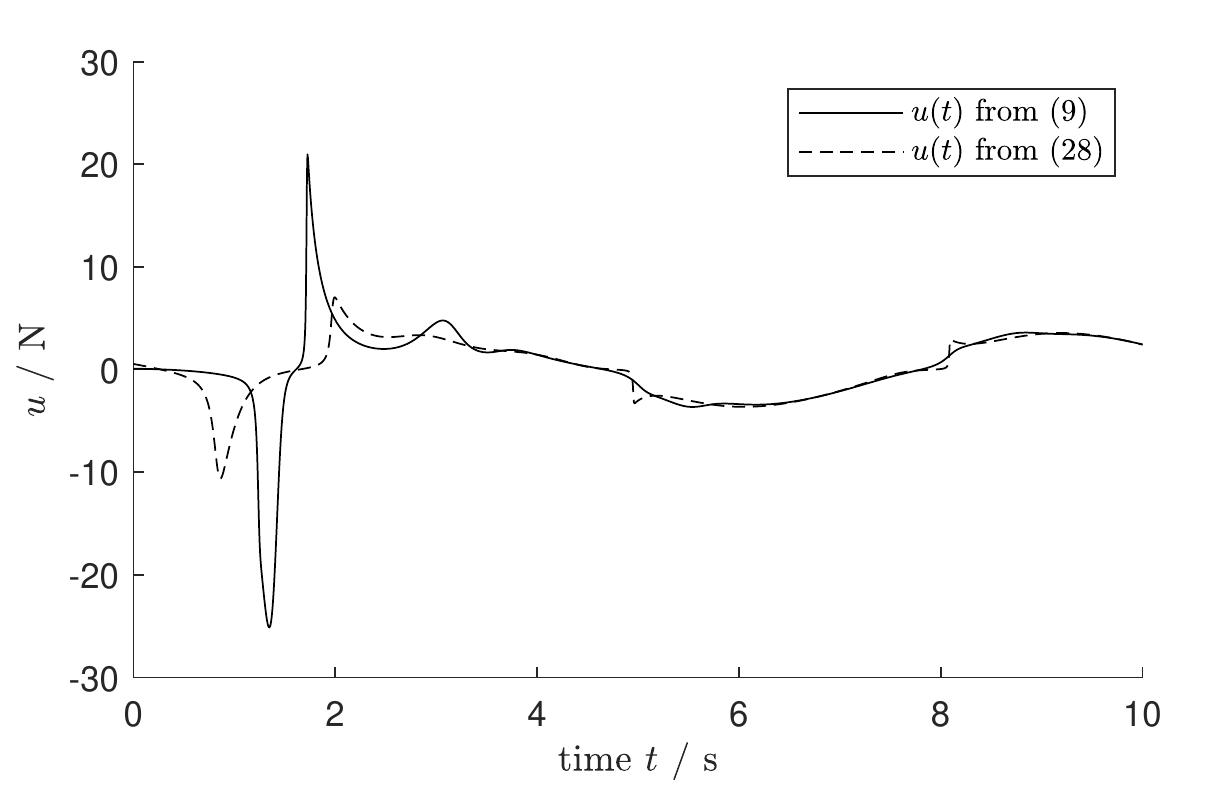}
\label{fig:sim2-u}
}
\caption{Simulation, under controllers~\eqref{eq:FC} and~\eqref{eq:FC-BLR-3}, of system~\eqref{mass.on.car.equ} with $\vartheta=0$.}
\label{fig:sim-rd3}
\end{figure}

\noindent The simulation suggests that the controllers are broadly similar in performance. While
controller~\eqref{eq:FC} requires more input action than controller~\eqref{eq:FC-BLR-3},
the latter exhibits a significantly higher level of complexity, which makes it more difficult to implement (this issue becomes
even more severe for relative degrees higher than three).

\subsection{{\bf Nonlinear MIMO system}}\label{Ssec:NonlMimo}

\begin{figure}[h!tb]
%
%
%
\begin{center}
\includegraphics[width=6cm]{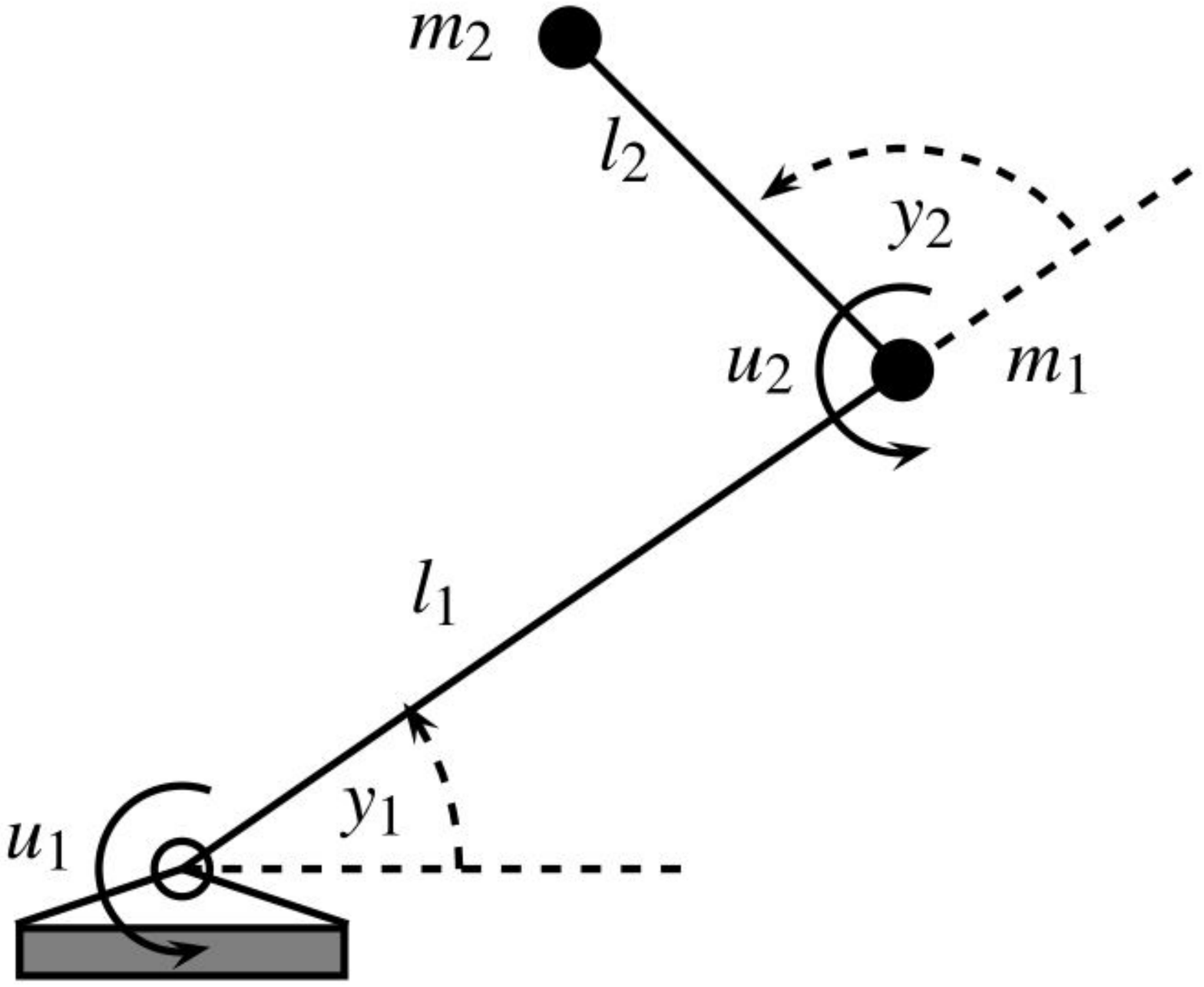}
\end{center}
\vspace{-0.5cm}
\caption{Planar rigid revolute joint robotic manipulator.}
\label{Fig:robot-model}
\end{figure}

As a nonlinear multi-input, multi-output example we consider the robotic manipulator from~\cite[Ch.~13]{Hack17} as depicted in Fig.~\ref{Fig:robot-model}. It is planar, rigid, with revolute joints and has two degrees of freedom.

The two joints are actuated by~$u_1$ and~$u_2$ (in \si{\newton\metre}). The links are assumed to be massless and have lengths~$l_1$ and~$l_2$ (in \si{\metre}), resp., with point masses~$m_1$ and~$m_2$ (in \si{\kilo\gram})  attached to their ends. The two outputs are the joint angles~$y_1$ and~$y_2$ (in \si{\radian}) and the equations of motion are given by (see also~\cite[p.~259]{SponHutc06})
\begin{equation}\label{eq:robot}
M(y(t))\ddot{y}(t)+C(y(t),\dot{y}(t))\dot{y}(t)+G(y(t))=u(t)
\end{equation}
with initial value $(y(0),\dot{y}(0))=\left(0,0\right)$, inertia matrix
\[
M:\R^2\to \R^{2\times 2}, \ (y_1,y_2)\mapsto
\begin{bmatrix}
m_1l_1^2+m_2(l_1^2+l_2^2+2l_1l_2\cos(y_2)) & m_2(l_2^2+l_1l_2\cos(y_2))\\
m_2(l_2^2+l_1l_2\cos(y_2)) & m_2l_2^2
\end{bmatrix}
\]
centrifugal and Coriolis force matrix
\[
C:\R^4\to \R^{2 \times 2},\ (y_1,y_2,v_1,v_2)\mapsto \begin{bmatrix}
-2m_2l_1l_2\sin(y_2)v_1 & -m_2l_1l_2\sin(y_2)v_2\\
-m_2l_1l_2\sin(y_2)v_1 & 0
\end{bmatrix},
\]
and gravity vector
\[
G:\R^2\to \R^2,\ (y_1,y_2)\mapsto g \begin{pmatrix}
m_1l_1\cos(y_1)+m_2(l_1\cos(y_1)+l_2\cos(y_1+y_2)) \\
m_2l_2\cos(y_1+y_2)
\end{pmatrix},
\]
where $g=\SI{9.81}{\metre\per{\second^2}}$ is the acceleration of gravity. Multiplying~\eqref{eq:robot} with~$M(y(t))^{-1}$, which is pointwise positive definite, from the left we see that the resulting system is of the form~\eqref{eq:nonlSys} and satisfies the positive-definite high-gain property, hence it belongs to~$\cN^{2,2}$.

For the simulation, we choose the parameters~$m_1=m_2=1$, $l_1=l_2=1$ and the reference signal $y_{\rm ref}\colon t\mapsto (\sin t,\,\sin 2t)$. We compare the controller~\eqref{eq:FC} to the multivariate version of~\eqref{eq:FC-BLR-2} from~\cite{BergLe18}, that is
\begin{equation}\label{eq:FC-BLR-MIMO}
u(t) = -\alpha\big(\varphi_1(t)^2 \|w_1(t)\|^2\big)\, w_1(t),\quad\text{with}\quad
  w_1(t) = \dot e(t) + \alpha\big(\varphi(t)^2 \|e(t)\|^2\big)\, e(t),
\end{equation}
where $\alpha(s) = 1/(1-s)$ for~$s\in[0,1)$. We choose $\varphi(t) = (4e^{-2t}+0.1)^{-1} = \varphi_1(t)$ for $t\ge 0$.

\begin{figure}[h!tb]
  \centering
  \subfloat[Fig.~\ref{fig:robot}a: Funnel and first tracking error components]
{
\centering
\hspace{-2mm}
  \includegraphics[width=0.53\textwidth]{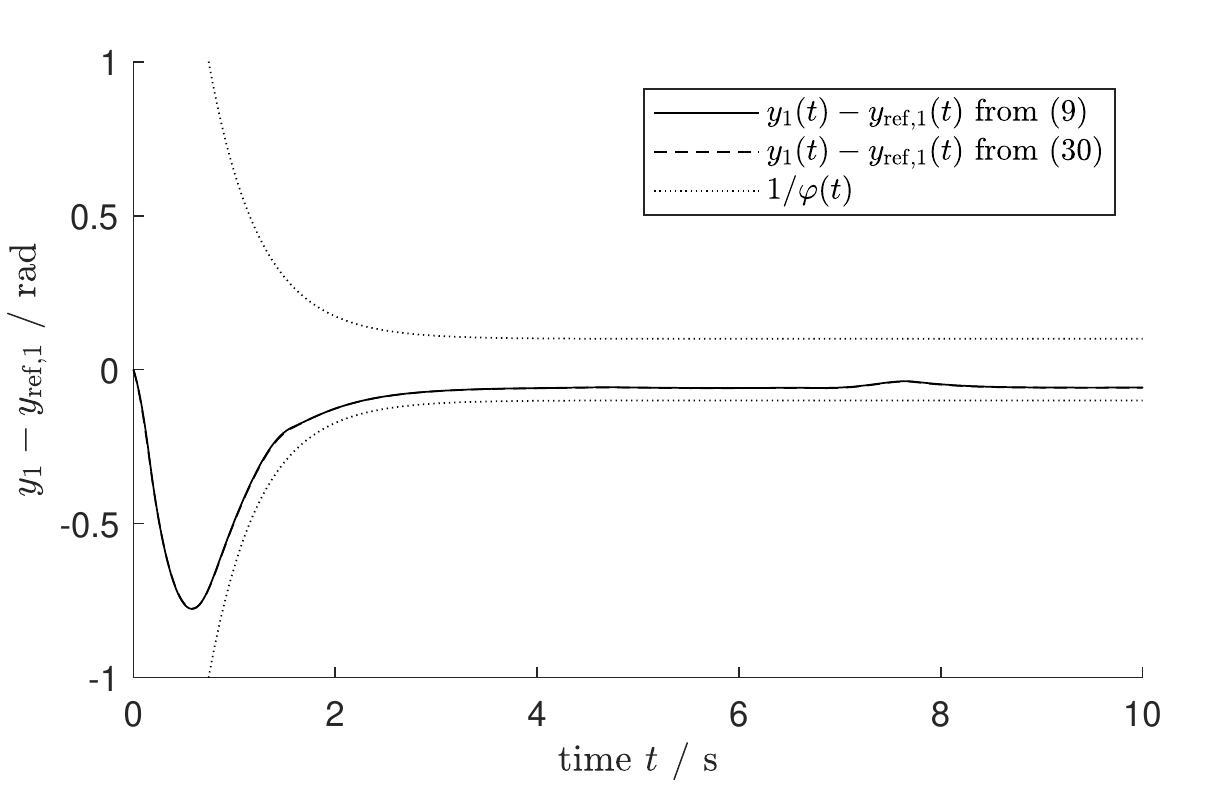}
\label{fig:robot-error1}
}
\subfloat[Fig.~\ref{fig:robot}b: Funnel and second tracking error com-\newline ponents]
{
\centering
\hspace{-5mm}
  \includegraphics[width=0.53\textwidth]{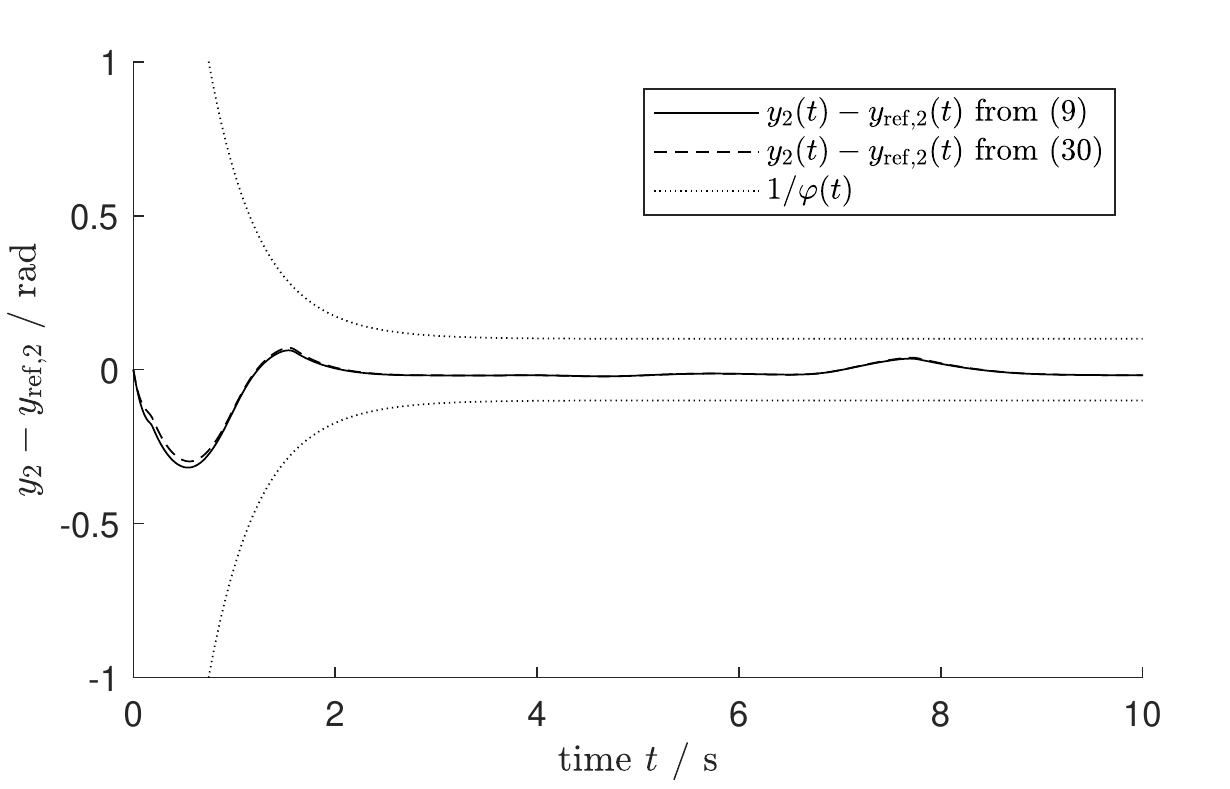}
\label{fig:robot-error2}
}\\
\subfloat[Fig.~\ref{fig:robot}c: First input components]
{
\centering
\hspace{-2mm}
  \includegraphics[width=0.53\textwidth]{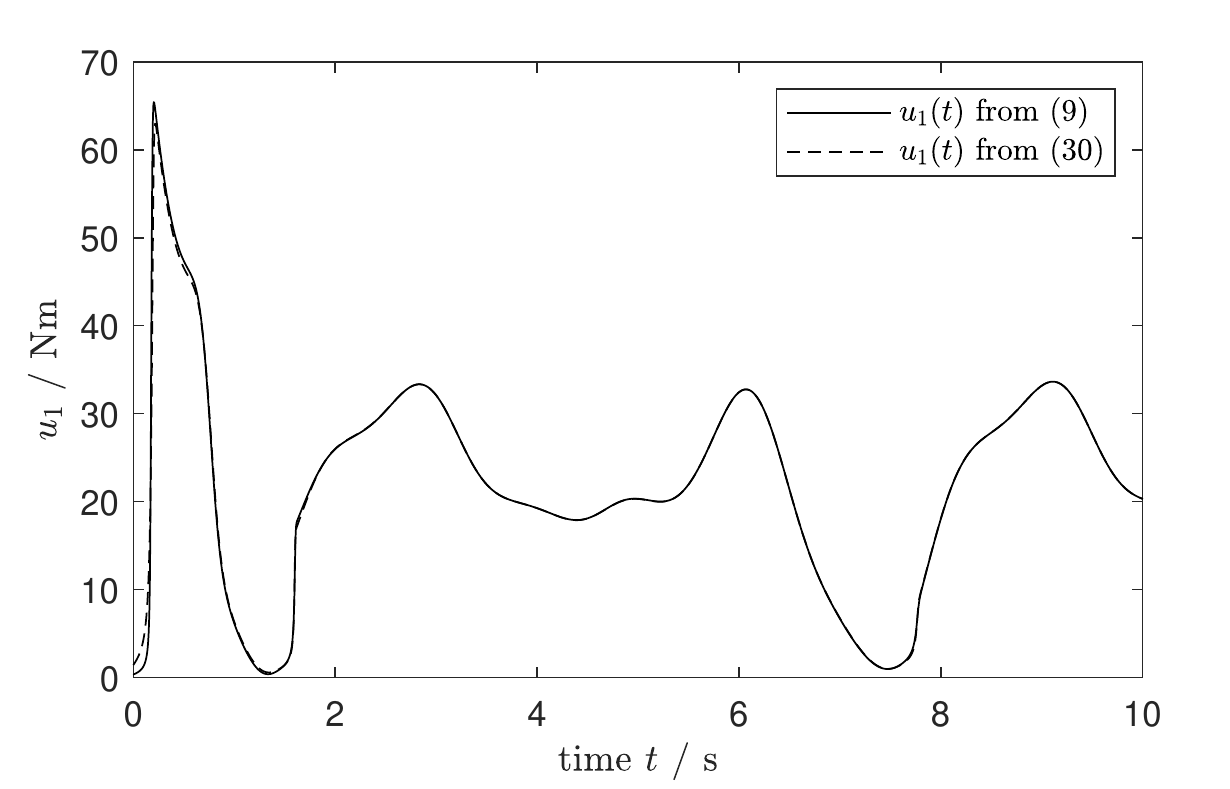}
\label{fig:robot-input1}
}
\subfloat[Fig.~\ref{fig:robot}d: Second input components]
{
\centering
\hspace{-5mm}
  \includegraphics[width=0.53\textwidth]{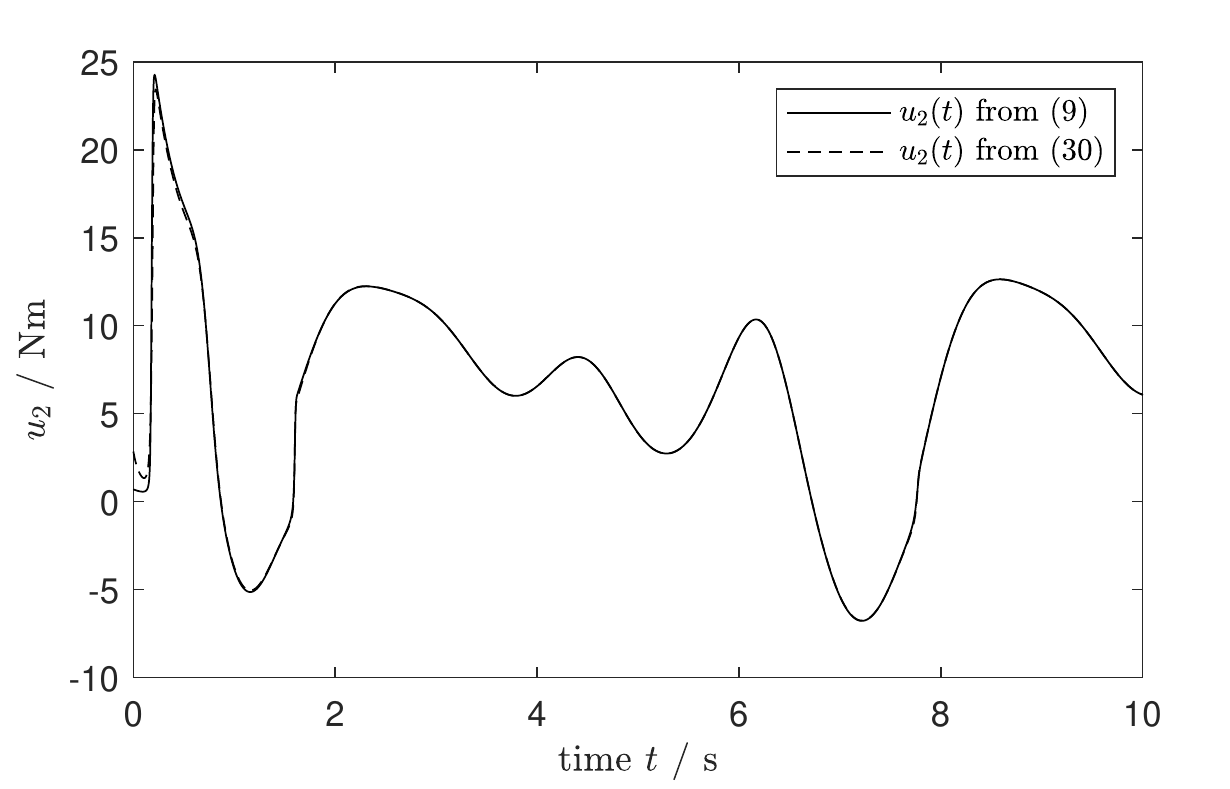}
\label{fig:robot-input2}
}
\caption{Simulation of the controllers~\eqref{eq:FC} and~\eqref{eq:FC-BLR-MIMO} applied to~\eqref{eq:robot}.}
\label{fig:robot}
\end{figure}

\noindent
The simulation of the controllers~\eqref{eq:FC} and~\eqref{eq:FC-BLR-MIMO} applied to~\eqref{eq:robot} over the time interval~$[0,10]$ is depicted in Fig.~\ref{fig:robot}. It can be seen that for this example both controllers exhibit a nearly identical performance.

\begin{remark}
  A closer look at the simulations reveals that the controller performance of~\eqref{eq:FC} differs from that of the controller presented in~\cite{BergLe18} for the example in Subsection~\ref{Ssec:Mass_on_car}, while it is practically identical for the example in Subsection~\ref{Ssec:NonlMimo}. Since the different dimensions of input/output spaces ($m=1$ compared to~$m=2$) is probably not the reason here, the presumable cause seems to be the internal dynamics. System~\eqref{mass.on.car.equ} has two-dimensional internal dynamics in Case~1 ($r=2$) and one-dimensional internal dynamics in Case~2 ($r=3$), while system~\eqref{eq:robot} has trivial internal dynamics. This seems to suggest that the controllers exhibit a different behaviour in the presence of non-trivial internal dynamics.
\end{remark}

\subsection{{\bf {{A nonlinear system with dead-zone input}}}}\label{Ssec:Ex-deadzone}

{{To demonstrate that the controller~\eqref{eq:FC} can achieve asymptotic tracking and is feasible when the control direction is unknown,
we {treat} a system with dead-zone input and also investigate the case {wherein} $\dot y_{\rm ref}(t)$ is not available for feedback. }}
{Specifically, we consider
\begin{equation}\label{eq:ex22rev}
\left.
\begin{aligned}
&\dot\xi_1(t) =\big(1+\xi_1(t)^2\big)\xi_2(t),\\
&\dot \xi_2(t) = \alpha_1\xi_1(t) +\alpha_2\xi_2(t) + \alpha_3\eta(t) +\beta\big(u(t)\big),
\\
&\dot\eta (t) =-\eta(t)^2\big( \alpha_4 \xi_1(t)+\alpha_5 \xi_2(t)+\eta(t)\big),\\
&(\xi_1(0),\xi_2(0),\eta(0))=(\xi_1^0,\xi_2^0,\eta^0)\in\R^3
\end{aligned}
\right\}~~\text{with output $y(t)=\xi_1(t),$}
\end{equation}
with real parameters $\alpha_i\in\R$, $i=1,\ldots,5$, and with a dead-zone input function} $\beta$ as in Section~\ref{Ssec:DeadZone}. {{We show that}} system~\eqref{eq:ex22rev} belongs to the class of systems $\cN^{1,2}$. {
Introducing the function
\[
g\colon \R^2\times\R\to \R,~(x,z)=(x_1,x_2,z)\mapsto -z^2\big(\alpha_4x_1+\alpha_5 (1+x_1^2)^{-1}x_2+z\big),
\]
and writing $y_1(t)=\xi_1(t)=y(t)$ and $y_2(t)=\big(1+\xi_1(t)^2\big)\xi_2(t)=\dot y(t)$, the third equation in~\eqref{eq:ex22rev} may be expressed in the form
\begin{equation}\label{ivp-id}
\dot\eta (t)=g\big(y_1(t),y_2(t),\eta (t)\big),~~\eta (0)=\eta^0
\end{equation}
which, viewed in isolation as a system with independent inputs $(y_1,y_2)$, generates a controlled flow~$\Xi$.  In particular,
for $\eta^0\in\R$ and $(y_1,y_2)\in L^\infty_{\text{\rm{loc}}}(\R_{\ge 0},\R^2)$,
the initial-value problem \eqref{ivp-id}
has unique maximal solution $\eta(\cdot) = \Xi(\cdot\,;y_1,y_2,\eta^0):[0,\omega)\to\R$, $0 <\omega \leq\infty$.
Also, writing $\alpha :=|\alpha_4|+|\alpha_5|$, we have
\[
\forall\, (x,z)\in\R^2\times\R:\ z\, g(x,z)\leq -z^4+ \alpha |z|^3\|x\|\leq -\quarter z^4 +\quarter (\alpha \|x\|)^4,
\]
wherein Young's inequality has been used.
Therefore,~$V\colon z\mapsto \half z^2$  is an ISS-Lyapunov function for system~\eqref{ivp-id} which, in consequence, is input-to-state stable, see~\cite[Rem.~2.4 \& Lem.~2.14]{SontWang95b}.
Therefore, for all $c_0 >0$, there exists $c_1 >0$ such that for all $\eta^0\in\R$ and all $(y_1,y_2)\in L^\infty_{\text{loc}}(\R_{\ge 0},\R^{2})$ we have
\[
 \|\eta^0\|+\esup_{t\geq 0}\|(y_1(t),y_2(t))\|  \leq c_0 \quad \implies\quad
    \sup\nolimits_{t\geq 0}\, \|\Xi (t;y_1,y_2,\eta^0)\| \leq c_1.
\]
The above property ensures that solutions of \eqref{ivp-id} are globally defined: specifically, for each $\eta^0\in\R$ and $(y_1,y_2)\in L^\infty_{\text{loc}}(\R_{\ge 0},\R^{2})$,
the unique maximal solution of~\eqref{ivp-id} has interval of existence $[0,\infty)$. Therefore, with each fixed $\eta^0\in\R$, we may associate an operator
\[
\fT\colon C(\R_{\ge 0},\R^{2}) \to L^\infty_{\text{loc}}(\R_{\ge 0},\R^3),\quad
(y_1,y_2)\mapsto \big( y_1,y_2, \, \Xi(\cdot\,;y_1,y_2,\eta^0)  \big).
\]
Clearly, $\fT$ is causal, i.e.,~(TP1) holds.  Moreover, the above ISS property of~$\Xi$ ensures that  properties~(TP2) and~(TP3) hold. Therefore, $\fT\in \TTT{2}{3}{0}$. Defining
\[
    f:\R^4\to\R,\ (x_1,x_2,z,u)\mapsto \frac{2x_1 x_2^2}{1+x_1^2}  +\alpha_2x_2 + (1+x_1^2) \big(\alpha_1 x_1 + \alpha_3 z + \beta(u)\big),
\]
it is readily verified that~\eqref{eq:ex22rev} is equivalent to
\[
    \ddot y(t) = f\big(\fT(y,\dot y)(t), u(t)\big), ~~y(0)=\xi_1^0,~,\dot y(0)=\big(1+(\xi_1^0)^2\big)\xi_2^0.
\]
Furthermore, by applying the findings of Sections~\ref{Ssec:InputNonl} and~\ref{Ssec:DeadZone} we have that~$f$ satisfies the high-gain property and hence $(0,f,\fT)\in\cN^{1,2}$. Therefore, feasibility of funnel control follows from Theorem~\ref{Thm:FunCon-Nonl}.}

{For the simulation, we (randomly) select
\[
\alpha_1=\alpha_3=\alpha_5=1~~\text{and}~~\alpha_2=-\alpha_4=-2,
\]
and the dead-zone input function as
\[
   \beta : \R\to \R, \quad v\mapsto \left\{
    \begin{array}{rcl} v-1, && v\ge 1,\\ 0,&& -1<v<1,\\ v+1, && v\le -1.\end{array}\right.
\]
}
The initial values are chosen as $\xi_1(0)= \xi_2(0)=\eta(0)=0$ and the reference signal is $y_{\rm ref}\colon t\mapsto \cos t$. For the funnel controller~\eqref{eq:FC} we choose the design parameters $\alpha\colon s\mapsto 1/(1-s)$ and $N\colon s\mapsto s \sin s$; the latter choice is based on the assumption that the exact shape of~$\beta$ (and in particular the control direction) is unknown to the controller.

We consider two different cases: If information of the instantaneous signals~$\dot y_{\rm ref}(t)$ are available to the controller, then we choose $\hat r = 2 = r$ and an unbounded funnel function $\varphi\colon t\mapsto t^2$. If information of~$\dot y_{\rm ref}(t)$ is not available, then we choose $\hat r = 1 < 2$ and a bounded funnel function $\varphi\colon t\mapsto \big(2 e^{-t} + 0.01\big)^{-1}$.

\captionsetup[subfloat]{labelformat=empty}
\begin{figure}[h!tb]
  \centering
  \subfloat[Fig.~\ref{fig:sim4}a: Funnel and tracking error for $\hat r=1$]
{
\centering
\hspace{-2mm}
  \includegraphics[width=0.52\textwidth]{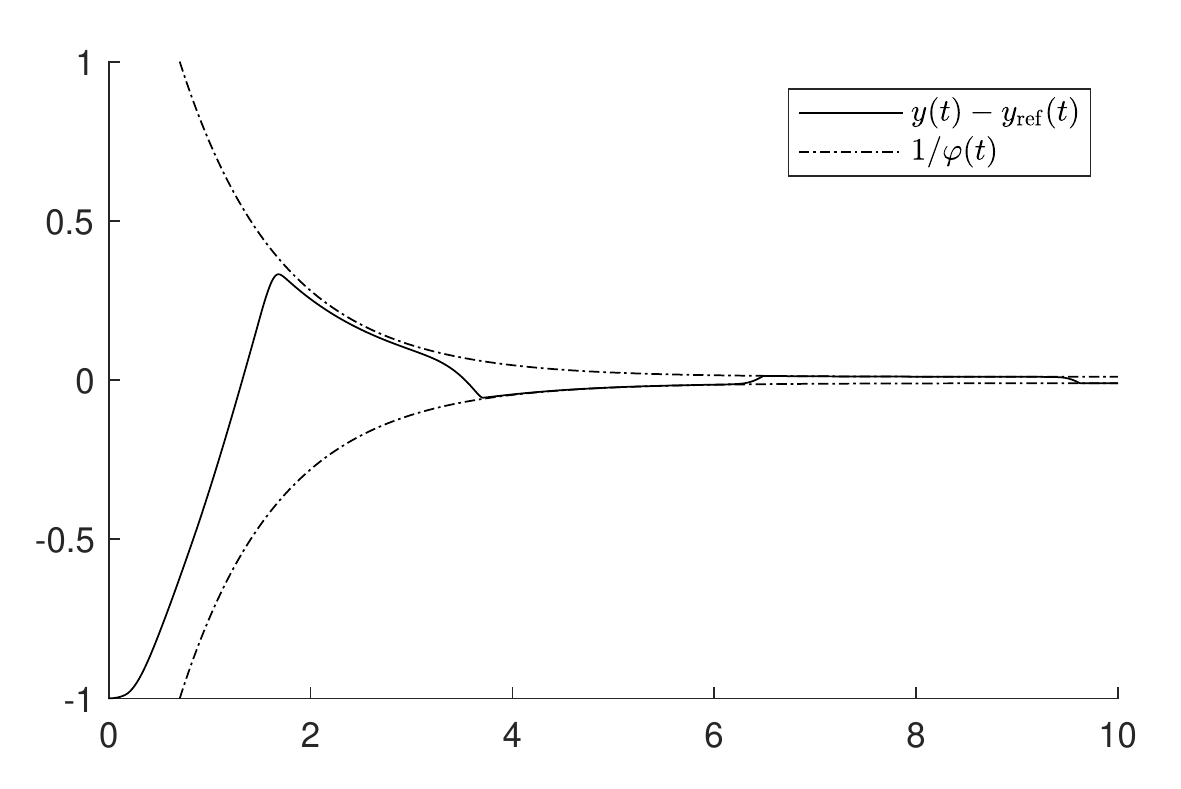}
\label{fig:sim4-e1}
}
\subfloat[Fig.~\ref{fig:sim4}b: Funnel and tracking error for $\hat r=2$]
{
\centering
\hspace{-4mm}
  \includegraphics[width=0.52\textwidth]{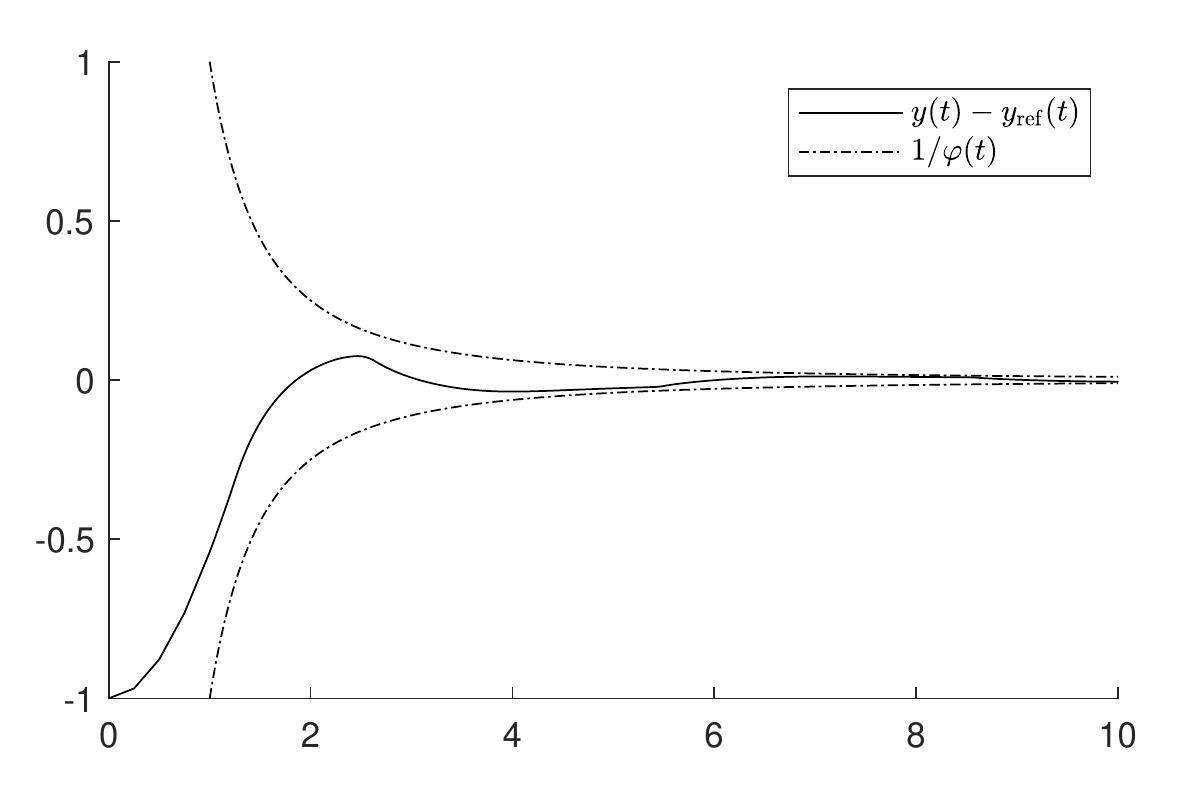}
\label{fig:sim4-e2}
}\\
\subfloat[Fig.~\ref{fig:sim4}c: Input function for $\hat r=1$]
{
\centering
\hspace{-2mm}
  \includegraphics[width=0.52\textwidth]{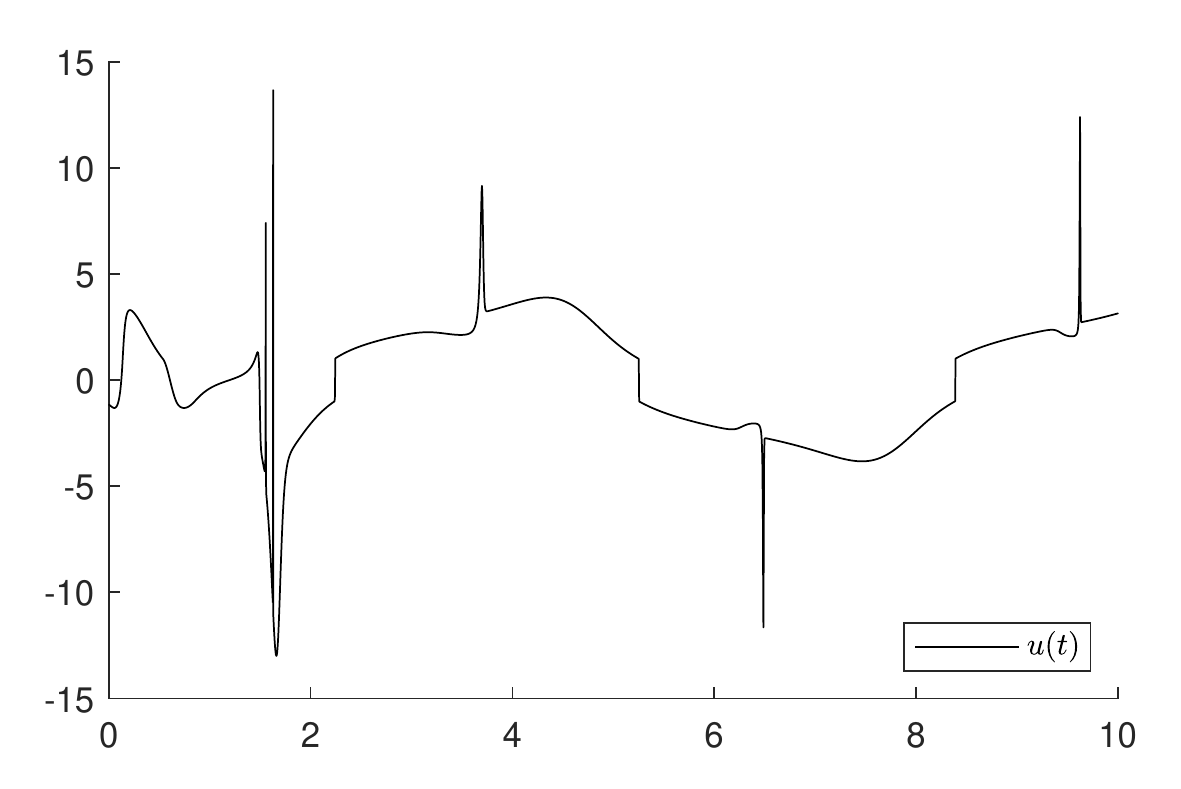}
\label{fig:sim4-u1}
}
\subfloat[Fig.~\ref{fig:sim4}d: Input function for $\hat r=2$]
{
\centering
\hspace{-4mm}
  \includegraphics[width=0.52\textwidth]{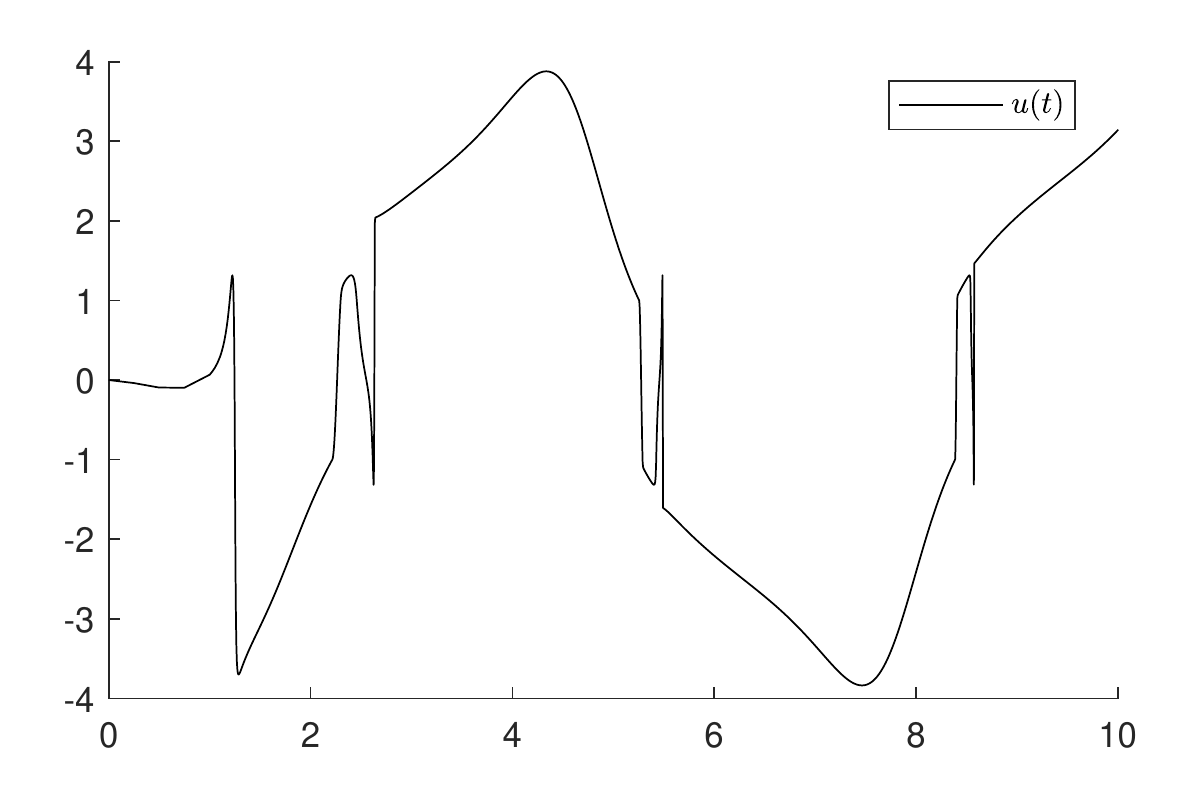}
\label{fig:sim4-u2}
}
\caption{Simulation of system~\eqref{eq:ex22rev} under control~\eqref{eq:FC} in the cases $\hat r=1$ and $\hat r = 2$.}
\label{fig:sim4}
\end{figure}

\noindent
The simulation of the controller~\eqref{eq:FC} applied to~\eqref{eq:ex22rev} in the cases $\hat r=1$ and $\hat r = 2$ is depicted in Fig.~\ref{fig:sim4}. The ``jumps'' in the input~$u$ are due to the dead-zone induced by the function~$\beta$. Comparing Figs.~\ref{fig:sim4}d and~\ref{fig:sim4}c a degradation in performance may be observed.  However, this is not surprising in view of the enhanced information available for feedback in
case $\hat r=2$. We may also observe, that in the latter case asymptotic tracking is achieved.

\section{Conclusion}\label{Sec:Concl}

{ An asymptotic and non-asymptotic tracking control
objective has been achieved
for a large class of nonlinear systems with ``higher relative degree''
described by functional differential equations that satisfy a high-gain property.
A feedback strategy has been developed which is} simple in the sense
of funnel control and as ``simple'' as one may expect  for higher relative degree.
We believe that the present { contribution is somehow ``definitive'' in the context of
 funnel control for} nonlinear systems whose internal dynamics satisfy a BIBO property { (viewed as a generalization of the minimum phase condition for linear systems.)} First results on funnel control for systems which are not minimum phase are given in~\cite{Berg20} for uncertain linear systems and in~\cite{BergLanz20} for a nonlinear robotic manipulator.

{In the present paper we did not treat} funnel control for systems described by partial differential
equations. {This is however, a very important field and in fact  very different.}
 On the one hand, there are systems which have a well-defined relative degree and exhibit infinite-dimensional internal dynamics, see e.g.~\cite{BergPuch20a}. Such systems are susceptible to funnel control with the control laws presented in the present paper; for instance, a linearized model of a moving water tank, where sloshing effects appear, is discussed in~\cite{BergPuch19}. On the other hand, not even every linear infinite-dimensional system has a well-defined relative degree, in which case the results presented here cannot be applied. For such systems,  the feasibility of funnel control has to be investigated directly for the (nonlinear) closed-loop system, see e.g.~\cite{ReisSeli15b} for a boundary controlled heat equation,~\cite{PuchReis19pp} for a general class of boundary control systems,~\cite{BergBrei19} for the monodomain equations (which represents defibrillation processes of the human heart) {and~\cite{Berg20pp} for the Fokker-Planck equation corresponding to the Ornstein-Uhlenbeck process.}

One {important} {problem remains}:  non-derivative funnel control, {that is, when}
only  the output~$y$ is available  for feedback,
 but not its first~$r-1$ derivatives~$\dot y,\ldots,y^{(r-1)}$.
First results on this have been {obtained} in~\cite{IlchRyan06b,IlchRyan07}
using a backstepping approach. However, these results {necessitate a level of controller complexity which, on the evidence of numerical simulation,
can lead to practical performance drawbacks.}
An attempt to overcome these {backstepping-induced drawbacks  through the adoption of pre-compensators can be found in
~\cite{BergReis18a,BergReis18b}
but only for systems with relative degree at most three: the higher relative degree case remains open, even in the context of
single-input, single-output linear systems with positive high-frequency gain and asymptotically stable zero dynamics.}

\begin{appendices}

\section{Proofs}\label{Sec:App}


\begin{proof}[\bf Proof of Theorem \ref{Thm:FunCon-Nonl}]
%
%
For $k=1,\ldots,r$, we define
\begin{align*}
&\pi_k:\R_{\ge 0}\times\R^{rm}\to\R^{km},\\
&(t,\xi)=(t,\xi_1,\ldots,\xi_r)\mapsto \left\{\begin{array}{l}
\varphi(t)  \Big( \xi_1-y_{\text{\rm ref}}(t),\ldots, \xi_k-y_{\text{\rm ref}}^{(k-1)}(t)\Big),\quad k=1,\ldots,\hat r,
\\[1ex]
\varphi (t)\Big(  \xi_1-y_{\text{\rm ref}}(t),\ldots, \xi_{\hat r}-y_{\text{\rm ref}}^{(\hat r-1)}(t),  \xi_{\hat r+1},\ldots ,\xi_k\Big),
\\
\hfill
k=\hat r+1,\ldots ,r.
\end{array}\right.
\end{align*}
The proof now proceeds in several steps.
\\[1ex]
{\it Step~1.}\
We recast the feedback-controlled system in the form of an initial-value problem to which a variant of an
extant existence theory applies. Set~$n=rm$ and
\[
\cD:=\setdef{(t,\xi)\in\R_{\ge 0}\times\R^n }{ \pi_r (t,\xi)\in\cD_r },
\]
which is non-empty and relatively open, and define $\rho\colon\cD\to\cB$ by $\rho := \rho_r\circ\pi_r$.   Introducing the function $F\colon \cD\times\R^q\to\R^n$ given by
\[
(t,\xi,\eta)=(t,\xi_1,\ldots,\xi_r,\eta)\mapsto F(t,\xi,\eta):=\begin{pmatrix} \xi_2\\\vdots\\\xi_r\\f\big(d(t),\eta,(N\circ \alpha)(\|\rho(t,\xi)\|^2)\,\rho(t,\xi)\big)\end{pmatrix}
\]
and writing
\[
x(t)=\begin{pmatrix}y(t)\\\vdots\\y^{(r-1)}(t)\end{pmatrix}
\]
we see that the (formal) control~\eqref{eq:FC} may be expressed as
\[
u(t)=(N\circ \alpha)(\|\rho(t,x(t))\|^2)  \, \rho (t,x(t)).
\]
The feedback-controlled initial-value problem~\eqref{eq:nonlSys}~\&~\eqref{eq:FC}
may now be formulated as
\begin {equation}\label{eq:IVP-CL}
\dot x(t)=F\big(t,x(t),\fT(x)(t)\big),\quad x|_{[-h,0]}=x^0\in \cC([-h,0],\R^n),
\end{equation}
where
\[
x^0(t):=\begin{pmatrix} y^0(t)\\\vdots\\(y^0)^{(r-1)}(t)\end{pmatrix},\quad t\in [-h,0].
\]
A continuous function~$x\in \cC(I,\R^n)$ on an interval of the form~$I=[-h,\tilde\omega]$, $0 < \tilde\omega <\infty$, or of the form $[-h,\omega)$, $0 <\omega \leq\infty$,  is a {\it solution}
of~\eqref{eq:IVP-CL},
if~$x|_{[-h,0]}=x^0$, $(t,x(t))\in\cD$ for all $t\in I\backslash [-h,0)$ and
\begin{equation}\label{solution}
\forall\,t\in I,~t\ge 0:\ x(t)=x(t_0)+\int_{0}^t F\big(s,x(s),\fT(x)(s)\big)\, \dd s.
\end{equation}
A solution is {\em maximal}, if it has no right extension that is also a solution.
Since~$\fT$ is an operator with domain~$\cCC$, some care is required in interpreting the above notion of a solution~$x\in C(I,\R^n)$ when~$I$ is a bounded interval of
the form~$I=[-h,\tilde\omega]$ or~$I=[-h,\omega)$.  Let~$I$ be any such interval and write~$J:= I\backslash [-h,0)$.  Let~$x\in C(I,\R^n)$ and, for each~$\tau\in J$, define $x_\tau\in C([-h,\infty),\R^n)$ by
\[
x_\tau(t):=\left\{\begin{array}{l} x(t),~~t\in [-h,\tau]
\\x(\tau), ~~t >\tau .
\end{array}\right.
\]
With~$\fT\in\cTT$ we may associate $\tilde \fT\colon C(I,\R^n)\to \cL_{\rm loc}^\infty(J,\R^q)$ defined by the property
\[
\forall~\tau\in J:\ \tilde \fT(x)|_{[0,\tau]}=\fT(x_\tau)|_{[0,\tau]}.
\]
The causality property~(P1) of~$\fT\in\cTT$ ensures that~$\tilde \fT$ is well defined. Replacing~$\fT$ by~$\tilde \fT$ in~\eqref{solution} we arrive at the correct interpretation of a
solution.  However, for simplicity, we will not distinguish notationally between an operator~$\fT\in\cTT$ and its ``localization''~$\tilde \fT$.

It is readily verified that~$F$ has the following properties: If~$I\subset\R_{\ge 0}$ is a compact interval and~$K_n\subset\R^n$, $K_q\subset\R^q$ are compact with~$I\times K_n\subset\cD$, then
\begin{enumerate}[(a)]
\item $F(t,\cdot,\cdot)\colon K_n\times K_q\to \R^n$ is continuous for all $t\in I$;
\item $F(\cdot,v,w)\colon I\to\R^n$ is measurable for all $(v,w)\in K_n\times K_q$;
\item there exists $\hat f \in (0,\infty)$ such that $ \|F(t,v,w)\|\leq\hat f $ for almost all $t\in I$ and all $(v,w)\in K_n\times K_q$.
\end{enumerate}
Invoking~\eqref{ic}, we see that~$(0,x^0(0))\in\cD$. An application of a variant (a straightforward modification tailored to the current context) of~\cite[Thm.~B.1]{IlchRyan09}
yields the existence of a maximal solution $x:[-h,\omega)\to\R^n$, $0 <\omega\leq \infty$, of~\eqref{eq:IVP-CL} and so
\[
\cG=\textrm{graph}\big(x|_{[0,\omega)}\big)\subset\cD.
\]
Moreover, the closure of~$\cG$ is not a compact subset of~$\cD$.
\\[1ex]
{\it Step~2.}\
 Before embarking on the proof proper, we record some preliminary observations and definitions. Since~$(t,x(t))\in\cD$ for all~$t\in [0,\omega)$, we have $\pi_k(t,x(t))\in\cD_k =\text{dom}(\rho_k)$, $k=1,\ldots,r$. Introduce continuous functions
\[
e_k\colon [0,\omega)\to\cB,~~\alpha_k\colon [0,\omega)\to [1,\infty),~~\gamma_k\colon [0,\omega)\to\R^m,  ~~k=1,\ldots,r,
\]
given by
\[
e_k(t):= (\rho_k\circ\pi_k)(t,x(t)),~~\alpha_k(t):=\alpha(\|e_k(t)\|^2),~~\gamma_k(t):= \gamma(e_k(t))=\alpha_k(t)  \, e_k(t),
\]
where~$\gamma$ is given by~\eqref{eq:fcts-FC2}, and, for later notational consistency, we also write~$\gamma_0(\cdot):=0$.
Clearly, \begin{equation}\label{eq:ekfun}
\forall\, k=1,\ldots,r\ \forall\,t\in [0,\omega):\ \|e_k(t)\| < 1.
\end{equation}
In particular, for~$k=1$ we have $\|e_1(t)\|=\varphi(t)\|e(t)\| <1$
for all~$t\in [0,\omega)$ and so the tracking error $e(\cdot )=y(\cdot)-y_{\text{ref}}(\cdot)$
evolves in the funnel~$\cF_{\varphi}$.

Observe that the continuous control function~$u$ may be expressed as
\begin{equation}\label{eq:controlform}
u(t)= N(\alpha_r(t)) \, e_r(t),\quad t\in [0,\omega),
\end{equation}
and, for all~$t\in [0,\omega)$ and $k=1,\ldots,r$, we have by definition of~$\rho_k$ in~\eqref{eq:fcts-FC3}
\begin{equation}\label{eq:outputder}
e_{k}(t)-\gamma_{k-1}(t)
= \left\{\begin{array}{lll}
\varphi(t)\ e^{(k-1)}(t), && \text{if}~k\leq\hat r
\\
\varphi (t)\ y^{(k-1)}(t), && \text{otherwise}.
\end{array}\right.\quad
\end{equation}
We also record that
\begin{equation}\label{eq:dotgamk}
\dot\alpha_k(t)=-2\alpha^\prime \big(\|e_k(t)\|^2\big) \langle e_k(t),\dot e_k(t)\rangle
\qquad \text{for a.a. $t\in [0,\omega)$},~~k=1,\ldots, r.
\end{equation}
Define functions $\psi_k\colon [0,\infty)\to\R^m$, $k=1,\ldots,r$, as follows
\begin{align*}
&\hat r =r~\implies~ \psi_k(\cdot ) :=0,~k=1,\ldots, r
\\[1ex]
&\hat r < r ~\implies~\psi_k(t):=\left\{\begin{array}{cll}
0, && \text{if}~k < \hat r
\\[1ex]
-\varphi(t)y_{\rm{ref}}^{(\hat r)}(t), && \text{if}~k=\hat r
\\[1ex]
\dphi (t)y_{\rm{ref}}^{(k-1)}, &&\text{if}~ \hat r < k \leq r-1
\\[1ex]
\dphi (t)y_{\rm{ref}}^{(r-1)}(t) +\varphi(t) y_{\rm{ref}}^{(r)}(t), && \text{if}~k=r.\end{array}\right.
\end{align*}
By choice of the design parameters as in~\eqref{eq:fcts-FC}, $\varphi$ is bounded (and so~$\dphi$ is essentially bounded
{by the definition of~$\Phi$})
if~$\hat r < r$.  Therefore, we may infer the existence of~$\psi^*\in\R$ (with~$\psi^*=0$ if~$\hat r=r$) such that
\begin{equation}\label{psikbound}
\|\psi_k(t)\|\leq\psi^*~~\text{for a.a.}~t\in [0,\infty),~~k=1,\ldots,r.
\end{equation}
Observe that, for almost all $t\in [0,\omega)$,
\begin{equation}\label{eq:errors}
\left.
\begin{array}{l}
 \dot e_k(t)=\dphi(t)e^{(k-1)}(t)+e_{k+1}(t)-\gamma_k(t)+\dot\gamma_{k-1}(t)+\psi_k(t),~~k=1,\ldots,r-1
 \\[1ex]
 \dot e_r(t)= \dphi(t) e^{(r-1)}(t)+\varphi(t) e^{(r)}(t)+\dot\gamma_{r-1}(t)+\psi_r(t)
 \end{array} \right\}
\end{equation}
which, if $r=1=\hat r$, collapses to the tautology: $\dot e_1(t)=(\varphi e)^{(1)}(t)$ for a.a. $t\in [0,\omega)$.

Arbitrarily fix $\tau\in (0,\omega)$.  By continuity, there exists $\theta\in (0,\infty)$ such that
\begin{equation}\label{mubound}
\forall\,t\in [0,\tau]:\ \big(1+\varphi(t)\big)\sum_{k=1}^r \|e^{(k-1)}(t)\|\leq \theta
\end{equation}
and so, by properties of~$\Phi$, there exists~$c >0$ such that
\begin{equation}\label{b1}
\|\dphi(t)e^{(k-1)}(t)\|\leq c \big(1+\varphi(t)\big)\|e^{(k-1)}(t)\|\leq c\,\theta~~\text{for a.a.\, $t\in [0,\tau]$},~~k=1,\ldots,r.
\end{equation}
Again by properties of~$\Phi$, the following are well defined:
\[
\sup_{t\in [\tau,\infty)} \left(\frac{1}{\varphi (t)}\right) =:\lambda > 0\quad\text{and}\quad \esup_{t\in [\tau,\infty)} \left({{\frac{|\dot\varphi(t)|}{\varphi(t)}}}\right)=: \mu  \ge 0.
\]
For $k\in\{1,\ldots,r\}$ and invoking~\eqref{eq:ekfun},~\eqref{eq:outputder} and~\eqref{psikbound},   we find
\begin{align*}
{\rm (a)}\quad\|\dphi(t)e^{(k-1)}(t)\|&\leq \mu\|\varphi(t)e^{(k-1)}(t)\|\leq \mu\big(1+\|\gamma_{k-1}(t)\|\big)
\\
&~~\hspace{3.5cm}\text{for a.a.}\, t\in [\tau,\omega),~~\text{if}~k \leq \hat r
\\[1ex]
{\rm(b)}\quad \|\dphi(t)e^{(k-1)}(t)\|&\leq \mu\|\varphi (t)y^{(k-1)}(t)\|+\|\dphi(t)y_{\rm{ref}}^{(k-1)}\|
\\
&\leq  \mu\big(1+\|\gamma_{k-1}(t)\|\big)+\psi^*
\\
&\hspace{3.5cm}~~\text{for a.a.}\, t\in [\tau,\omega),~~\text{if}~k> \hat r
\end{align*}
and so, {\it a fortiori}, we have
\begin{equation}\label{b3}
\|\dphi(t)e^{(k-1)}(t)\|\leq \mu\big(1+\|\gamma_{k-1}(t)\|\big)+\psi^*~~\text{for a.a.}\, t\in [\tau,\omega),~~k=1,\ldots,r.
\end{equation}
We complete the preliminaries by writing
\[
\hat\varepsilon_k := \max_{t\in [0,\tau]}\|e_k(t)\|^2 < 1,~~ k=1,\ldots,r.
\]
{\it Step~3.} \
Assume that $r\ge 2$, otherwise proceed to Step~5.
Let~$\varepsilon_1$ be the unique point of ($0,1)$ such that $\varepsilon_1\alpha (\varepsilon_1)=1+\mu+2\psi^*$ and $\varepsilon_1^* := \max\left\{\hat\varepsilon_1\,,\,\varepsilon_1\right\} < 1$.  We will show that
\begin{equation}\label{e1bound}
\forall\, t\in [0,\omega):\ \|e_1(t)\|^2\leq \varepsilon_1^*.
\end{equation}
Suppose that this claim is false.  Then $\|e_1(s)\|^2 >\varepsilon_1^*$ for some $s\in (0,\omega)$.   Since $\|e_1(t)\|^2\leq\hat\varepsilon_1 \leq \varepsilon_1^*$ for all $t\in [0,\tau]$,
we have~$\tau < s$ and so we may define
\[
\sigma:=\max\setdef{t\in [\tau,s) }{ \|e_{1} (t)\|^2=\varepsilon_1^*}.
\]
Clearly,
\[
\forall\,t\in [\sigma,s]:\ \varepsilon_1\leq\varepsilon_1^*\leq \|e_1(t)\|^2,
\]
whence, by monotonicity of~$\alpha$,
\[
\forall\,t\in [\sigma,s]:\ \alpha(\varepsilon_1 )\leq\alpha\big(\|e_1(t)\|^2\big)
=\alpha_1(t).
\]
Therefore,
\begin{equation}\label{bound1}
\forall\,t\in [\sigma,s]:\ \alpha_1(t)\|e_1(t)\|^2 \geq \varepsilon_1\alpha (\varepsilon_1)=1+\mu+2\psi^*
\end{equation}
which, by the first of relations~\eqref{eq:errors} in conjunction with
\eqref{eq:ekfun} and~\eqref{b3} (and recalling $\gamma_0(\cdot) = 0$), gives
\begin{align*}\label{eq:ineq}
\half\ddts\|e_1(t)\|^2&= \langle e_{1}(t), \dot e_{1}(t)\rangle\nonumber
\\
& =\langle \dphi(t) e_1(t), e(t)\rangle +  \langle e_{1}(t),e_{2}(t)\rangle -\alpha_{1}(t)\|e_1(t)\|^2+  \langle e_1(t),\psi_1(t)\rangle \nonumber
\\
& < 1+\mu+2\psi^*  - \alpha_1(t)\|e_1(t)\|^2 \leq 0
\end{align*}
for almost all $t\in [\sigma,s]$ and so $\|e_1(s)\|^2 < \|e_1(\sigma)\|^2$,  whence the contradiction
\[
\varepsilon_1^* < \|e_1(s)\|^2 < \|e_1(\sigma)\|^2 =\varepsilon_1^*.
\]
Therefore~\eqref{e1bound} holds.
\\[1ex]
{\it Step~4.} \
For notational convenience,  write
\[
\cW_1:= \cW^{1,\infty}([0,\omega),\R)\quad\text{and}\quad \cW_m := \cW^{1,\infty}([0,\omega),\R^m).
\]
We show by induction that
\begin{equation}\label{eq:induction-gamma-e-nu}
(\alpha_k,e_k,\gamma_k)\in\cW_1\times\cW_m\times\cW_m
\quad \text{for}\quad
k=1,\ldots,r-1.
\end{equation}
This step is vacuous in the case~$r=1$.
Let~$k=1$.  By~\eqref{e1bound}, we see that~$e_1$ is bounded by~$\sqrt{\varepsilon_1^*}$, $\alpha_1$ is bounded by~$\alpha(\varepsilon_1^*)$ and
that~$\gamma_1$ is bounded by $\sqrt{\varepsilon_1^*}\alpha(\varepsilon_1^*)$.
Recalling that~$\gamma_0(\cdot) = 0$, essential boundedness of~$\dot e_1$  follows by the first of relations~\eqref{eq:errors} together
with~\eqref{eq:ekfun}, \eqref{psikbound}, \eqref{b1}, \eqref{b3}.
Invoking~\eqref{eq:dotgamk}, we may conclude essential boundedness of~$\dot\alpha_1$.
Essential boundedness of $\dot\gamma_1 =\alpha_1\dot e_1+\dot\alpha_1 e_1$ then follows.  Therefore,
$(\alpha_1,e_1,\gamma_1)\in\cW_1\times\cW_m\times\cW_m$.

 Now assume that $k\in \{2,\ldots,r-1\}$ and
\[
\big(\alpha_j,e_j,\gamma_j\big)\in \cW_1\times\cW_m\times\cW_m,
\qquad j=1,\ldots,k-1.
\]
Set
\[
\beta := \max\left\{\psi^*\,,\,{\textstyle{\sup_{t\in [0,\omega)}}}\|\gamma_{k-1}(t)\|\,,\,\esup_{t\in [0,\omega)}\|\dot\gamma_{k-1}(t)\|\right\}<\infty .
\]
By~\eqref{eq:ekfun},~\eqref{eq:errors} and~\eqref{b3}, we have
\begin{align}\label{eq:ineq2}
\langle e_{k}(t), \dot e_{k}(t)\rangle & =\dot\varphi (t)  \langle e_k(t),e^{(k-1)}(t)\rangle +  \langle e_{k}(t), e_{k+1}(t)\rangle \nonumber
\\
&\quad+  \langle e_k(t),\big(\dot\gamma_{k-1}(t)+\psi_k(t)\big)\rangle -\alpha_{k}(t)\|e_k(t)\|^2\nonumber
\\
&< 1+3\beta +\mu(1 + \beta) - \alpha_k(t)\|e_k(t)\|^2
\end{align}
for almost all $t\in [\tau,\omega)$. Let~$\varepsilon_k$ be the unique point of~$(0,1)$  such that $\varepsilon_k\alpha (\varepsilon_k)=1+3\beta +\mu(1+\beta)$ and define
\[
\varepsilon_k^* := \max\left\{\hat\varepsilon_k\,,\,\varepsilon_k\right\} > 0.
\]
We first show that
\begin{equation}\label{ekbound}
\forall\, t\in [0,\omega):\ \|e_k(t)\|^2\leq\varepsilon_k^*
\end{equation}
by the contradiction argument of Step~3 ({\em mutatis mutandis}).
Suppose that~\eqref{ekbound} is false.
Then $\|e_k(s)\|^2 >\varepsilon_k^*$ for some~$s\in (0,\omega)$.   Since $\|e_k(t)\|^2\leq\hat\varepsilon_k \leq \varepsilon_k^*$ for all~$t\in [0,\tau]$,
we have~$\tau < s$ and so we may define
$\sigma:=\max\setdef{t\in [\tau,s)}{\|e_{k} (t)\|^2=\varepsilon_k^*}$.
The counterpart of~\eqref{bound1} now follows:
\[
\forall\,t\in [\sigma,s]:\ \alpha_k(t)\|e_k(t)\|^2 \geq \varepsilon_k\alpha (\varepsilon_k)=1+3\beta+\mu(1+\beta)
\]
which, in conjunction with~\eqref{eq:ineq2}, gives $\half\ddts \|e_k(t)\|^2 < 0$ for almost all $t\in [\sigma,s]$, whence the contradiction
\[
\varepsilon_k^* < \|e_k(s)\|^2 < \|e_k(\sigma)\|^2 =\varepsilon_k^*.
\]
Therefore,~\eqref{ekbound} holds which, in turn, implies that~$\alpha_{k}$ is bounded (by~$\alpha (\varepsilon_k^*)$) and that~$\gamma_k=\alpha_ke_k$ is bounded (by
$\sqrt{\varepsilon_k^*}\, \alpha (\varepsilon_k^*)$).
By boundedness of~$e_{k+1}$, $\gamma_k$ and essential boundedness of~$\dot\gamma_{k-1}$,
it follows from ~\eqref{eq:errors}, together
with~\eqref{b1} and~\eqref{b3}, that~$\dot e_k$ is essentially bounded and so~$e_k\in\cW_m$.
Invoking~\eqref{eq:dotgamk},
we may now infer essential boundedness of~$\dot\alpha_k$.  Therefore, $\alpha_k\in\cW_1$.  Finally, since
$\dot\gamma_k=\alpha_k \dot e_k + \dot\alpha_k e_k$, we have essential boundedness of~$\dot\gamma_k$ and so~$\gamma_k\in\cW_m$.  In summary, we have shown that,
for~$k\in\{2,\ldots,r-1\}$,
\begin{equation*}
\big(\alpha_j,e_j,\gamma_j\big)\in \cW_1\times\cW_m\times\cW_m,~j=1,\ldots,k-1~~\Longrightarrow~~ (\alpha_k,e_k,\gamma_k)\in\cW_1\times\cW_m\times\cW_m,
\end{equation*}
and so, by induction, we conclude~\eqref{eq:induction-gamma-e-nu}.
\\[1ex]
{\it Step~5.}  \
Our next goal is to prove boundedness of the solution~$x$.  Recalling that $y_{\rm ref}\in\cW^{r,\infty}(\R_{\ge 0},\R^m)$, it suffices to show that the output error~$e$ and its derivatives $\dot e,\ldots, e^{(r-1)}$ are bounded on~$[0,\omega)$.  By~\eqref{mubound}, we already know that
\[
\forall\,k=1,\ldots,r~~ \forall\,t\in [0,\tau]:\ \|e^{(k-1)}(t)\|\leq \theta,
\]
and so it remains to show that~$e^{(k-1)}$ is bounded on $[\tau,\omega)$, $k=1,\ldots,r$.  Since $\varphi(t)e(t)=e_1(t)\in\cB$ for all $t\in [0,\omega)$, we have
\[
\forall\, t\in [\tau,\omega):\ \|e(t)\|\leq \frac{1}{\varphi(t)}\leq \lambda.
\]
By boundedness of the functions~$\gamma_k$ (Step~4), there exists~$\gamma^*>0$ such that
\begin{equation}\label{nustar}
\forall\,k=2,\ldots,r~~ \forall\,t\in [0,\omega):\ \|\gamma_{k-1}(t)\|\leq \gamma^*.
\end{equation}
Let~$k\in\{2,\ldots,r\}$.  By
\eqref{eq:outputder}, we have
\[
\forall\,t\in [\tau,\omega):\ \|e^{(k-1)}(t)\|\leq \lambda\left(1+\gamma^*\right) +\sup_{t\geq \tau}\|y^{(k-1)}_{\rm{ref}}(t)\| < \infty.
\]
This completes Step~5.
\\[1ex]
{\it Step~6.}  \
We prove boundedness
of~$\alpha_r{\colon [0,\omega) \to [1,\infty),\ t\mapsto \alpha(\|e_r(t)\|^2)}$ { together with an
immediate consequence thereof:}
\begin{equation}\label{erbound}
\exists\,  \varepsilon_r^*\in (0,1) \
\forall\,t\in [0,\omega):\ \|e_r(t)\|^2 \leq \varepsilon_r^*.
\end{equation}
By boundedness of~$x$ (Step~5) and property~(TP3)
of the operator class~$\cTT$, there exists compact~$K_q\subset\R^q$ such that~$\fT( x)(t)\in K_q$ for almost all $t\in [0,\omega)$.
Since $d\in\cL^\infty(\R_{\ge 0},\R^p)$, there exists  compact~$K_p\subset\R^p$ such that $d(t)\in K_p$ for almost all~$t\in [0,\omega)$.
By the high-gain property, there exists $v^*\in (0,1)$ such that the  continuous function
\[
\chi\colon\R\to\R,~~s \mapsto \min\setdef{  \langle v,f(\delta,z,-s v)\rangle }{ (\delta,z,v)\in K_p \times K_q   \times  A_m }
\]
is unbounded from above, where, for notational convenience, we have introduced the compact annulus
\[
A_m := \setdef{v\in\R^m }{ v^*\leq \|v\| \leq 1}.
\]
Choose a real sequence
$(s_j)$ such that  the sequence $\big(\chi (s_j)\big)$  is unbounded,
positive, and strictly increasing.
By surjectivity and continuity of~$N$, for every~$a\in\R_{\ge 0}$ and every $b\in\R$, the set $\setdef{\kappa >a }{ N(\kappa )=b}$ is non-empty. Choose
$\kappa_1> \alpha ({(1-v^*)^2})+\alpha_r (0)$
such that~$N(\kappa_1)=s_1$ and
define the strictly increasing sequence $(\kappa_j)$ by the recursion
\[\kappa_{j+1}:= \inf\setdef{\kappa >\kappa_j}{N(\kappa)=s_{j+1}}.\]
Observe that
\[
    \lim_{j\to \infty} \chi (N(\kappa_j))= \lim_{j\to \infty} \chi (s_j) = \infty.
\]
Seeking a contradiction, suppose
that~$\alpha_r(\cdot)$
is not bounded.  Then, since $\kappa_{j+1} >\kappa_1 >\alpha_r(0)$ for all $j\in\N$,  the sequence $(\tau_j)$  in $(0,\omega)$ defined by
\[
\tau_{j}= \inf\setdef{t\in [0,\omega)}{~\alpha_r (t)=\kappa_{j+1}},\quad j\in\N_0,
\]
is well-defined and strictly increasing with $N(\alpha_r(\tau_j))=N(\kappa_{j+1})=s_{j+1}$ for each~$j\in\N_0$.   Now, define the sequence~$(\sigma_j)$ in $(0,\omega)$ by
\[
\sigma_j = \sup\setdef{t\in [\tau_{j-1},\tau_j]}{~ \chi(N(\alpha_r(t)))=\chi(s_j)},\quad j\in\N.
\]
Since the sequence~$\big(\chi (s_j)\big)$ is strictly increasing, we have
\[
\forall\,j\in\N:\ \chi(N(\alpha_r(\sigma_j)))=\chi (s_j)
<\chi (s_{j+1})=\chi(N(\alpha_r(\tau_j))),
\]
and so
\begin{equation}\label{chiineq}
\forall\,j\in\N\ \forall\,t\in (\sigma_j,\tau_j]:\  \sigma_j <\tau_j\quad\text{and}\quad \chi(N(\alpha_r(\sigma_j)))=\chi(s_j) < \chi(N(\alpha_r(t))).
\end{equation}
Next, suppose that, for some~$j\in\N$, there exists~$t\in [\sigma_j,\tau_j]$ such that~$e_r(t)\not\in A_m$. We first show that $\alpha_r(t)\ge \kappa_j$.
If $\alpha_r(t)<\kappa_j$, then $\alpha_r(\tau_j) = \kappa_{j+1} > \kappa_j$ and continuity of~$\alpha_r$ imply that there exists $\tilde t\in(\sigma_j,\tau_j)$ such that $\alpha_r(\tilde t) = \kappa_j$, thus
\[
    \chi(N(\alpha_r(\tilde t)))= \chi(N(\kappa_j))=\chi(s_j),
\]
which contradicts the definition of~$\sigma_j$. Therefore, $\alpha_r(t)\ge \kappa_j$ which, together with the supposition~$\|e_r(t)\| < {1-v^*}$,
leads to the contradiction:
\[
\alpha ({(1-v^*)^2}) < \kappa_1 \leq  \kappa_j \leq \alpha_r(t)=\alpha\big(\|e_r(t)\|^2\big) < \alpha ({(1-v^*)^2}).
\]
Therefore,
\begin{equation}\label{erKm}
\forall\,j\in\N\ \forall \,t\in [\sigma_j,\tau_j]:\ e_r(t)\in A_m,
\end{equation}
which, in conjunction with the facts that~$d(t)\in K_p$ and~$(\fT x)(t)\in K_q$ for almost all~$t\in [0,\omega)$ and invoking~\eqref{chiineq}, yields
\begin{align}\label{chiineq2}
\langle e_r(t), f(d(t),(&\fT x)(t),u(t))\rangle = -\langle -e_r(t), f\big(d(t),(\fT x)(t), -N(\alpha_r(t))(-e_r(t))\big)\rangle\nonumber
\\[1ex]
& \leq -\min\setdef{\langle v,f\big(\delta,z,-N(\alpha_r (t))v\big)\rangle }{ (\delta,z,v)\in K_p\times K_q\times  A_m}\nonumber
\\[1ex]
&
=-\chi\big(N(\alpha_r(t))\big) \leq -\chi (s_{j})
\end{align}
for all $j\in\N$ and almost all $t\in [\sigma_j,\tau_j]$. By~\eqref{b1}, \eqref{b3} and~\eqref{nustar},
\[
\|\dphi(t)e^{(r-1)}(t)\|\leq c\,\theta + \mu (1+\gamma^*)+\psi^*=: \theta^*~~\text{for a.a.}\, t\in [0,\omega).
\]
Since $e^{(r)}(t)=f(d(t),\fT(x)(t),u(t))-y_{\textrm{ref}}^{(r)}(t)$ for almost all~$t\in [0,\omega)$ and recalling the last of relations~\eqref{eq:errors}, we have
\[
\dot e_r (t)  =\varphi (t)\big(f(d(t),\fT(x)(t),u(t))-y_{\textrm{ref}}^{(r)}(t)\big) + \dphi(t)e^{(r-1)}(t)+\dot\gamma_{r-1}(t) +\psi_r(t)
\]
for almost all $t\in [0,\omega)$. By~{\eqref{eq:induction-gamma-e-nu},}
$\dot\gamma_{r-1}$ is essentially bounded and, since $y_{\text{\rm ref}} \in \cW^{r,\infty}(\R_{\ge 0},\R^m)$,
we have essential boundedness of~$y_{\textrm{ref}}^{(r)}$.
Write
\[
c_1:=\theta^* +\psi^*+\esup_{t\in [0,\omega)}\|\dot\gamma_{r-1}(t)\|\quad\text{and}\quad c_2:= \esup_{t\ge 0} \|y_{\textrm{ref}}^{(r)}(t)\|.
\]
Invoking~\eqref{eq:ekfun}, \eqref{psikbound} and~\eqref{chiineq2}, we arrive at
\[
\half\ddts \|e_r(t)\|^2 \leq c_1 -\varphi (t)\big(\chi (s_j)-c_2\big)
\]
for all $j\in\N$ and almost all $t\in [\sigma_j,\tau_j]$. By properties of~$\varphi\in\Phi$ and noting that~$\sigma_1 >0$, we have $\inf_{t\in [\sigma_1,\infty)}\varphi (t) >0$.   Since~$\chi (s_j)\to \infty$ as~$j\to\infty$, we
may choose~$j$ sufficiently large so that {$c_1 -\varphi (t)\big(\chi (s_j)-c_2\big) < 0$}
for almost all $t\in [\sigma_j,\tau_j]$,
in which case we have $\|e_r (\tau_j) \|^2 <\|e_r(\sigma_j)\|^2$ and so
\[
\alpha_r(\tau_j)=\alpha\big(\|e_r(\tau_j)\|^2\big) < \alpha\big(\|e_r(\sigma_j)\|^2\big)=\alpha_r(\sigma_j)
\]
which is impossible since, by definition of~$\tau_j$, we have $\alpha_r(t) < \alpha_r(\tau_j)$ for all~$t\in [0,\tau_j)$.  Therefore, our original supposition that~$\alpha_r$ is unbounded is false.
{This proves~\eqref{erbound}
and completes the proof of Step~6.}
\\[1ex]
{\it Step~7.}  \
We prove Assertion\,\ref{item-main-1} of the theorem.  Recalling inequalities~\eqref{e1bound}, \eqref{ekbound} and~\eqref{erbound} of Steps~3, 4 and~6, we have
\begin{equation}\label{ekeps}
\|e_k(t)\|\leq \varepsilon:=\sqrt{\max\{\varepsilon_1^*,\ldots,\varepsilon_r^*\}}< 1
\end{equation}
for all $t\in [0,\omega)$ and all $k=1,\ldots,r$. Define
\[
\widehat\cD_r :=\setdef{(\eta_1,\ldots,\eta_r)\in \R^{rm}}{\|\rho_k(\eta_1,\ldots,\eta_k)\|\leq \varepsilon,~k=1,\ldots,r},
\]
which is evidently a compact subset of~$\cD_r$ as in~\eqref{eq:fcts-FC3}.  Since $e_k(t) =(\rho_k\circ \pi_k)(t,x(t))$ for all $t\in [0,\omega)$, $k=1,\ldots,r$, it follows that
$\pi_r(t,x(t))\in\widehat\cD_r$ for all~$t\in [0,\omega)$.  Suppose that~$\omega <\infty$.  Then
\[
\forall\, t\in [0,\omega):\ (t,x(t))\in \widehat\cD :=\setdef{(s,\xi)\in [0,\omega]\times \R^{rm} }{ \pi_r(s,\xi)\in\widehat\cD_r}\subset \cD.
\]
By compactness of~$\widehat\cD$ it follows that the closure of $\textrm{graph}\big(x|_{[0,\omega)}\big)$ is a compact subset of~$\cD$, which contradicts the findings of Step~1.
Therefore,~$\omega =\infty$.
\\[1ex]
{\it Step~8.} \
We complete the proof by establishing Assertions\,\ref{item-main-2},~\ref{item-main-3} and~\ref{item-main-4}.
Assertion\,\ref{item-main-2} is a direct consequence of Assertion\,\ref{item-main-1} and the results of Steps~5~\&~6.
Recalling that~$e_1=\varphi e$, we may infer Assertion\,\ref{item-main-3} from~\eqref{ekeps} and Assertion\,\ref{item-main-1}.
Assertion\,\ref{item-main-4} follows by Assertion\,\ref{item-main-1}
and~\ref{item-main-3}, together with~\eqref{eq:ekfun},~\eqref{eq:outputder} and~\eqref{nustar}.
\\[1ex]
{\it Step 9.}
Assume that the negative-definite (respectively, positive-definite) high-gain
property is known to hold. Steps~1--5 of the above are
unaffected by this assumption.  Step~6 is readily modified as follows.  By the assumption, there exists a positive
(respectively, negative)
real sequence~$(s_j)$ such that  the sequence~$\big(\chi (s_j)\big)$  is unbounded,
positive, and strictly increasing.
{ Replacing ~$N$ by $s\mapsto s$
(respectively, by $s\mapsto -s$)}, the remaining arguments of
Step~6 apply {\em mutatis mutandis} to conclude boundedness of $\alpha_r$. Steps~7 and~8 then follow as before. This completes the proof of the theorem.
\end{proof}

\begin{proof}[\bf Proof of Corollary~\ref{Prop:FunCon-Bounds}]
Let $y_{\textrm{ref}}\in W^{r,\infty}(\R_{\ge 0},\R^m)$ and $y^0\in W^{r,\infty}([-h,0],\R^m)$.  By Theorem~\ref{Thm:FunCon-Nonl},
the feedback-controlled system~\eqref{eq:nonlSys}~\& \eqref{eq:FC}
has a solution, every solution can be maximally extended and every maximal solution is global. Let $y\colon [-h,\infty)\to\R^m$
be any such global solution.
In the following, we adopt the notation introduced in the proof of Theorem~\ref{Thm:FunCon-Nonl} and recall that, for all $k=1,\ldots,\hat r-1$,
$\psi_k (\cdot)=0$,  and for all $t\ge 0$,
\begin{align*}
\|e_k (t)\|&<1,\quad \|\gamma_k(t)\|  =  \alpha(\|e_k(t)\|^2)\ \|e_k(t)\|,
\\[.8ex]
\|\dot\gamma_k(t)\|
&=\left\|2\alpha^\prime(\|e_k(t)\|^2)\, \langle e_k(t),\dot e_k(t)\rangle \, e_k(t)
    +\alpha(\|e_k(t)\|^2) \, \dot e_k(t)\right\|
\\[.8ex]
&\leq\tilde\alpha (\|e_k(t)\|^2)\|\dot e_k(t)\|.
\end{align*}
Invoking \eqref{eq:ekfun}, \eqref{eq:outputder} and \eqref{eq:errors}, with the convention that $\gamma_0(\cdot )\equiv 0\equiv\dot\gamma_0(\cdot)$,
we have, for almost all~$t\ge0$,
\[
\begin{array}{rcl}
\|\dot e_k(t)\|&=& \left\|(\dphi(t)/\varphi(t)) \
\big(e_k(t)-\gamma_{j-1}(t)\big)+e_{k+1}(t)+\dot\gamma_{k-1}(t)-\gamma_k(t)\right\|
\\[.8ex]
 &\leq & M_k(t)+\|\gamma_k(t)\|,
\\[.8ex]
\langle e_k(t),\dot e_k(t)\rangle
&\leq &
M_k(t)-\alpha_k(t)\ \|e_k(t)\|^2,
\\[.8ex]
M_k(t)&:=& 1+\mu_0 \big(\|e_k(t)\|+ \|\gamma_{k-1}(t)\|\big)  +\|\dot\gamma_{k-1}(t)\|.
\end{array}
\]
Setting $k=1$, we have
\[
\langle e_1(t),\dot e_1(t)\rangle \leq \mu_0+1 -\alpha_1(t)\ \|e_1(t)\|^2~~\text{for a.a. $t\geq 0$}.
\]
With $e_1^0$ and $c_1$ as in~\eqref{ckbounds},
the argument used in
Step~3 of the proof of Theorem~\ref{Thm:FunCon-Nonl}  applies,
{\it mutatis mutandis}, to conclude that $\|e_1(t)\| \leq c_1$ for all $t\geq 0$.

With $\mu_1=1+\mu_0c_1$ as in~\eqref{ckbounds} we have, for almost all $t\ge 0$,
\[
\|\gamma_1(t)\|
\leq c_1 \alpha(c_1^2) \quad \text{and} \quad \|\dot\gamma_1(t)\|
\leq \tilde\alpha (c_1^2) \big(\mu_1+c_1\alpha(c_1^2)\big),
\]
wherein we have used the facts that $\alpha$ and $\tilde\alpha$ are non-decreasing functions (monotonicity of the latter being assured by
the assumption of monotonicity of $\alpha^\prime$).
Now set $k=2$, in which case we have
\[
M_2(t)\leq 1+\mu_0\big(1+c_1\alpha(c_1^2)\big)+\tilde\alpha (c_1^2)
     \big(\mu_1+c_1\alpha(c_1^2)\big)=\mu_2~~\text{for a.a. $t\geq 0$}.
\]
With $e_2^0$ and $c_2$ as in~\eqref{ckbounds}
the argument used in Step 4 of the proof of Theorem~\ref{Thm:FunCon-Nonl}  applies,
{\it mutatis mutandis}, to conclude that  $\|e_2(t)\| \leq c_2$ for all $t\geq 0$.
Iterating this process, we arrive at
\[
\forall\, k=1,\ldots,\hat r -1\ \forall\,t\geq 0:\ \|e_k(t)\| \leq c_k
\]
To complete the proof, simply note that, for all $t\geq 0$,
\[
\begin{array}{l}
\varphi(t)\|e(t)\|=\|e_1(t)\|\leq {c_1}\quad\text{and}
\\[.8ex]
\varphi (t)\|e^{(k)}(t)\|=\|e_{k+1}(t)-\gamma_{k}(t)\|\leq {c_{k+1}}+c_{k} \alpha(c_{k}^2),~~k=1,\ldots,\hat r-2.  \hspace*{1.4cm}\qedhere
\end{array}
\]
\end{proof}


\end{appendices}


%
%

\bibliographystyle{spmpsci}      

\begin{thebibliography}{10}
\providecommand{\url}[1]{{#1}}
\providecommand{\urlprefix}{URL }
\expandafter\ifx\csname urlstyle\endcsname\relax
  \providecommand{\doi}[1]{DOI~\discretionary{}{}{}#1}\else
  \providecommand{\doi}{DOI~\discretionary{}{}{}\begingroup
  \urlstyle{rm}\Url}\fi

\bibitem{BechRovi08}
Bechlioulis, C.P., Rovithakis, G.A.: Robust adaptive control of feedback
  linearizable {MIMO} nonlinear systems with prescribed performance.
\newblock {IEEE} Trans. Autom. Control \textbf{53}(9), 2090--2099 (2008)

\bibitem{BechRovi14}
Bechlioulis, C.P., Rovithakis, G.A.: A low-complexity global approximation-free
  control scheme with prescribed performance for unknown pure feedback systems.
\newblock Automatica \textbf{50}(4), 1217--1226 (2014)

\bibitem{Berg14c}
Berger, T.: Zero dynamics and stabilization for linear {DAE}s.
\newblock In: S.~Sch\"ops, A.~Bartel, M.~G\"unther, E.J.W. ter Maten, P.C.
  M\"uller (eds.) Progress in Differential-Algebraic Equations,
  Differential-Algebraic Equations Forum, pp. 21--45. Springer-Verlag,
  Berlin-Heidelberg (2014)

\bibitem{Berg20pp}
Berger, T.: Funnel control of the Fokker-Planck equation for a multi-dimensional Ornstein-Uhlenbeck process (2020).
\newblock Submitted for publication, preprint available at
  \url{https://arxiv.org/abs/2005.13377v2}

\bibitem{Berg20}
Berger, T.: Tracking with prescribed performance for linear non-minimum phase
  systems.
\newblock Automatica \textbf{115}, Article 108909 (2020).

\bibitem{BergBrei19}
Berger, T., Breiten, T., Puche, M., Reis, T.: Funnel control for the monodomain
  equations with the {FitzHugh-Nagumo} model (2019).
\newblock Submitted for publication, preprint available at
  \url{https://arxiv.org/abs/1912.01847}

\bibitem{BergLanz20}
Berger, T., Lanza, L.: Output tracking for a non-minimum phase robotic
  manipulator (2020).
\newblock To appear in Proceedings of the MTNS 2020, preprint available at
  \url{https://arxiv.org/abs/2001.07535}

\bibitem{BergLe18}
Berger, T., L{\^e}, H.H., Reis, T.: Funnel control for nonlinear systems with
  known strict relative degree.
\newblock Automatica \textbf{87}, 345--357 (2018).

\bibitem{BergOtto19}
Berger, T., Otto, S., Reis, T., Seifried, R.: Combined open-loop and funnel
  control for underactuated multibody systems.
\newblock Nonlinear Dynamics \textbf{95}, 1977--1998 (2019).

\bibitem{BergPuch19}
Berger, T., Puche, M., Schwenninger, F.L.: Funnel control for a moving water
  tank (2020).
\newblock Submitted for publication, preprint available at
  \url{https://arxiv.org/abs/1902.00586v2}

\bibitem{BergPuch20a}
Berger, T., Puche, M., Schwenninger, F.L.: Funnel control in the presence of
  infinite-dimensional internal dynamics.
\newblock Syst. Control Lett. \textbf{139}, Article 104678 (2020)

\bibitem{BergRaue18}
Berger, T., Rauert, A.L.: A universal model-free and safe adaptive cruise
  control mechanism.
\newblock In: Proceedings of the MTNS 2018, pp. 925--932. Hong Kong (2018)

\bibitem{BergRaue20}
Berger, T., Rauert, A.L.: Funnel cruise control.
\newblock Automatica \textbf{119}, Article 109061 (2020)

\bibitem{BergReis14a}
Berger, T., Reis, T.: Zero dynamics and funnel control for linear electrical
  circuits.
\newblock J. Franklin Inst. \textbf{351}(11), 5099--5132 (2014)

\bibitem{BergReis18a}
Berger, T., Reis, T.: Funnel control via funnel pre-compensator for minimum
  phase systems with relative degree two.
\newblock {IEEE} Trans. Autom. Control \textbf{63}(7), 2264--2271 (2018).

\bibitem{BergReis18b}
Berger, T., Reis, T.: The {f}unnel {p}re-{c}ompensator.
\newblock Int. J. Robust \& Nonlinear Control \textbf{28}(16), 4747--4771
  (2018).

\bibitem{ByrnWill84}
Byrnes, C.I., Willems, J.C.: Adaptive stabilization of multivariable linear
  systems.
\newblock In: Proc. 23rd~{IEEE} Conf. Decis. Control, pp. 1574--1577 (1984)

\bibitem{GeHong04}
Ge, S.S., Hong, F., Lee, T.H.: Adaptive neural control of nonlinear time-delay
  systems with unknown virtual control coefficients.
\newblock {IEEE} Trans. Syst. Man Cyb. {B} \textbf{34}(1), 499--516 (2004)

\bibitem{GeWang02}
Ge, S.S., Wang, J.: Robust adaptive neural control for a class of perturbed
  strict feedback nonlinear systems.
\newblock {IEEE} Trans. Neural Netw. \textbf{13}(6), 1409--1419 (2002)

\bibitem{GeWang03}
Ge, S.S., Wang, J.: Robust adaptive tracking for time-varying uncertain
  nonlinear systems with unknown control coefficients.
\newblock {IEEE} Trans. Autom. Control \textbf{48}(8), 1463--1469 (2003)


\bibitem{Hack14}
Hackl, C.M.: Funnel control for wind turbine systems.
\newblock In: Proc. 2014 {IEEE} Int. Conf. Contr. Appl., Antibes, France, pp.
  1377--1382 (2014)

\bibitem{Hack15a}
Hackl, C.M.: Current {PI}-funnel control with anti-windup for synchronous
  machines.
\newblock In: Proc. 54th~{IEEE} Conf. Decis. Control, Osaka, Japan, pp.
  1997--2004 (2015)

\bibitem{Hack15b}
Hackl, C.M.: Speed funnel control with disturbance observer for wind turbine
  systems with elastic shaft.
\newblock In: Proc. 54th~{IEEE} Conf. Decis. Control, Osaka, Japan, pp.
  12005--2012 (2015)

\bibitem{Hack17}
Hackl, C.M.: Non-identifier Based Adaptive Control in Mechatronics~--~Theory and
  Application, \emph{Lecture Notes in Control and Information Sciences}, vol.
  466.
\newblock Springer-Verlag, Cham, Switzerland (2017)

\bibitem{HackHopf13}
Hackl, C.M., Hopfe, N., Ilchmann, A., Mueller, M., Trenn, S.: Funnel control
  for systems with relative degree two.
\newblock {SIAM} J. Control Optim. \textbf{51}(2), 965--995 (2013)

\bibitem{HuIsi18}
Huang, J., Isidori, A., Marconi, L., Mischiati, M., Sontag, E., Wonham, W.:
  Internal models in control, biology and neuroscience.
\newblock In: Proc. 57th~{IEEE} Conf. Decis. Control, Miami Beach, FL, USA, pp.
  5370--5390 (2018)

\bibitem{IlchRyan06a}
Ilchmann, A., Ryan, E.P.: Asymptotic tracking with prescribed transient
  behaviour for linear systems.
\newblock Int. J. Control \textbf{79}(8), 910--917 (2006)

\bibitem{IlchRyan09}
Ilchmann, A., Ryan, E.P.: Performance funnels and tracking control.
\newblock Int. J. Control \textbf{82}(10), 1828--1840 (2009)

\bibitem{IlchRyan02a}
Ilchmann, A., Ryan, E.P., Sangwin, C.J.: Systems of controlled functional
  differential equations and adaptive tracking.
\newblock {SIAM} J. Control Optim. \textbf{40}(6), 1746--1764 (2002).

\bibitem{IlchRyan02b}
Ilchmann, A., Ryan, E.P., Sangwin, C.J.: Tracking with prescribed transient
  behaviour.
\newblock ESAIM: Control, Optimisation and Calculus of Variations \textbf{7},
  471--493 (2002)

\bibitem{IlchRyan06b}
Ilchmann, A., Ryan, E.P., Townsend, P.: Tracking control with prescribed
  transient behaviour for systems of known relative degree.
\newblock Syst. Control Lett. \textbf{55}(5), 396--406 (2006)

\bibitem{IlchRyan07}
Ilchmann, A., Ryan, E.P., Townsend, P.: Tracking with prescribed transient
  behavior for nonlinear systems of known relative degree.
\newblock {SIAM} J. Control Optim. \textbf{46}(1), 210--230 (2007).

\bibitem{IlchTren04}
Ilchmann, A., Trenn, S.: Input constrained funnel control with applications to
  chemical reactor models.
\newblock Syst. Control Lett. \textbf{53}(5), 361--375 (2004)

\bibitem{IlchWirt13}
Ilchmann, A., Wirth, F.: On minimum phase.
\newblock Automatisierungstechnik \textbf{12}, 805--817 (2013)

\bibitem{Isid95}
Isidori, A.: Nonlinear Control Systems, 3rd edn.
\newblock Communications and Control Engineering Series. Springer-Verlag,
  Berlin (1995)

\bibitem{JianMare04}
Jiang, Z., Mareels, I.M.Y., Hill, D.J., Huang, J.: A unifying framework for
  global regulation via nonlinear output feedback: {F}rom {ISS} to i{ISS}.
\newblock {IEEE} Trans. Autom. Control \textbf{49}, 549--562 (2004)

\bibitem{LeeTren19}
Lee, J.G., Trenn, S.: Asymptotic tracking via funnel control.
\newblock In: Proc. 58th~{IEEE} Conf. Decis. Control, Nice, France (2019).
\newblock To appear

\bibitem{LibeTren13b}
Liberzon, D., Trenn, S.: The bang-bang funnel controller for uncertain
  nonlinear systems with arbitrary relative degree.
\newblock {IEEE} Trans. Autom. Control \textbf{58}(12), 3126--3141 (2013)

\bibitem{LogeMawb00}
Logemann, H., Mawby, A.: Low-gain integral control of infinite dimensional
  regular linear systems subject to input hysteresis.
\newblock In: F.~Colonius, U.~Helmke, D.~Pr{\"a}tzel-Wolters, F.R. Wirth (eds.)
  Advances in Mathematical Systems Theory, pp. 255--293. Birkh{\"a}user,
  Boston, Basel, Berlin (2000)

\bibitem{Mare84}
Mareels, I.M.Y.: A simple selftuning controller for stably invertible systems.
\newblock Syst. Control Lett. \textbf{4}(1), 5--16 (1984)

\bibitem{MillDavi91}
Miller, D.E., Davison, E.J.: An adaptive controller which provides an
  arbitrarily good transient and steady-state response.
\newblock {IEEE} Trans. Autom. Control \textbf{36}(1), 68--81 (1991)

\bibitem{Mors83}
Morse, A.S.: Recent problems in parameter adaptive control.
\newblock In: I.D. Landau (ed.) Outils et Mod\`{e}les Math\'{e}matiques pour
  l'Automatique, l'Analyse de Syst\`{e}mes et le Traitment du Signal, vol.~3,
  pp. 733--740. \'Editions du Centre National de la Recherche Scientifique
  (CNRS), Paris (1983)

\bibitem{Na13}
Na, J.: Adaptive prescribed performance control of nonlinear systems with
  unknown dead zone.
\newblock Int. J. Adapt. Contr. and Sign. Proc. \textbf{27}, 426--446 (2013)

\bibitem{Nuss83}
Nussbaum, R.D.: Some remarks on a conjecture in adaptive control.
\newblock Syst. Control Lett. \textbf{3}, 243--246 (1983)

\bibitem{PompAlfo14}
Pomprapa, A., Alfocea, S.R., G{\"o}bel, C., Misgeld, B.J., Leonhardt, S.:
  Funnel control for oxygenation during artificial ventilation therapy.
\newblock In: Proceedings of the 19th IFAC World Congress, pp. 6575--6580. Cape
  Town, South Africa (2014)

\bibitem{PompWeye15}
Pomprapa, A., Weyer, S., Leonhardt, S., Walter, M., Misgeld, B.: Periodic
  funnel-based control for peak inspiratory pressure.
\newblock In: Proc. 54th~{IEEE} Conf. Decis. Control, Osaka, Japan, pp.
  5617--5622 (2015)

\bibitem{PuchReis19pp}
Puche, M., Reis, T., Schwenninger, F.L.: Funnel control for boundary control
  systems.
\newblock Evol. Eq. Control Th., to appear (2020), doi: 10.3934/eect.2020079

\bibitem{ReisSeli15b}
Reis, T., Selig, T.: Funnel control for the boundary controlled heat equation.
\newblock {SIAM} J. Control Optim. \textbf{53}(1), 547--574 (2015)

\bibitem{RyanSang01}
Ryan, E.P., Sangwin, C.J.: Controlled functional differential equations and
  adaptive stabilization.
\newblock Int. J. Control \textbf{74}(1), 77--90 (2001)

\bibitem{RyanSang09}
Ryan, E.P., Sangwin, C.J., Townsend, P.: Controlled functional differential
  equations: approximate and exact asymptotic tracking with prescribed
  transient performance.
\newblock ESAIM: Control, Optimisation and Calculus of Variations \textbf{15},
  745--762 (2009)

\bibitem{SeifBlaj13}
Seifried, R., Blajer, W.: Analysis of servo-constraint problems for
  underactuated multibody systems.
\newblock Mech. Sci. \textbf{4}, 113--129 (2013)

\bibitem{SenfPaug14}
Senfelds, A., Paugurs, A.: Electrical drive {DC} link power flow control with
  adaptive approach.
\newblock In: Proc. 55th Int. Sci. Conf. Power Electr. Engg. Riga Techn. Univ.,
  Riga, Latvia, pp. 30--33 (2014)


\bibitem{SontWang95b}
Sontag, E.D., Wang, Y.: On characterizations of the input-to-state stability
  property.
\newblock Syst. Control Lett. \textbf{24}(5), 351--359 (1995)

\bibitem{SponHutc06}
Spong, M.W., Hutchinson, S., Vidyasagar, M.: Robot Modeling and Control.
\newblock John Wiley and Sons Inc. (2006)

\bibitem{TaoKoko96}
Tao, G., Kokotovic, P.V.: Adaptive Control of Systems with Actuator and Sensor
  Nonlinearities.
\newblock Wiley, New York (1996)

\bibitem{TaoLewi01}
Tao, G., Lewis, F.L.: Adaptive Control of Nonsmooth Dynamic Systems.
\newblock Springer-Verlag, London (2001)

\bibitem{TrenStoo01}
Trentelman, H.L., Stoorvogel, A.A., Hautus, M.L.J.: Control Theory for Linear
  Systems.
\newblock Communications and Control Engineering. Springer-Verlag, London
  (2001).

\bibitem{Twain}
Twain, M.: Life on the Mississippi.
\newblock Osgood, Boston (1883)

\bibitem{VergDima19}
Verginis, C.K., Dimarogonas, D.D.: Asymptotic stability of uncertain
  {L}agrangian systems with prescribed transient response.
\newblock In: Proc. 58th~{IEEE} Conf. Decis. Control, Nice, France, pp. 7037--7042 (2019)

\bibitem{VergDima20}
Verginis, C.K., Dimarogonas, D.D.: Asymptotic tracking of second-order
  nonsmooth feedback stabilizable unknown systems with prescribed transient
  response.
\newblock {IEEE} Trans. Autom. Control (2021), doi: 10.1109/TAC.2020.3015785.
\newblock To appear

\bibitem{Wonh79}
Wonham, W.M.: Linear Multivariable Control: A Geometric Approach, 2nd edn.
\newblock Springer-Verlag, Heidelberg (1979)

\bibitem{Ye01}
Ye, X.: Adaptive nonlinear output-feedback control with unknown high-frequency
  gain sign.
\newblock {IEEE} Trans. Autom. Control \textbf{46}(1), 112--115 (2001)

\end{thebibliography}

\end{document}